\documentclass[preprint]{imsart}

\RequirePackage[OT1]{fontenc}
\RequirePackage{amsthm,amsmath,amssymb}
\RequirePackage[numbers]{natbib}
\RequirePackage[colorlinks,citecolor=blue,urlcolor=blue]{hyperref}

\usepackage{enumitem}
\usepackage{graphicx,psfrag,epsf}
\usepackage{epstopdf}
\usepackage{epsfig}
\usepackage[titletoc]{appendix} 
\usepackage{bbm}	
\usepackage{color}
\usepackage{accents} 
\usepackage{mathrsfs} 
\usepackage{url}
\usepackage{mathtools} 

\theoremstyle{plain}
\newtheorem{lemma}{Lemma}
\newtheorem{theorem}{Theorem}

\newtheorem{corollary}{Corollary}

\theoremstyle{definition} 
 
\newtheorem{example}{Example}
\newtheorem{remark}{Remark}

\newcommand\ind{\overset{\scriptscriptstyle\rm ind}{\sim}}
\newcommand\iid{\overset{\scriptscriptstyle\rm iid}{\sim}}
\newcommand{\N}{\mathrm{N}}
\newcommand{\Cov}{\mathrm{Cov}}

\makeatletter
\newcommand{\mylabel}[2]{#2\def\@currentlabel{#2}\label{#1}}
\makeatother

\begin{document}

\begin{frontmatter}
\title{Unified Bayesian theory of sparse linear regression with nuisance parameters\thanks{Research is partially supported by a Faculty Research and Professional Development Grant from College of Sciences of North Carolina State University.}
}
\runtitle{Bayesian sparse linear regression}

\begin{aug}
\author{\fnms{Seonghyun} \snm{Jeong}\ead[label=e1]{sjeong@yonsei.ac.kr}}

\address{Department of Statistics and Data Science, Department of Applied Statistics\\ Yonsei University, Seoul 03722, South Korea\\
	\printead{e1}
}

\author{\fnms{Subhashis} \snm{Ghosal}\ead[label=e2]{sghosal@ncsu.edu}}

\runauthor{S. Jeong and S. Ghosal}

\affiliation{North Carolina State University}

\address{Department of Statistics\\
North Carolina State University, Raleigh, NC 27607, USA\\
\printead{e2}
}

\end{aug}


\begin{abstract}
We study frequentist asymptotic properties of Bayesian procedures for high-dimensional Gaussian sparse regression when unknown nuisance parameters are involved. Nuisance parameters can be finite-, high-, or infinite-dimensional.
A mixture of point masses at zero and continuous distributions is used for the prior distribution on sparse regression coefficients, and appropriate prior distributions are used for nuisance parameters.
The optimal posterior contraction of sparse regression coefficients, hampered by the presence of nuisance parameters, is also examined and discussed.
It is shown that the procedure yields strong model selection consistency. A Bernstein-von Mises-type theorem for sparse regression coefficients is also obtained for uncertainty quantification through credible sets with guaranteed frequentist coverage. Asymptotic properties of numerous examples are investigated using the theories developed in this study.
\end{abstract}



\tableofcontents
\end{frontmatter}

\section{Introduction}
\label{c2sec:intro}

While Bayesian model selection for classical low-dimensional problems has a long history, sparse estimation in high-dimensional regression was studied much later; see \citet{bondell2012consistent}, \citet{johnson2012bayesian}, and \citet{narisetty2014bayesian} for consistent Bayesian model selection methods in high-dimensional linear models. 
Extensive theoretical investigations, however, have been carried out only very recently. Since the pioneering work of \citet{castillo2015bayesian}, frequentist asymptotic properties of Bayesian sparse regression have been discovered under various settings, and there is now a substantial body of literature  \citep[e.g.,][]{martin2017empirical,atchade2017contraction,song2017nearly,belitser2017empirical,rockova2018bayesian,bai2019spike,chae2016bayesian,ning2018bayesian,gao2015general,jeong2020posterior,jeong2020posterior-b}.

Most of the existing studies deal with sparse regression setups without nuisance parameters and there are only a few exceptions. An unknown variance parameter, the simplest type of nuisance parameters, was incorporated for high-dimensional linear regression in \citet{song2017nearly} and \citet{bai2019spike}. In these studies, the optimal properties of Bayesian procedures are characterized with continuous shrinkage priors. For more involved models, \citet{chae2016bayesian} adopted a nonparametric approach to estimate unknown symmetric densities in sparse linear regression. \citet{ning2018bayesian} considered a sparse linear model for vector-valued response variables with unknown covariance matrices.

Although nuisance parameters may not be of primary interest, modeling frameworks require the complete description of their roles as they explicitly parameterize models. Therefore, one may want to achieve optimal estimation properties for sparse regression coefficients, no matter what a nuisance parameter is. It may also be of interest to examine posterior contraction of nuisance parameters as a secondary objective.
Despite these facts, however, there have not been attempts to consider a general class of high-dimensional regression models with nuisance parameters. In this study, we consider a general form of Gaussian sparse regression in the presence of nuisance parameters, and establish a theoretical framework for Bayesian procedures.

We formulate a general framework to treat sparse regression models in a unified way as follows. Let $\eta$ be possibly an infinite-dimensional nuisance parameter taking values in a set $\mathbb{H}$. For each $\eta\in \mathbb{H}$ and an integer $m_i\in\{1,\dots,\overline m\}$ for some $\overline m\ge1$, suppose that there are a vector $\xi_{\eta,i}\in\mathbb{R}^{m_i}$ and a positive definite matrix $\Delta_{\eta,i}\in\mathbb{R}^{m_i\times m_i}$ which define a regression model for a vector-valued response variable $Y_i\in\mathbb{R}^{m_i}$ against covariates $X_i\in\mathbb{R}^{m_i\times p}$ given by
\begin{align}
Y_i&=X_i\theta+\xi_{\eta,i}+\varepsilon_i,\quad\varepsilon_i\ind \N_{m_i}(0,\Delta_{\eta,i}),
\quad i=1,\dots,n,
\label{c2eqn:model}
\end{align}
where $\theta\in\mathbb{R}^p$ is a vector of regression coefficients. Here $m_i$ (and $\overline m$) can increase with $n$. We consider the high-dimensional situation where $p>n$, but $\theta$ is assumed to be sparse, with many coordinates zero. The form in \eqref{c2eqn:model} clearly includes sparse linear regression with unknown error variances. Our main interest lies in more complicated setups. As will be shortly discussed in Section~\ref{sec:example}, many interesting examples belong to form \eqref{c2eqn:model}.

In this paper, we develop a unified theory of posterior asymptotics in the high-dimensional sparse regression models described by form \eqref{c2eqn:model}. To the best of our knowledge, there is no study thus far considering a general modeling framework of sparse regression as in \eqref{c2eqn:model}, even from the frequentist perspective. The results on complicated high-dimensional regression models are only available at model-specific levels and cannot be universally used for different model classes.
On the other hand, our approach is a unified theoretical treatment of the general model structure in \eqref{c2eqn:model} under the Bayesian framework. We establish general theorems on nearly optimal posterior contraction rates, a Bernstein-von Mises theorem via shape approximation to the posterior distribution of $\theta$, and model selection consistency. 

The general theory of posterior contraction using the canonical root-average-squared Hellinger metric on the joint density \citep{ghosal2007convergence} is not very useful in this context, since to recover rates in terms of the metric of interest on the regression coefficients, some boundedness conditions are needed \citep{jeong2020posterior}. To deal with this issue, we construct an exponentially powerful likelihood ratio test in small pieces that are sufficiently separated from the true parameters in terms of the average R\'enyi divergence of order $1/2$ (which coincides with the average negative log-affinity). This test provides posterior contraction relative to the corresponding divergence. The posterior contraction rates of $\theta$ and $\eta$ can then be recovered in terms of the metrics of interest under mild conditions on the parameter space. 
Due to a nuisance parameter $\eta$, the resulting posterior contraction for $\theta$ may be suboptimal.
Conditions for the optimal posterior contraction will also be examined.
Our results show that the obtained posterior contraction rates are adaptive to the unknown sparsity level. 

For a Bernstein-von Mises theorem and selection consistency, stronger conditions are required than those used for posterior contraction, in line with the existing literature 
\citep[e.g.,][]{castillo2015bayesian,martin2017empirical}.
As pointed out by \citet{chae2016bayesian}, the Bernstein-von Mises theorems for finite dimensional parameters in classical semiparametric models \citep[e.g.,][]{castillo2012semiparametric} may not be directly useful in the high-dimensional context. We thus directly characterize a version of the Bernstein von-Mises theorem for model \eqref{c2eqn:model}.
The key idea is to find a suitable orthogonal projection that satisfies some required conditions, which is typically straightforward if the support of a prior for $\xi_{\eta,i}$ is a linear space. The complexity of the space of covariance matrices, measured by its metric entropy, also has an important role in deriving the Bernstein-von Mises theorem and selection consistency. 
Combining these two leads to a single component of normal distributions for an approximation, which enables to correctly quantify remaining uncertainty on the parameter through the posterior distribution.

\subsection{Sparse linear regression with nuisance parameters}
\label{sec:example}

As briefly discussed above, the form in \eqref{c2eqn:model} is general and includes many interesting statistical models. Here we provide specific examples belonging to \eqref{c2eqn:model} in detail.
In Section~\ref{c2sec:app}, these examples will be used to apply the main results developed in this study.


\begin{example}[Multiple response models with missing components] We consider a general multiple response model with missing values, which is very common in practice. Suppose that for each $i$, a vector of $\overline m$ responses  with covariance matrix $\Sigma$ are supposed to be observed, but for the $i$th group (or subject) only $m_i$ entries are actually observed with the rest missing. Letting $Y_i\in\mathbb{R}^{m_i}$ be the $i$th observation and $Y_i^{\rm aug}\in\mathbb{R}^{\overline m}$ be the augmented vector of $Y_i$ and missing entries, we can write $Y_i=E_i^T Y_i^{\rm aug}$ and $\Cov(Y_i)=E_i^T\Sigma E_i$, where $E_i\in\mathbb{R}^{\overline m\times m_i}$ is the submatrix of the $\overline m\times \overline m$ identity matrix with the $j$th column included if the $j$th element of $Y_i^{\rm aug}$ is observed, $j=1,\ldots,{\overline m}$. Assuming that the mean of $Y_i$ is only $X_i\theta$ for covariates $X_i\in\mathbb{R}^{m_i\times p}$ and sparse coefficients $\theta\in\mathbb{R}^p$ with $p>n$, the model of interest can be written as $Y_i=X_i\theta+\varepsilon_i$, $\varepsilon_i\ind \N_{m_i}(0,E_i^T\Sigma E_i)$, $i=1,\dots,n$. The model belongs to the class described by \eqref{c2eqn:model} with $\xi_{\eta,i}=0_{m_i}$ and $\Delta_{\eta,i}=E_i^T\Sigma E_i$ for $\eta=\Sigma$.
	\label{c2exm:missing}
\end{example}

\begin{example}[Multivariate measurement error models] Suppose that a scalar response variable $Y_i^\ast\in\mathbb{R}$ is connected to fixed covariates  $X_i^\ast\in\mathbb{R}^{p}$ with $p>n$ and random covariates $Z_i\in\mathbb{R}^q$ with fixed $q\ge1$, through the following linear additive relationship:
	$Y_i^\ast=\alpha+X_i^{\ast T}\theta+Z_i^T\beta+\varepsilon_i^\ast$, $Z_i\iid \N_q(\mu,\Sigma)$, $\varepsilon_i^\ast\iid \N(0,\sigma^2)$, $i=1,\dots,n$.
	While $X_i^\ast$ is fully observed without noise, we observe a surrogate $W_i$ of $Z_i$ as $W_i=Z_i+\tau_i$, $\tau_i\iid \N_q(0,\Psi)$, where to ensure identifiability, $\Psi$ is assumed to be known. 
	This type of model is called a measurement error model or an errors-in-variables model; see \citet{fuller1987measurement} and \citet{carroll2006measurement} for a complete overview.
	By direct calculations, the joint distribution of $(Y_i^\ast,W_i)$ is given by
	\begin{align*}
	\begin{pmatrix}
	Y_i^\ast \\ W_i
	\end{pmatrix}
	\ind \N_{q+1} \left(
	\begin{pmatrix}
	\alpha+X_i^{\ast T}\theta+\mu^T\beta \\ \mu
	\end{pmatrix}
	,
	\begin{pmatrix}
	\beta^T\Sigma\beta+\sigma^2 & \beta^T\Sigma \\ \Sigma\beta & \Sigma+\Psi
	\end{pmatrix}\right).
	\end{align*}
	By writing $Y_i=(Y_i^\ast,W_i^T)^T\in\mathbb R^{q+1}$, $X_i=(X_i^{\ast},0_{p\times q})^T\in\mathbb R^{(q+1)\times p}$, $\xi_{\eta,i}=(\alpha+\mu^T\beta,\mu^T)^T\in\mathbb R^{q+1}$, and $\Delta_{\eta,i}=\left(\begin{smallmatrix}
	\beta^T\Sigma\beta+\sigma^2 & \beta^T\Sigma \\ \Sigma\beta & \Sigma+\Psi
	\end{smallmatrix}\right)\in\mathbb R^{(q+1)\times(q+1)}$ with $\eta=(\alpha,\beta,\mu,\sigma^2,\Sigma)$, the model is of form \eqref{c2eqn:model} with $m_i=q+1$.
	\label{c2exm:mem}
\end{example}

\begin{example}[Parametric correlation structure]
	For $m_i\ge1$, $i=1,\dots,n$, suppose that we have a response variable $Y_i\in \mathbb{R}^{m_i}$ and covariates $X_i\in\mathbb R^{m_i\times p}$ with $p>n$.
	We consider a standard regression model given by 
	$Y_i=X_i\theta+\varepsilon_i$, $\varepsilon_i \ind \N_{m_i}(0,\Sigma_i)$, $i=1,\dots,n$,
	but $m_i$ is considered to be possibly increasing.
	For a known parametric correlation structure $G_i$ and a fixed dimensional Euclidean parameter $\alpha$, we model the covariance matrix as $\Sigma_i=\sigma^2 G_i(\alpha)$ using a variance parameter $\sigma^2$ and a correlation matrix $G_i(\alpha)\in\mathbb{R}^{m_i\times m_i}$. 
	Examples of $G_i$ include first order autoregressive and moving average correlation matrices.
	The model belongs to \eqref{c2eqn:model} by writing $\xi_{\eta,i}=0_{m_i}$ and $\Delta_{\eta,i}=\sigma^2 G_i(\alpha)$ with $\eta=(\alpha,\sigma^2)$.
	\label{c2exm:parcor}
\end{example}

\begin{example}[Mixed effects models]
	For $m_i\ge1$, $i=1,\dots,n$, consider a response variable $Y_i
	\in\mathbb{R}^{m_i}$ and covariates $X_i\in\mathbb{R}^{m_i\times p}$ with $p>n$ and $Z_i\in\mathbb{R}^{m_i\times  q}$ with fixed $q\ge 1$.
	A mixed effect model given by $Y_i=X_i\theta+Z_ib_i+\varepsilon_i^\ast$, $b_i\iid {\rm N}_q(0,\Psi)$, $\varepsilon_i^\ast\ind {\rm N}_{m_i}(0,\sigma^2 I_{m_i})$, $i=1,\dots,n$, where $\Psi\in\mathbb{R}^{q\times q}$ is a positive definite matrix. Then the marginal law of $Y_i$ is given by $Y_i=X_i\theta+\varepsilon_i$, $\varepsilon_i\ind {\rm N}_{m_i}(0,\sigma^2I_{m_i}+Z_i\Psi Z_i^T)$. We assume that $\sigma^2$ is known. The model belongs to \eqref{c2eqn:model} by letting $\xi_{\eta,i}=0_{m_i}$ and $\Delta_{\eta,i}=\sigma^2I_{m_i}+Z_i\Psi Z_i^T$ with $\eta=\Psi$.
	\label{c2exm:mm}
\end{example}

\begin{example}[Graphical structure with sparse precision matrices]
	For a response variable $Y_i\in \mathbb{R}^{\overline m}$ and covariates $X_i\in \mathbb{R}^{\overline m\times p}$ with increasing $\overline m\ge 1$ and $p>n$, consider a model given by $Y_i=X_i\theta+\varepsilon_i$, $\varepsilon_i\iid \N_{\overline m}(0,\Omega^{-1})$, $i=1,\dots,n$,
	where $\theta$ is a sparse coefficient vector and the precision matrix $\Omega\in\mathbb{R}^{\overline m\times \overline m}$ is a positive definite matrix. Along with $\theta$, we also impose sparsity on the off-diagonal entries of $\Omega$, which accounts for a graphical structure between observations. More precisely, if an off-diagonal entry is zero, it implies the conditional independence between the two concerned entries of $\varepsilon_i$ given the remaining ones, and we suppose that most off-diagonal entries are actually zero, even though we do not know their locations. The model is then seen to be a special case of \eqref{c2eqn:model} by letting $\xi_{\eta,i}=0_{\overline m}$ and $\Delta_{\eta,i}=\Omega^{-1}$ with $\eta=\Omega$.
	\label{c2exm:graph}
\end{example}

\begin{example}[Nonparametric heteroskedastic regression models]
	For a response variable $Y_i\in\mathbb{R}$ and a row vector of covariates $X_i\in\mathbb{R}^{1\times p}$,	a linear regression model with a nonparametric heteroskedastic error is given by $Y_i=X_i\theta+\varepsilon_i$, $\varepsilon_i\ind \N(0,v(z_i))$, $i=1,\dots,n$,
	where $\theta$ is a sparse coefficient vector, $v:[0,1]\mapsto(0,\infty)$ is a univariate variance function, and $z_i\in[0,1]$ is a one-dimensional variable associated with the $i$th observation that controls the variance of $Y_i$ through the variance function $v$.
	Then the model belongs to \eqref{c2eqn:model} by letting $\xi_{\eta,i}=0$ and $\Delta_{\eta,i}=v(z_i)$ with $\eta=v$.
	\label{c2exm:het}
\end{example}

\begin{example}[Partial linear models] 
	Consider a partial linear model given by $	Y_i=X_i\theta+g(z_i)+\varepsilon_i$, $\varepsilon_i\iid \N(0,\sigma^2)$, $i=1,\dots,n$,
	where $Y_i\in\mathbb{R}$ is a response variable, $X_i\in\mathbb{R}^{1\times p}$ is a row vector of covariates with $p>n$, $\theta\in\mathbb{R}^p$ is a sparse coefficient vector, $g:[0,1]\mapsto\mathbb{R}$ is a univariate function, and $z_i\in[0,1]$ is a  scalar predictor. This model is expressed in form \eqref{c2eqn:model} by writing $\xi_{\eta,i}=g(z_i)$ and $\Delta_{\eta,i}=\sigma^2$ with $\eta=(g,\sigma^2)$.
	\label{c2exm:plm}
\end{example}

\subsection{Outline}

The rest of this paper is organized as follows. In Section~\ref{c2sec:setup}, some notations are introduced and a prior distribution on sparse regression coefficients is specified. Sections~\ref{c2sec:pcr}--\ref{c2sec:bvm} provide our main results on the posterior contraction, the Bernstein-von Mises phenomenon, and selection consistency of the posterior distribution. In Section~\ref{c2sec:app}, our general theorems are applied to the examples considered above to derive the posterior asymptotic properties in each specific example. All technical proofs are provided in Appendix.

\section{Setup, notations, and prior specification}
\label{c2sec:setup}

\subsection{Notation}
Here we describe the notations we use throughout this paper.
For a vector $\theta=(\theta_j)\in\mathbb{R}^p$ and a set $S\subset\{1,\dots,p\}$ of indices, we write $S_\theta=\{j:\theta_j\ne0\}$ to denote the support of $\theta$, $s \coloneqq |S|$ (or $s_\theta \coloneqq |S_\theta|$) to denote the cardinality of $S$ (or $S_\theta$), and $\theta_S=\{\theta_j: j\in S\}$ and $\theta_{S^c}=\{\theta_j: j\notin S\}$ to separate components of $\theta$ using $S$.
In particular, the support of the true parameter $\theta_0$ and its cardinality are written as $S_0$ and $s_0 \coloneqq |S_0|$, respectively. The notation $\lVert\theta\rVert_q=(\sum_j|\theta_j|^q)^{1/q}$, $1\le q < \infty$, stands for the $\ell_q$-norm and $\lVert\theta\rVert_\infty=\max_j|\theta_j|$ denotes the maximum norm. 
We write $\rho_{\min}(A)$ and $\rho_{\max}(A)$ for the minimum and maximum eigenvalues of a square matrix $A$, respectively. 
For a matrix $X=(\!(x_{ij})\!)$, let $\lVert X\rVert_{\rm sp}=\rho_{\max}^{1/2}(X^T X)$ stand for the spectral norm and $\lVert X \rVert_{\rm F}=(\sum_{i,j}x_{ij}^2)^{1/2}$ stand for the Frobenius norm of $X$.  We also define a matrix norm $\lVert X\rVert_\ast=\max_j \lVert X_{\cdot j}\rVert_2$ for $X_{\cdot j}$ the $j$th column of $X$, which is used for compatibility conditions. 
The column space of $X$ is denoted by ${\rm span}(X)$.
For further convenience, we write $\varsigma_{\min}(X)=\rho_{\min}^{1/2}(X^T X)$ for the minimum singular value of $X$.
The notation $X_S$ means the submatrix of $X$ with columns chosen by $S$.
For sequences $a_n$ and $b_n$, 
$a_n\lesssim b_n$ (or $b_n\gtrsim a_n$) stands for $a_n\le C b_n$ for some constant $C>0$ independent of $n$, and $a_n\asymp b_n$ means $a_n\lesssim b_n\lesssim a_n$. These inequalities are also used for relations involving constant sequences.

For given parameters $\theta$ and $\eta$, we write the joint density as $p_{\theta,\eta}=\prod_{i=1}^n p_{\theta,\eta,i}$ for $p_{\theta,\eta,i}$ the density of the $i$th observation vector $Y_i$.
In particular, the true joint density is expressed as $p_0=\prod_{i=1}^n p_{0,i}$ for 
$p_{0,i} \coloneqq p_{\theta_0,\eta_0,i}$ with the true parameters $\theta_0$ and $\eta_0$. 
The notation $\mathbb{E}_0$ denotes the expectation operator with the true density $p_{0}$.
For two probability measures $P$ and $Q$, let $\lVert P-Q\rVert_{\rm TV}$ denote the total variation between $P$ and $Q$.
For two $n$-variate densities $f \coloneqq \prod_{i=1}^n f_i$ and $g \coloneqq \prod_{i=1}^n g_i$ of independent variables,
denote the average R\'enyi divergence (of order $1/2$) by $R_n(f,g)=-n^{-1}\sum_{i=1}^n \log \int \sqrt{f_i g_i}$.

For any $\eta_1,\eta_2\in\mathbb{H}$, we define  $d_n^2(\eta_1,\eta_2)=d_{A,n}^2(\eta_1,\eta_2)+d_{B,n}^2(\eta_1,\eta_2)$
for the two squared pseudo-metrics:
\begin{align*}
d_{A,n}^2(\eta_1,\eta_2)=\frac{1}{n}\sum_{i=1}^n \lVert\xi_{\eta_1,i}-\xi_{\eta_2,i}\rVert_2^2,\quad d_{B,n}^2(\eta_1,\eta_2)=\frac{1}{n}\sum_{i=1}^n \lVert\Delta_{\eta_1,i}-\Delta_{\eta_2,i}\rVert_{\rm F}^2.
\end{align*} 
For compatibility conditions, the uniform compatibility number $\phi_1$ and the smallest scaled singular value $\phi_2$ are defined as
\begin{align*}
\phi_1(s)=\inf_{\theta:1\le |S_\theta|\le s}\frac{\lVert X\theta\rVert_2 |S_\theta|^{1/2}}{\lVert X\rVert_{\ast}\lVert \theta\rVert_1},\quad\phi_2(s)=\inf_{\theta:1\le |S_\theta|\le s}\frac{\lVert X\theta\rVert_2 }{\lVert X\rVert_{\ast}\lVert \theta\rVert_2}.
\end{align*}
We write $Y^{(n)}=(Y_1^T,\dots,Y_n^T)^T$ for the observation vector, $n_\ast=\sum_{i=1}^n m_i$ for the dimension of $Y^{(n)}$, and $\Theta=\mathbb{R}^p$ for the parameter space of $\theta$.
Lastly, for a (pseudo-)metric space $({\cal F}, d)$, let $N(\epsilon,{\cal F},d)$ denote the $\epsilon$-covering number, the minimal number of $\epsilon$-balls that cover $\cal F$.

\subsection{Prior for the high-dimensional coefficients}
\label{c2sec:prior}
In this subsection, we specify a prior distribution for the high-dimensional regression coefficients $\theta$. A prior for $\eta$ should satisfy the conditions required for the main results, so its specific characterization is deferred to Section~\ref{c2sec:pcr}. On the other hand,  the prior for $\theta$ specified here is always good for our purposes and satisfies all requirements.

We first select a dimension $s$ from a prior $\pi_p$, and then randomly choose $S\subset\{1,\dots,p\}$ for given $s$.
A nonzero part $\theta_S$ of $\theta$ is then selected from a prior $g_S$ on $\mathbb{R}^s$ while $\theta_{S^c}$ is fixed to zero. 
The resulting prior specification for $(S,\theta)$ is formulated as
\begin{align}
(S,\theta)\mapsto\frac{\pi_p(s)}{\binom{p}{s}}g_S(\theta_S)\delta_0(\theta_{S^c}),
\label{c2eqn:prior}
\end{align}
where $\delta_0$ is the Dirac measure at zero on $\mathbb{R}^{p-s}$ with suppressed dimensionality. For the prior $\pi_p$ on the model dimensions, we consider a prior satisfying the following: for some constants $A_1,A_2,A_3,A_4>0$, 
\begin{align}
A_1p^{-A_3}\pi_p(s-1)\le\pi_p(s)\le A_2p^{-A_4}\pi_p(s-1),\quad s=1,\dots,p.
\label{c2eqn:prdim}
\end{align}
Examples of priors satisfying \eqref{c2eqn:prdim} can be found in \citet{castillo2012needles} and \citet{castillo2015bayesian}. 
For the prior $g_S$, the $s$-fold product of the exponential power density is considered, where the regularization parameter is allowed to vary with $p$ and $\lVert X \rVert_{\ast}$, i.e.,
\begin{align}
g_S(\theta_S)=\prod_{j\in S}\frac{\lambda}{2}\exp\left(-\lambda |\theta_j|\right),\quad
\frac{\lVert X \rVert_{\ast}}{L_1 p^{L_2}}\le\lambda\le \frac{L_3\lVert X\rVert_\ast}{\sqrt{n}},
\label{c2eqn:lambda}
\end{align}
for some constants $L_1,L_2,L_3>0$. The order of $\lambda$ is important in that it determines the boundedness requirement of the true signal $\theta_0$ (see condition \ref{c2con:trueparcon1} below). A particularly interesting case is obtained when $\lambda$ is set to the lower bound $\lVert X \rVert_{\ast}/(L_1p^{L_2})$. Then the boundedness condition becomes very mild by choosing $L_2$ sufficiently large. When $\lambda$ is set to the upper bound, the boundedness condition is still reasonably mild.  However, it can actually be relaxed if the true signal is known to be small enough, though we do not pursue this generalization in this study. In Section~\ref{c2sec:bvm}, we shall see that values of $\lambda$ that do not increase too fast are in fact necessary for a distributional approximation and selection consistency.

\begin{remark}
Since some size restriction on $\theta_0$ will be made unlike \citet{castillo2015bayesian},
we note that the use of the Laplace density is not necessary and other prior distributions may also be used for $\theta$. For example, normal densities can be used for $g_S$ to exploit semi-conjugacy. However, if its precision parameter is fixed independent of $n$, a normal prior requires a stronger restriction on the true signal than \ref{c2con:trueparcon1} below.
To achieve the nearly optimal posterior contraction, other densities with similar tail properties should also work with appropriate modifications for the true signal size (see, e.g., \citet{jeong2020posterior}).
Instead of the spike-and-slab prior in \eqref{c2eqn:prior} and \eqref{c2eqn:prdim}, a class of continuous shrinkage priors may also be used at the expense of substantial modifications in the technical details \citep{song2017nearly}. In this paper, we only consider the prior in \eqref{c2eqn:prior}--\eqref{c2eqn:lambda}.

\end{remark}


\section{Posterior contraction rates}
\label{c2sec:pcr}
The prior for a nuisance parameter $\eta$ should be chosen to complete the prior specification.
Once we assign the prior for the full parameters, the posterior distribution $\Pi(\cdot\,|\,Y^{(n)})$ is defined by Bayes' rule. How the prior for $\eta$ is chosen is crucial to obtain desirable asymptotic properties of the posterior distribution. 
In this subsection, we shall examine such conditions on the prior distribution for a nuisance parameter and study the posterior contraction rates for both $\theta$ and $\eta$. 

The prior for $\eta$ is put on a subspace ${\cal H}\subset\mathbb{H}$. In many instances, we take ${\cal H}={\mathbb H}$, especially when a nuisance parameter is finite dimensional, but the flexibility of a subspace may be beneficial in infinite-dimensional situations.
We need to choose $\cal H$ to satisfy certain conditions.

\begin{enumerate}[label=(C\arabic*),leftmargin=3.0\parindent, itemsep=0.0cm]
	\item \label{c2con:maxfrob} There exists a nondecreasing sequence $a_n=o(n)$ such that
	\begin{align*}
	a_n \max_{1\le i\le n}\lVert\Delta_{\eta',i}-\Delta_{\eta_0,i}\rVert_{\rm F}^2& \eqqcolon e_n\rightarrow 0,\quad \text{for some $\eta'\in{\cal H}$},\\
	\max_{1\le i\le n}\lVert\Delta_{\eta_1,i}-\Delta_{\eta_2,i}\rVert_{\rm F}^2&\le a_n d_{B,n}^2(\eta_1,\eta_2),\quad \eta_1,\eta_2\in{\cal H}.
	\end{align*}
	\item \label{c2con:thm1prcon} For some sequence $\bar\epsilon_n$ such that $a_n\bar\epsilon_n^2\rightarrow0$ and $n\bar\epsilon_n^2\rightarrow\infty$ with $a_n$ satisfying \ref{c2con:maxfrob},
	\begin{align*}
	\log\Pi\left(\eta\in{\cal H}:d_n(\eta,\eta_0)\le \bar\epsilon_n\right) \gtrsim- n\bar\epsilon_n^2.
	\end{align*}
\end{enumerate}
The first condition of \ref{c2con:maxfrob} implies that we have a good approximation to the true parameter value in the parameter set $\cal H$. This holds trivially if there exists $\eta'\in{\cal H}$ such that $\Delta_{\eta',i}=\Delta_{\eta_0,i}$ for every $i\le n$,
which is obviously true if $\eta_0\in{\cal H}$. The second condition of \ref{c2con:maxfrob} means that in $\cal H$, the maximum Frobenius norm of the difference between covariance matrices can be controlled by the average Frobenius norm multiplied by the sequence $a_n$. Clearly, this holds with $a_n=1$ if $\Delta_{\eta,i}$ is the same for every $i\le n$. By the triangle inequality, we see that \ref{c2con:maxfrob} implies that 
\begin{align}
\max_{1\le i\le n}\lVert\Delta_{\eta,i}-\Delta_{\eta_0,i}\rVert_{\rm F}^2\lesssim e_n+a_n d_{B,n}^2(\eta,\eta_0),\quad \eta\in{\cal H},
\label{c2eqn:maxfrob2}
\end{align}
which is used throughout the paper. 
Condition \ref{c2con:thm1prcon} is typically called the prior concentration condition, which requires a prior to put sufficient mass around the true parameter $\eta_0$, measured by the pseudo-metric $d_n$. As in other infinite-dimensional situations, such a closeness is translated into the closeness in terms of the Kullback-Leibler divergence and variation (see Lemma~\ref{c2lmm:1} in Appendix for more details).

As noted in Section~\ref{c2sec:intro}, the true parameters should be restricted to certain norm-bounded subset of the parameter space. This is clarified as follows.
\begin{enumerate}[label=(C\arabic*),leftmargin=3.0\parindent,resume, itemsep=0.0cm]
	\item \label{c2con:trueparcon1} The true signal satisfies $\lVert\theta_0\rVert_\infty\lesssim \lambda^{-1}\log p$.
	\item \label{c2con:trueparcon2} The eigenvalues of the true covariance matrix satisfy 
	\begin{align*}
	1\lesssim\min_{1\le i\le n}\rho_{\min}(\Delta_{\eta_0,i})\le \max_{1\le i\le n}\rho_{\max}(\Delta_{\eta_0,i})\lesssim 1.
	\end{align*}
\end{enumerate}
Condition \ref{c2con:trueparcon1} is required to apply the general strategy for posterior contraction to our modeling framework containing nuisance parameters.
More specifically, the condition is imposed such that the prior assigns sufficient mass on a Kullback-Leibler neighborhood of $\theta_0$.
If nuisance parameters are not present, one can directly handle the model and such a restriction may be removed \citep[e.g.,][]{castillo2015bayesian,gao2015general}.
One may refer to \citet{song2017nearly}, \citet{ning2018bayesian}, and \citet{bai2019spike} for conditions similar to ours, where a variance parameter stands for a nuisance parameter. 
Still, the condition is mild if $\lambda$ is chosen to decrease at an appropriate order. In particular, if $\lambda$ is matched to the lower bound $1/(L_1p^{L_2})$, the condition becomes $\lVert \theta_0\rVert_\infty\lesssim (p^{L_2}\log p)/\lVert X\rVert_\ast$ which is very mild if $L_2$ is sufficiently large.
Even if the upper bound $L_3\lVert X\rVert_\ast/\sqrt{n}$ is chosen, the condition is not  restrictive as the right hand side of the condition can be made nondecreasing as long as $\lVert X\rVert_\ast$ is increasing at a suitable order.
Condition \ref{c2con:trueparcon2} implies that the eigenvalues of the true covariance matrix are bounded below and above. The lower and upper bounds are required for a lot of technical details, including the construction of an exponentially powerful test in Lemma~\ref{c2lmm:test} in Appendix.

\begin{remark}
Condition \ref{c2con:trueparcon1} is actually stronger than what it needs to be, but is adopted for the ease of interpretation. For Theorem~\ref{c2thm:thetacon} below to hold, it suffices if we have $\lambda\lVert \theta_0\rVert_1\le (s_0\log p)\vee n\bar\epsilon_n^2$ for $\bar\epsilon_n$ satisfying \ref{c2con:thm1prcon}. For the optimal posterior contraction in Theorem~\ref{c2thm:thmmar} below, a slightly stronger bound is needed: $\lambda\lVert \theta_0\rVert_1\le s_0\log p$ (see Lemma~\ref{lmm:n2} and its proof in Appendix).
\end{remark}

\subsection{R\'enyi posterior contraction  and recovery}
\label{sec:renyi}

The goal of this subsection is to study posterior contraction of $\theta$ relative to the $\ell_1$- and $\ell_2$-metrics. To do so, we derive the posterior contraction rate with respect to the average R\'enyi divergence $R_n(f,g)$, and then the rates for $\theta$ relative to more concrete metrics will be recovered from the R\'enyi contraction.

To proceed, we first need to examine a dimensionality property of the support of $\theta$.
The following theorem shows that the posterior distribution is concentrated on models of relatively small sizes.
\begin{theorem}[Dimension]
	Suppose that {\rm\ref{c2con:maxfrob}--\ref{c2con:trueparcon2}} are satisfied.
	Then for $s_\star  \coloneqq s_0\vee (n\bar\epsilon_n^2/\log p)$, there exists a constant $K_1$ such that
	\begin{align*}
	\mathbb{E}_0\Pi\left(\theta:s_\theta >K_1 s_\star \,\big|\,Y^{(n)}\right)\rightarrow 0.
	\end{align*}
	\label{c2thm:dimen}
\end{theorem}
Compared to the literature
\citep[e.g.,][]{castillo2015bayesian,martin2017empirical,belitser2017empirical}, 
the rate in Theorem~\ref{c2thm:dimen} is floored by the extra term $n\bar\epsilon_n^2/\log p$. This arises from the presence of a nuisance parameter in the model formulation. 
To minimize its impact, a prior on $\eta$ should be chosen such that \ref{c2con:thm1prcon} holds for as small $\bar\epsilon_n$ as possible; a suitable choice induces the (nearly) optimal contraction rate.


Using the basic results in Theorem~\ref{c2thm:dimen}, the next theorem obtains the rate at which the posterior distribution contracts at the truth with respect to the average R\'enyi divergence. The theorem requires additional assumptions on a prior.
\begin{enumerate}[label=(C\arabic*),leftmargin=3.0\parindent,resume, itemsep=0.0cm]
	\item \label{c2con:thm2} For $s_\star  \coloneqq s_0\vee (n\bar\epsilon_n^2/\log p)$ with $\bar\epsilon_n$ satisfying \ref{c2con:thm1prcon}, a sufficiently large $B>0$, and some sequences $\gamma_n$ and $\epsilon_n\ge\sqrt{s_\star \log (p\vee\overline m\vee \gamma_n) /n}$ satisfying $\epsilon_n^2/\overline m\rightarrow 0$, there exists a subset ${\cal H}_{n}\subset{\cal H}$ such that 
	\begin{align}
	\min_{1\le i\le n}\inf_{\eta\in{\cal H}_{n}}\rho_{\min}(\Delta_{\eta,i})&\ge \frac{1}{\gamma_n},\label{c2eqn:thm2c1}\\
	\log N\left(\frac{1}{6\overline m\gamma_n n^{3/2}},{\cal H}_{n},d_n\right)&\lesssim n\epsilon_n^2,\label{c2eqn:thm2c2}\\
	e^{Bs_\star \log p}\Pi({\cal H}\setminus{\cal H}_{n})&\rightarrow 0.\label{c2eqn:thm2c3}
	\end{align}
\end{enumerate}
The above conditions are related to the classical ones in the literature (e.g., see Theorem 2.1. of \citet{ghosal2000convergence}).
Condition \eqref{c2eqn:thm2c1} requires that for every $i\le n$, the minimum eigenvalue of $\Delta_{\eta,i}$ is not too small on a sieve ${\cal H}_n$. Although $\gamma_n$ can be any positive sequence, a sequence increasing exponentially fast makes the entropy in \eqref{c2eqn:thm2c2} too large, resulting in a suboptimal rate $\epsilon_n$. If $\gamma_n$ can be chosen to be smaller than $p$ and $\overline m$, then this does not lead to any deterioration of the rate in $\epsilon_n$.
The entropy condition \eqref{c2eqn:thm2c2} is actually stronger than needed. Scrutinizing the proof of the theorem, one can see that the entropy appearing in the theorem is obtained using pieces that are smaller than those giving the exponentially powerful test in Lemma~\ref{c2lmm:test} in Appendix. However, the covering number with those pieces looks more complicated and the form in \eqref{c2eqn:thm2c2} suffices for all examples in the present paper. Lastly, condition \eqref{c2eqn:thm2c3} implies that the outside of a sieve ${\cal H}_n$ should possess sufficiently small prior mass to kill the factor $s_\star \log p$ arising from the lower bound of the denominator of the posterior distribution.
In fact, conditions similar to \ref{c2con:thm1prcon}, \eqref{c2eqn:thm2c2} and \eqref{c2eqn:thm2c3} are also required for the prior of $\theta$. By reading the proof, it is easy to see that the prior \eqref{c2eqn:prior} explicitly satisfies the analogous conditions on an appropriately chosen sieve.

\begin{theorem}[Contraction rate, R\'enyi]
	Suppose that {\rm\ref{c2con:maxfrob}--\ref{c2con:thm2}} are satisfied.
	Then there exists a constant $K_2$ such that
	\begin{align*}
	{\mathbb E}_0\Pi\left((\theta,\eta):R_n(p_{\theta,\eta},p_0)>K_2\epsilon_n^2 \,\big |\,Y^{(n)}\right)\rightarrow 0.
	\end{align*}
	\label{c2thm:helcon}
\end{theorem}
We want to sharpen the rate $\epsilon_n\ge\sqrt{s_\star \log (p\vee\overline m\vee \gamma_n) /n}$ as much as possible.
In most instances, $\gamma_n$ can be chosen such that $\log \gamma_n\lesssim \log p$. This is trivially satisfied if $\gamma_n$ is some polynomial in $n$ as in the examples in this paper. If $p$ is known to increase much faster than $n$, e.g., $\log p\asymp n^{c}$ for some $c\in(0,1)$, then $\gamma_n$ need not be a polynomial in $n$ and the condition can be met more easily with a sequence that grows even faster. 
Note also that we typically have $\log \overline m\lesssim \log p$ in most cases.
These postulates lead to $\epsilon_n\ge\sqrt{(s_\star \log p) /n}$.
Indeed, it is often possible to choose $\epsilon_n=\sqrt{(s_\star \log p)/n}$, which is commonly guaranteed by choosing an appropriate sieve ${\cal H}_n$ and a prior. 
The condition will be made precise in \ref{c2con:thm2mod1} below for recovery and we only consider the  situation that $\epsilon_n=\sqrt{(s_\star \log p)/n}$ in what follows.

Although Theorem~\ref{c2thm:helcon} provides the basic results for posterior contraction, it does not give precise interpretations for the parameters $\theta$ and $\eta$ themselves, because of the abstruse expression of the average R\'enyi divergence. The contraction rates with respect to more concrete metrics are recovered under some additional conditions. 
Under the additional assumption $a_n\epsilon_n^2\rightarrow0$, it can be shown that Theorem~\ref{c2thm:dimen} and Theorem~\ref{c2thm:helcon} explicitly imply that for the set
\begin{align*}
{\cal A}_n =\bigg\{ &(\theta,\eta)\in \Theta\times{\cal H} :  s_\theta\le K_1s_\star , \\
& \qquad\qquad\frac{1}{n}\sum_{i=1}^n \lVert  X_i(\theta-\theta_0)+\xi_{\eta,i}-\xi_{\eta_0,i} \rVert_2^2+ d_{B,n}^2(\eta,\eta_0)\le M_1 \epsilon_n^2\bigg\},
\end{align*}
with a sufficiently large constant $M_1$, the posterior mass of ${\cal A}_n$ goes to one in probability (see the proof of Theorem~\ref{c2thm:thetacon}).
To complete the recovery, we need to separate the sum of squares of the mean into $\lVert X(\theta-\theta_0) \rVert_2$ and $nd_{A,n}^2(\eta,\eta_0)$, which requires an additional condition. The conditions required for the  recovery are clarified as follows.

\begin{enumerate}[label=(C\arabic*),leftmargin=3.0\parindent,resume, itemsep=0.0cm]
	\item[\mylabel{c2con:thm2mod1}{(C5{$^\ast$})}] 
	While $\log \overline m\lesssim \log p$, \ref{c2con:thm2} holds for $\gamma_n$ and $\epsilon_n=\sqrt{(s_\star \log p)/n}$ such that $\log \gamma_n\lesssim\log p$ and $a_n\epsilon_n^2\rightarrow 0$ with $a_n$ satisfying \ref{c2con:maxfrob}.
	\item \label{c2con:sep}  For $s_\star $ satisfying \ref{c2con:thm2mod1}, there exists $\eta_\ast\in\mathbb{H}$ such that
	\begin{align*}
	\liminf_{n\ge1}\inf_{(\theta,\eta)\in {\cal A}_n}\frac{\sum_{i=1}^n (\theta-\theta_0)^T X_i^T(\xi_{\eta,i}-\xi_{\eta_\ast,i})}{\lVert X(\theta-\theta_0)\rVert_2^2 +n d_{A,n}^2(\eta,\eta_\ast)}&>-\frac{1}{2},\\
	d_{A,n}(\eta_\ast,\eta_0)&\lesssim \sqrt{\frac {s_\star \log p}{n}},
	\end{align*}
	where $\epsilon_n$ in $\mathcal A_n$ satisfies $\epsilon_n=\sqrt{(s_\star \log p)/n}$.
\end{enumerate}

By expanding the quadratic term for the mean in ${\cal A}_n$,
one can see that the separation is possible if \ref{c2con:sep} is satisfied.
Clearly, \ref{c2con:sep} is trivially satisfied if the model has only $X\theta$ for its mean, in which we take $\xi_{\eta,i}-\xi_{\eta_\ast,i}=\xi_{\eta_\ast,i}-\xi_{\eta_0,i}=0$ for every $i\le n$. 
In many cases where there exists $\eta'\in{\cal H}$ such that $d_{A,n}(\eta',\eta_0)=0$, we can often take $\eta_\ast = \eta'$ for the second inequality of \ref{c2con:sep} to hold automatically.

The following theorem shows that the posterior distribution of $\theta$ and $\eta$ contracts at their respective true values at some rates, relative to more easily comprehensible metrics than the average R\'enyi divergence. 
In the expressions, if $K_1s_\star +s_0<1$, the compatibility numbers should be understood be equal to 1 for interpretation.
\begin{theorem}[Recovery]
	Suppose that {\rm\ref{c2con:maxfrob}--\ref{c2con:trueparcon2}}, {\rm\ref{c2con:thm2mod1}}, and {\rm\ref{c2con:sep}} are satisfied.
	Then, there exists a constant $K_3$ such that
	\begin{align}
	\begin{split}
	{\mathbb E}_0\Pi\left(\theta:\lVert \theta-\theta_0\rVert_1>\frac{K_3s_\star \sqrt{\log p }}{\phi_1(K_1s_\star +s_0)\lVert X\rVert_{\ast}} \,\bigg|\,Y^{(n)} \right)&\rightarrow 0,\\
	{\mathbb E}_0\Pi\left(\theta:\lVert \theta-\theta_0\rVert_2> \frac{K_3\sqrt{s_\star \log p }}{\phi_2(K_1s_\star +s_0)\lVert X\rVert_{\ast}} \,\bigg|\,Y^{(n)} \right)&\rightarrow 0,\\
	{\mathbb E}_0\Pi\left(\theta:\lVert X(\theta-\theta_0)\rVert_2>K_3 \sqrt{s_\star \log p } \,\big |\,Y^{(n)} \right)&\rightarrow 0,\\
	{\mathbb E}_0\Pi\left(\eta:d_n(\eta,\eta_0)>K_3\sqrt{\frac{s_\star \log p}{n} } \,\bigg|\,Y^{(n)} \right)&\rightarrow 0.
	\end{split}
	\label{c2eqn:thm3}
	\end{align}
	\label{c2thm:thetacon}
\end{theorem}
The thresholds for contraction depend upon the compatibility conditions, which make their implication somewhat vague.
As $K_1s_\star +s_0$ is much smaller than  $n_\ast$, it is not unreasonable to assume that $\phi_1(K_1s_\star +s_0)$ and $\phi_2(K_1s_\star +s_0)$ are bounded away from zero, whence the compatibility number is removed from the rates. 
We refer to Example~7 of \citet{castillo2015bayesian} for more discussion.
In the next subsection, we will see that one of these restrictions is actually necessary for shape approximation or selection consistency.

\begin{remark}
	The separation condition \ref{c2con:sep} can be left as an assumption to be satisfied, but can also be verified by a stronger condition on the design matrix without resorting to the values of the parameters. Suppose that for some integer $q\ge1$, there exists a matrix $Z_i\in\mathbb{R}^{m_i\times q}$ such that $\xi_{\eta,i}=Z_i h(\eta)$ for every $\eta\in{\cal H}$, with some map $h:{\cal H}\mapsto\mathbb{R}^{q}$.
	Since we can write $\xi_{\eta,i}-\xi_{\eta_\ast,i}=Z_i(h(\eta)-h(\eta_\ast))$ for any $\eta,\eta_\ast\in{\cal H}$, the Cauchy-Schwarz inequality indicates that the first inequality of \ref{c2con:sep} is implied by
	\begin{align*}
	\liminf_{n\ge 1}\inf_{(\theta,\eta)\in\Theta\times{\cal H}:s_\theta\le K_1 s_\star }\frac{(\theta-\theta_0)^T X^T Z(h(\eta)-h(\eta_\ast))}{\lVert X(\theta-\theta_0)\rVert_2\lVert Z(h(\eta)-h(\eta_\ast))\rVert_2}>-1,
	\end{align*}
	for $Z=(Z_1^T,\dots,Z_n^T)^T$.
	The left hand side is always between $-1$ and $1$ by the Cauchy-Schwarz inequality, and is exactly equal to $-1$ or $1$ if and only if the two vectors are linearly dependent. 
	A sufficient condition for the preceding display is thus $\min\{\varsigma_{\min}([X_S,Z]):{s\le K_1 s_\star +s_0}\}\gtrsim 1$ since the linear dependence cannot happen under such a condition 
	due to the inequality $s_{\theta-\theta_0}\le s_\theta+s_0\le K_1s_\star +s_0$ for $\theta$ such that $s_\theta\le K_1 s_\star $.
	This sufficient condition is not restrictive at all if  $q=o(n)$ as we already have $K_1s_\star +s_0 = o(n)$.
	Since there typically exists $\eta_\ast\in{\cal H}$ satisfying the second inequality of \ref{c2con:sep} as long as $\cal H$ provides a good approximation for the true parameter $\eta_0$, condition \ref{c2con:sep} can be easily satisfied if the sufficient condition is met.
	\label{rmk:1}
\end{remark}

Notwithstanding the lack of formal study of minimax rates with additional complications, we still want to match our rates for $\theta$ with those in simple linear regression, which we call the ``optimal'' rates.
In this sense, Theorem~\ref{c2thm:thetacon} only provides the suboptimal rates for $\theta$ if $s_0=o(s_\star)$.
Although the theorem gives the optimal results if $s_0\log p\gtrsim n\bar\epsilon_n^2$, it is practically hard to check this condition as $s_0$ is unknown. If $s_0$ is known to be nonzero, the desired conclusion is trivially achieved as soon as $n\bar\epsilon_n^2/\log p\lesssim 1$. The following corollary, however, shows that the optimal rates are still available  even if $s_0=0$, with restrictions on $\bar\epsilon_n$ and the prior.

\begin{corollary}[Optimality under restriction]
	For $\bar\epsilon_n$ satisfying the conditions for Theorem~\ref{c2thm:thetacon}, we have the following assertions.
	\begin{enumerate}[label={\rm (\alph*)}]
		\item \label{c2cor:opt-a} Assume that $n\bar\epsilon_n^2/\log p\rightarrow 0$. Then, Theorems~\ref{c2thm:dimen} and \ref{c2thm:thetacon} hold for $s_\star$ replaced by $s_0$.
		\item \label{c2cor:opt-b} Assume that $n\bar\epsilon_n^2/\log p\lesssim 1$.  Then, Theorems~\ref{c2thm:dimen} and \ref{c2thm:thetacon} hold for $s_\star$ replaced by $s_0$ if either $A_4$ in \eqref{c2eqn:prdim} is chosen large enough or $s_0>0$.
	\end{enumerate} 
	\label{c2cor:opt}
\end{corollary}
The corollary is useful in limited situations, especially when a parametric rate is available for a nuisance parameter. 
Even if $n\bar\epsilon_n^2=\log n$, we need to further assume that $\log n=o(\log p)$, i.e., the ultra high-dimensional setup, to conclude that \ref{c2cor:opt-a} holds, while we can always apply \ref{c2cor:opt-b} because $\log n\lesssim \log p$. 
Although assertion \ref{c2cor:opt-b} holds for any $s_0\ge0$ if  $A_4$ is chosen sufficiently large, its specific threshold is not directly available. Indeed, by carefully reading the proof of Theorem~\ref{c2thm:dimen} together with Lemma~\ref{c2lmm:1} in Appendix, one can see that the threshold depends on unknown constant bounds for the eigenvalues of the true covariance matrix in \ref{c2con:trueparcon2}. Still, \ref{c2cor:opt-b} holds for any $A_4>0$ if $s_0>0$. We believe that the assumption $s_0>0$ is very mild, and hence simply apply \ref{c2cor:opt-b} with this assumption to conclude the optimal contraction for models with finite dimensional nuisance parameters.
The optimal rates can still be achieved for any $s_0\ge 0$ by verifying the conditions in the following subsection. With finite dimensional nuisance parameters, we do not pursue this direction as it seems an overkill considering the mildness of the assumption  $s_0>0$, though those conditions are actually required for the Bernstein-von Mises theorem and selection consistency in Section~\ref{c2sec:bvm}.

In semiparametric situations with high- or infinite-dimensional nuisance parameters, none of \ref{c2cor:opt-a} and \ref{c2cor:opt-b} generally works unless $p$ increases sufficiently fast. Still, the optimal rates can be achieved under stronger conditions using the semiparametric theory,
as the following subsection provides.

\subsection{Optimal posterior contraction for $\theta$}
\label{sec:marcon}
Recall that only suboptimal rates may be available from Theorem~\ref{c2thm:thetacon} if $s_0\log p\lesssim n\bar\epsilon_n^2$.
In many semiparametric situations, however, it is often possible to obtain parametric rates for finite dimensional parameters under stronger conditions, even when there are infinite-dimensional nuisance parameters in a model \citep{bickel2012semiparametric,castillo2012semiparametric}. It has also been shown that a similar argument holds in some high-dimensional semiparametric regression models \citep{chae2016bayesian}.
Therefore, it is naturally of interest to examine under what conditions we can replace $s_\star$ by $s_0$ in the rates for $\theta$, even if $s_0\log p\lesssim n\bar\epsilon_n^2$.
Similar to other semiparametric settings \citep{bickel2012semiparametric,chae2016bayesian}, this can be established by the semiparametric theory, but requires stronger conditions than those in traditional fixed dimensional parametric cases because of the high-dimensions of the parameters in our setup.

To proceed, some additional conditions are required for technical reasons, which are made  for the size of $\bar\epsilon_n$ as the optimal rates are automatically attained if $s_0\log p\gtrsim n\bar\epsilon_n^2$. Still, in a practical sense, the conditions almost always need to be verified to reach the optimal rates, since only oracle rates are generally available and we do not know which term is greater. 

In what follows, we write $\bar s_\star\coloneqq n\bar\epsilon_n^2/\log p$ for $\bar\epsilon_n$ satisfying the conditions of Theorem~\ref{c2thm:thetacon} through the definition of $\epsilon_n$.
We first assume the following condition on the uniform compatibility number.
\begin{enumerate}[label=(C\arabic*),leftmargin=3.0\parindent,resume, itemsep=0.0cm]
	\item \label{c2con:combo1} 
	For a sufficiently large $M$, the uniform compatibility number $\phi_1(M\bar s_\star +s_0)$ is bounded away from zero.
\end{enumerate}
This condition is weaker than assuming that the smallest scaled singular value $\phi_2(M\bar s_\star +s_0)$ is bounded away from zero, as we have $\phi_1(s)\ge \phi_2(s)$ for any $s>0$ by the Cauchy-Schwarz inequality.
We will also resort on a slightly stronger condition with respect to $\phi_1$ for a distributional approximation in the following section.
In this sense, our condition is weaker than those for Theorem~4 of \citet{castillo2015bayesian}.
Condition \ref{c2con:combo1} is not restrictive as \ref{c2con:thm2mod1} requires $s_\star =o(n)$; we again refer to Example~7 of \citet{castillo2015bayesian}.

To precisely describe other conditions, hereafter we use the following additional notations.
We write
\begin{align*}
\tilde X=\left(
\begin{matrix}
\Delta_{\eta_0,1}^{-1/2}X_1 \\ \vdots \\ \Delta_{\eta_0,n}^{-1/2}X_n
\end{matrix}\right)\in\mathbb R^{n_\ast\times p}
,\quad \tilde \xi_\eta=\left(
\begin{matrix}
\Delta_{\eta_0,1}^{-1/2}\xi_{\eta,1} \\ \vdots  \\ \Delta_{\eta_0,n}^{-1/2}\xi_{\eta,n}
\end{matrix}\right)\in\mathbb R^{n_\ast},
\end{align*}
and $\tilde\Delta_\eta$ to denote the collection of $\Delta_{\eta,i}$ for $i=1,\dots,n$.
In particular, $\tilde X_S\in\mathbb R^{n_\ast\times |S|}$ denotes the submatrix of $\tilde X$ with columns chosen by an index set $S$. 
We also define the following neighborhoods of the true parameters: for $\bar s_\star $ and $\bar\epsilon_n$ satisfying \ref{c2con:thm2mod1}, and sufficiently large constants $\tilde M_1$ and $\tilde M_2$,
\begin{align}
\begin{split}
\widetilde\Theta_n&=\left\{\theta\in\Theta: s_\theta\le K_1 \bar s_\star , \, \lVert X(\theta-\theta_0)\rVert_2\le  \tilde M_1 \sqrt{n}\bar\epsilon_n\right\},\\
\widetilde{\cal H}_n&=\left\{\eta\in{\cal H} : d_n(\eta,\eta_0)\le \tilde M_2  \bar\epsilon_n \right\}.
\end{split}
\label{c2eqn:subsets0}
\end{align}
Combined by other conditions,
Theorem~\ref{c2thm:thetacon} implies that the posterior probabilities of these neighborhoods tend to one in probability if $s_0\log p\lesssim n\bar\epsilon_n^2$. We need some bounding conditions on these neighborhoods, which will be specified below.

Let $\Phi(\eta)=(\tilde\xi_\eta,\tilde\Delta_\eta)$ for any given $\eta\in\mathcal H$.
For a given $\theta$, we choose a bijective map $\eta\mapsto\tilde\eta_n(\theta,\eta):\mathcal H\mapsto \mathcal H$ such that
$\Phi(\tilde\eta_n(\theta,\eta))=(\tilde\xi_\eta+H\tilde X(\theta-\theta_0),\tilde\Delta_\eta)$ for some orthogonal projection $H$ which may depend on the true parameter values, but not on $\theta$ and $\eta$. The projection $H$ plays a key role here and for a distributional approximation in the following section, and thus should be appropriately chosen to satisfy the followings.

\begin{enumerate}[label=(C\arabic*),leftmargin=3.0\parindent,resume, itemsep=0.0cm]
	\item \label{c2con:thm4asm00} The orthogonal projection $H$ satisfies
	\begin{align*}
	\frac{1}{(s_0\vee 1)\log p}\sup_{\eta\in\widetilde{\cal H}_n}\lVert (I-H)(\tilde\xi_{\eta}-\tilde\xi_{\eta_0})\rVert_2^2&\rightarrow 0,\\
	\min_{S:s\le K_1\bar s_\star }\inf_{v\in\mathbb{R}^s : \lVert v\rVert_2=1}\frac{\lVert(I-H) \tilde X_S v\rVert_2}{\lVert \tilde X_S v\rVert_2}&\gtrsim 1.
	\end{align*}
	\item \label{c2con:thm4asm10} The conditional law $\Pi_{n,\theta}$ of $\tilde\eta_n(\theta,\eta)$ given $\theta$, induced by the prior, is absolutely continuous relative to its distribution $\Pi_{n,\theta_0}$ at $\theta=\theta_0$ (which is the same as the prior for $\eta$), and the Radon-Nikodym derivative $d\Pi_{n,\theta}/d\Pi_{n,\theta_0}$ satisfies 
	\begin{align*}
	\sup_{\theta\in\widetilde\Theta_n}\sup_{\eta\in\widetilde{\cal H}_n}\left\lvert\log \frac{d\Pi_{n,\theta}}{d\Pi_{n,\theta_0}}(\eta)\right\rvert\lesssim 1.
	\end{align*}
\end{enumerate}

By reading the proof, one can see that Theorem~\ref{c2thm:thmmar} below is based on the approximate likelihood ratio. The first condition of \ref{c2con:thm4asm00} is required to control the remainder of an approximation.
The second condition of \ref{c2con:thm4asm00} implies that $\lVert u\rVert_2\lesssim\lVert(I-H) u\rVert_2\le\lVert u\rVert_2$ for every $u\in{\rm span}(\tilde X_S)$ with $S$ such that $s\le K_1\bar s_\star $, as the second inequality trivially holds by the fact that $I-H$ is an orthogonal projection. 
The use of the shifting map $\eta\mapsto\tilde\eta_n(\theta,\eta)$ is justified by the condition  \ref{c2con:thm4asm10}, which implies that a shift in certain directions does not substantially affect the prior on $\eta$.
This is related in spirit to the absolute continuity condition in the semiparametric Bernstein-von Mises theorem (see, for example, Theorem 12. 8 of \citet{ghosal2017fundamentals}).
We will see that a distributional approximation also requires similar, but stronger conditions.

Lastly, the complexity of the neighborhood $\widetilde {\mathcal H}_n$ should also be controlled. Specifically, we make the following condition.
\begin{enumerate}[label=(C\arabic*),leftmargin=3.0\parindent,resume, itemsep=0.0cm]
	\item \label{c2con:thm4asm20} For $a_n$ and $e_n$ satisfying \ref{c2con:maxfrob}  and a sufficiently large $C>0$, 
	\begin{align*}
	&\sqrt{\frac{n\bar\epsilon_n^2(e_n+a_n\bar\epsilon_n^2)}{(s_0\vee 1)\log p}}+\sqrt{a_n}\int_0^{C\bar\epsilon_n}\sqrt{\log N(\delta,\widetilde {\mathcal H}_n,d_{B,n})} d\delta\rightarrow 0.
	\end{align*}
	\item \label{c2con:sepsp} The parameter space ${\mathcal H}$ is separable with the pseudo-metric $d_{B,n}$.
\end{enumerate}
Similar to \ref{c2con:thm4asm00}, these conditions are required to control the remainder of an approximation. The integral term comes from the expected supremum of a separable Gaussian process, exploiting the Gaussian likelihood of the model and the separability of $\widetilde{\cal H}_n$ with the standard deviation metric.
Condition \ref{c2con:sepsp} is crucial for this reason. 
Since we usually put a prior on $\eta$ in an explicit way, condition \ref{c2con:sepsp} is rarely violated in practice.
One may see a connection between the first term of \ref{c2con:thm4asm20} and the conditions for Corollary~\ref{c2cor:opt}.
The former easily tends to zero even if $n\bar\epsilon_n^2/\log p$ is increasing, due to the extra term $\bar\epsilon_n$ which commonly tends to zero in a polynomial order.
Note also that the term $s_0\vee 1$ appears in \ref{c2con:thm4asm00} and \ref{c2con:thm4asm20}. Although this gives sharper bounds, the conditions often need to be verified with $s_0\vee 1$ replaced by $1$ as $s_0$ is unknown.

Under the conditions specified above, we obtain the following theorem for the contraction rates for $\theta$ which do not depend on $\bar\epsilon_n$. The compatibility numbers below should be understood to be 1 if $s_0=0$.
\begin{theorem}[Optimal posterior contraction]
	Suppose that {\rm\ref{c2con:maxfrob}--\ref{c2con:trueparcon2}}, {\rm\ref{c2con:thm2mod1}}, and {\rm\ref{c2con:sep}--\ref{c2con:sepsp}} are satisfied.
	Then, there exist constants $K_4$ and $K_5$ such that
	\begin{align}
	\begin{split}
	{\mathbb E}_0\Pi\left(\theta: s_\theta> K_4s_0 \,\Big|\,Y^{(n)} \right)&\rightarrow 0,\\
	{\mathbb E}_0\Pi\left(\theta:\lVert \theta-\theta_0\rVert_1>\frac{K_5s_0 \sqrt{\log p }}{\phi_1((K_4+1)s_0)\lVert X\rVert_{\ast}} \,\bigg|\,Y^{(n)} \right)&\rightarrow 0,\\
	{\mathbb E}_0\Pi\left(\theta:\lVert \theta-\theta_0\rVert_2>\frac{K_5 \sqrt{s_0\log p }}{\phi_2((K_4+1)s_0)\lVert X\rVert_{\ast}} \,\bigg|\,Y^{(n)} \right)&\rightarrow 0,\\
	{\mathbb E}_0\Pi\left(\theta:\lVert X(\theta-\theta_0)\rVert_2>K_5 \sqrt{s_0 \log p } \,\big |\,Y^{(n)} \right)&\rightarrow 0.
	\end{split}
	\label{c2eqn:thmmar}
	\end{align}
	\label{c2thm:thmmar}
\end{theorem}

Similar to the paragraph followed by Theorem~\ref{c2thm:thetacon}, the compatibility numbers are easily bounded away from zero so that they can be removed from the expressions. These are actually weaker than before as $s_0\le s_\star$.
The simplified rates are then available  for ease of interpretation.

\begin{remark}
	In regression models where no additional mean part $\xi_{\eta,i}$ exists, conditions \ref{c2con:thm4asm00} and \ref{c2con:thm4asm10} are trivially satisfied by choosing the zero matrix for $H$.
	This is also true for \ref{c2con:thm4asm0} and \ref{c2con:thm4asm1} specified in the next section.
	\label{rmk:1.5}
\end{remark}

\begin{remark}
	Suppose that there exists a matrix $Z_i\in\mathbb{R}^{m_i\times q}$ such that $\xi_{\eta,i}=Z_i h(\eta)$ for every $\eta\in{\cal H}$ with some map $h:{\cal H}\mapsto\mathbb{R}^{q}$. Then, a general strategy to choose $H$ is to set $H=\tilde Z(\tilde Z^T \tilde Z)^{-1}\tilde Z^T$ for $\tilde Z=(Z_1^T \Delta_{\eta_0,1}^{-1/2},\dots,Z_n^T \Delta_{\eta_0,n}^{-1/2})^T$. In this case, by the triangle inequality, the first condition of \ref{c2con:thm4asm00} is satisfied if there exists $\eta_\ast\in{\cal H}$ such that $n d_{A,n}^2(\eta_\ast,\eta_0)/(s_0\log p)\rightarrow 0$. For \ref{c2con:thm4asm0} in the next section, this is replaced by $(s_\star^2\log p) n d_{A,n}^2(\eta_\ast,\eta_0)\rightarrow 0$.
	These are trivially the case if there exists $\eta'\in{\cal H}$ such that $d_{A,n}(\eta',\eta_0)=0$.
	Also similar to Remark~\ref{rmk:1}, a sufficient condition for the second line of \ref{c2con:thm4asm00} is $\min\{\varsigma_{\min}([X_S,Z]):{s\le K_1 \bar s_\star }\}\gtrsim 1$ as pre-multiplication of a positive definite matrix by $X_S$ and $Z$ is an isomorphism. This is also sufficient for  \ref{c2con:thm4asm0} in the next section with $\bar s_\star$ replaced by $s_\star$.
	\label{rmk:2}
\end{remark}

\begin{remark}
	In many instances, for every $\delta>0$ and $\zeta_n>0$, we typically have
	\begin{align*}
	\log N\left(\delta,\{ \eta\in\mathcal H : d_{B,n}(\eta,\eta_0)\le \zeta_n \},d_{B,n}\right) \le 0\vee r_n\log\left(\frac{b_n\zeta_n}{\delta}\right),
	\end{align*}
	for some sequences $r_n$ and $b_n$, especially when the part of $\eta$ involved with $d_{B,n}$ is an $r_n$-dimensional Euclidean parameter.
	Note that $	\int_0^{C\zeta_n}\sqrt{0\vee r_n\log(b_n\zeta_n/{\delta})}d\delta$ is equal to
	\begin{align*}
	&\int_0^{(C\wedge b_n)\zeta_n}\sqrt{r_n\log\left(\frac{b_n\zeta_n}{\delta}\right)}d\delta \\
	&\quad=(C\wedge b_n)\zeta_n\sqrt{r_n\log\left(\frac{b_n}{C\wedge b_n}\right)}+b_n \zeta_n\sqrt{r_n}\int_{\sqrt{\log (b_n/(C\wedge b_n))}}^\infty e^{-t^2}dt.
	\end{align*}
	If $b_n$ is increasing, the right hand side is bounded by a multiple of $\zeta_n \sqrt{r_n\log b_n}$ by the tail probability of a normal distribution, while it is bounded by a multiple of $\zeta_n b_n\sqrt{r_n}$ for nonincreasing $b_n$. 
	This simplification is useful to verify \ref{c2con:thm4asm20} in many applications, and can also be used for \ref{c2con:thm4asm2} in the next section.
	\label{rmk:3}
\end{remark}

\section{Bernstein-von Mises and selection consistency}
\label{c2sec:bvm}
An extremely important question is whether the true support $S_0$ is recovered with probability tending to one, which is the property called selection consistency.
We will show this based on a distributional approximation to the posterior distribution.
Combined with selection consistency, the shape approximation also leads to the product of a point mass and a normal distribution, which we call the Bernstein-von Mises theorem. This reduced approximate distribution enables us to correctly quantify the remaining uncertainty of the parameter through the posterior distribution.

\subsection{Shape approximation to the posterior distribution}

It is worth noting that selection consistency can often be verified without a distributional approximation.
For example, in sparse linear regression with scalar unknown variance $\sigma^2$, \citet{song2017nearly} deployed the marginal likelihood of the model support which can be obtained by integrating out $\theta$ and $\sigma^2$ from the likelihood using the inverse gamma kernel. 
In our general formulation, however, this approach is hard to implement due to the arbitrary structure of a nuisance parameter $\eta$. Indeed, the approach is not directly available even for a parametric covariance matrix with dimension $\overline m\ge2$.
In this sense, using a shape approximation could be a natural solution to the problem, which may require some extra conditions on the parameter space and on the priors for $\theta$ and $\eta$.

Recall that the results in Section~\ref{sec:marcon} are based on the semiparametric theory. In this section we will need very similar conditions as before, but the requirements are generally stronger, as the remainder of an approximation should be strictly manipulated. Since the setup is high-dimensional, our conditions are even more restrictive than those for semiparametric models with a fixed dimensional parametric segment \citep[e.g.,][]{castillo2012semiparametric}. One may refer to Section~3.3 of \citet{chae2016bayesian} for a relevant discussion.

Throughout this section, we only consider $s_\star$ that satisfies the conditions of Theorem~\ref{c2thm:thetacon}.
First of all, we make a modification of \ref{c2con:combo1}. The following condition is slightly stronger than \ref{c2con:combo1}, but is still not too restrictive as \ref{c2con:thm2mod1} requires $s_\star =o(n)$.
\begin{enumerate}[label=(C\arabic*),leftmargin=3.0\parindent,start=7, itemsep=0.0cm]
	\item[\mylabel{c2con:combo1r}{(C7{$^\ast$})}] 
	Condition~\ref{c2con:combo1} is satisfied with $\bar s_\star$ replaced by $s_\star$.
\end{enumerate}
The assumption on the prior for $\theta$ is made only through the regularization parameter $\lambda$. 
As in \citet{castillo2015bayesian}, $\lambda$ should not increase too fast and should satisfy $\lambda s_\star\sqrt{\log p}/\lVert X\rVert_\ast\rightarrow 0$.
In our setup, the range of $\lambda$ induces a sufficient condition for this: $s_\star^2 \log p=o(n)$.
Since this is weaker than the one that will be made later in this section, the ``small lambda regime'' is automatically met by a stronger condition for the entire procedure for a distributional approximation (see \ref{c2con:thm4asm2} below and the following paragraph).

For sufficiently large constants $\hat M_1$ and $\hat M_2$, we now define the neighborhoods,
\begin{align}
\begin{split}
\widehat\Theta_n&=\Big\{\theta\in\Theta: s_\theta\le K_1 s_\star , \, \lVert \theta-\theta_0\rVert_1\le  \hat M_1 s_\star\sqrt{\log p}/\lVert X\rVert_\ast\Big\},\\
\widehat{\cal H}_n&=\Bigg\{\eta\in{\cal H} : d_{A,n}(\eta,\eta_0)\le \hat M_2  s_\star\sqrt{\frac{\log p}{n}},  d_{B,n}(\eta,\eta_0)\le \hat M_2 \sqrt{\frac{s_\star\log p}{n}}\Bigg\}.
\end{split}
\label{c2eqn:subsets}
\end{align}
Note that $\widehat\Theta_n$ is defined with an $\ell_1$-ball, which makes it contract more slowly than $\widetilde\Theta_n$ in \eqref{c2eqn:subsets0} under \ref{c2con:combo1r}. This is due to technical reasons that for a distributional approximation, the $\ell_1$-ball should be directly manipulated in the complement of $\widehat\Theta_n$. The neighborhood $\widehat{\cal H}_n$ is also increased to be matched with $\widehat\Theta_n$. We leave more details on this to the reader; refer to the proof of Theorem~\ref{c2thm:bvm} below.

As in Section~\ref{sec:marcon}, we choose a bijective map $\eta\mapsto\tilde\eta_n(\theta,\eta)$ which gives rise to
$\Phi(\tilde\eta_n(\theta,\eta))=(\tilde\xi_\eta+H\tilde X(\theta-\theta_0),\tilde\Delta_\eta)$ for some orthogonal projection $H$.
Again, the orthogonal projection $H$ should be carefully chosen to satisfy some boundedness conditions. The conditions are similar to, but stronger than those in Section~\ref{sec:marcon}.
This is not only because of the increased neighborhoods $\widehat\Theta_n$ and $\widehat{\cal H}_n$, but also because the remainder of an approximation should be bounded on their complements. We precisely make the required conditions below.

\begin{enumerate}[label=(C\arabic*),leftmargin=3.0\parindent,resume, itemsep=0.0cm]
	\item[\mylabel{c2con:thm4asm0}{(C8{$^\ast$})}]  The orthogonal projection $H$ satisfies
	\begin{align*}
	s_\star^2 \log p\sup_{\eta\in\widehat{\cal H}_n}\lVert (I-H)(\tilde\xi_{\eta}-\tilde\xi_{\eta_0})\rVert_2^2&\rightarrow 0,\\
	\min_{S:s\le K_1 s_\star }\inf_{v\in\mathbb{R}^s : \lVert v\rVert_2=1}\frac{\lVert(I-H) \tilde X_S v\rVert_2}{\lVert \tilde X_S v\rVert_2}&\gtrsim 1.
	\end{align*}
	\item[\mylabel{c2con:thm4asm1}{(C9{$^\ast$})}] The conditional law $\Pi_{n,\theta}$ of $\tilde\eta_n(\theta,\eta)$ given $\theta$, induced by the prior, is absolutely continuous relative to its distribution $\Pi_{n,\theta_0}$ at $\theta=\theta_0$, and the Radon-Nikodym derivative $d\Pi_{n,\theta}/d\Pi_{n,\theta_0}$ satisfies 
	\begin{align*}
	\sup_{\theta\in\widehat\Theta_n}\sup_{\eta\in\widehat{\cal H}_n}\left\lvert\log \frac{d\Pi_{n,\theta}}{d\Pi_{n,\theta_0}}(\eta)\right\rvert\rightarrow 0.
	\end{align*}
	\item[\mylabel{c2con:thm4asm2}{(C10{$^\ast$})}] 
	For $a_n$ and $e_n$ satisfying \ref{c2con:maxfrob} and a sufficiently large $C>0$, 
	\begin{align*}
	s_\star \log p \Bigg\{&s_\star\sqrt{e_n+\frac{a_ns_\star \log p}{n}}\\
	&+
	\sqrt{a_n } \int_0^{C\sqrt{(s_\star \log p)/n}}\sqrt{\log N\left(\delta,\widehat{\cal H}_n,d_{B,n}\right)}d\delta\Bigg\}\rightarrow 0.
	\end{align*}
\end{enumerate}

Conditions \ref{c2con:thm4asm0}--\ref{c2con:thm4asm2} are required for similar reasons as in Section~\ref{sec:marcon}.
We mention that \ref{c2con:thm4asm2} is a sufficient condition for the small lambda regime, since its necessary condition is $s_\star^5 \log^3 p=o(n)$ that is stronger than $s_\star^2 \log p=o(n)$. This necessary condition for \ref{c2con:thm4asm2} is often a sufficient condition in many finite dimensional models.

We define the standardized vector,
\begin{align*}
U=\left(
\begin{matrix}
\Delta_{\eta_0,1}^{-1/2}(Y_1-X_1\theta_0-\xi_{\eta_0,1}) \\ \vdots \\  \Delta_{\eta_0,n}^{-1/2}(Y_n-X_n\theta_0-\xi_{\eta_0,n})
\end{matrix}\right)\in\mathbb R^{n_\ast}.
\end{align*}
Under the assumptions above, the posterior distribution of $\theta$ is approximated by $\Pi^\infty$ given by
\begin{align}
\Pi^\infty(\theta\in\cdot\,|\,Y^{(n)})&=\sum_{S: s\le K_1 s_\star } \hat w_S \left({\cal N}
_{\hat\theta_S,\tilde X_S^T(I-H)\tilde X_S}^S \otimes\delta_0^{S^c}\right)(\theta\in\cdot),
\label{c2eqn:norapp}
\end{align}
where $\mathcal N_{\mu,\Omega}^S$ is the Gaussian measure with mean $\mu$ and precision $\Omega$ on the coordinate $S$, $\delta_0^{S^c}$ is the Dirac measure at zero on $S^c$, $\hat\theta_S$ is the least squares solution
$\hat\theta_S=(\tilde X_S^T(I-H)\tilde X_S)^{-1}\tilde X_S^T(I-H)(U+\tilde X \theta_0)$, and the weights $\hat w_S$ satisfy
\begin{align*}
\hat w_S\propto\frac{{\pi_p(s)}}{{\binom p {s}}}\left(\frac{\lambda}{2}\right)^{s} (2\pi)^{s/2}\det\Big(\tilde X_S^T(I-H)\tilde X_S\Big)^{-1/2}\exp \bigg\{\frac{1}{2}\lVert(I-H)\tilde X_S\hat\theta_S \rVert_2^2\bigg\}.
\end{align*}
Another way to express $\Pi^\infty$, for any measurable ${\cal B}\subset \mathbb{R}^p$, is
\begin{align*}
\Pi^\infty(\theta\in {\cal B}\,|\,Y^{(n)})&=\frac{\sum_{S: s\le K_1 s_\star }{\pi_p(s)}{\binom p {s}}^{-1}\left({\lambda}/{2}\right)^{s}\int_{\cal B} \Lambda_n^\star(\theta) d\{\mathcal L(\theta_S) \otimes\delta_0(\theta_{S^c})\}}{\sum_{S: s\le K_1 s_\star }{\pi_p(s)}{\binom p {s}}^{-1}\left({\lambda}/{2}\right)^{s}\int_{\mathbb{R}^p} \Lambda_n^\star(\theta) d\{\mathcal L(\theta_S) \otimes\delta_0(\theta_{S^c})\}},
\end{align*}
where $\mathcal L$ denotes the Lebesgue measure and 
\begin{align}
\Lambda_n^\star(\theta)=\exp\left\{{-\frac{1}{2}\lVert (I-H)\tilde X(\theta-\theta_{0})\rVert_2^2+U^T(I-H)\tilde X(\theta-\theta_{0})}\right\}.
\label{c2eqn:applambda}
\end{align}
It can be easily checked that both the expressions are equivalent.
The results are summarized in the following theorem.

\begin{theorem}[Distributional approximation]
	Suppose that {\rm\ref{c2con:maxfrob}--\ref{c2con:trueparcon2}}, {\rm\ref{c2con:thm2mod1}},  {\rm\ref{c2con:sep}}, {\rm\ref{c2con:combo1r}--\ref{c2con:thm4asm2}}, and {\rm\ref{c2con:sepsp}} are satisfied for some orthogonal projection $H$.
	Then
	\begin{align}
	\mathbb{E}_0\left\lVert\Pi(\theta\in\cdot\,|\,Y^{(n)})-\Pi^\infty(\theta\in\cdot\,|\,Y^{(n)})\right\rVert_{\rm TV}\rightarrow 0.
	\label{c2eqn:bvm}
	\end{align}
	\label{c2thm:bvm}
\end{theorem}

\subsection{Model selection consistency}

The shape approximation to the posterior distribution facilitates obtaining the next theorem which shows that the posterior distribution is concentrated on subsets of the true support with probability tending to one.
The result is then used as the basis of selection consistency.
Similar to the literature, the theorem requires an additional condition on the prior as follows.

\begin{enumerate}[label=(C\arabic*),leftmargin=3.0\parindent,start=12, itemsep=0.0cm]
	\item \label{c2con:selpri} The prior satisfies $A_4>1$ and $s_\star \lesssim p^a$ for $a<A_4-1$.
\end{enumerate}

\begin{theorem}[Selection, no supersets]
	Suppose that {\rm\ref{c2con:maxfrob}--\ref{c2con:trueparcon2}}, {\rm\ref{c2con:thm2mod1}},  {\rm\ref{c2con:sep}}, {\rm\ref{c2con:combo1r}--\ref{c2con:thm4asm2}}, and {\rm\ref{c2con:sepsp}--\ref{c2con:selpri}}  are satisfied for some orthogonal projection $H$.
	Then
	\begin{align}
	{\mathbb E}_0\Pi\left(\theta:S_\theta\supset S_0, S_\theta\ne S_0\,|\,Y^{(n)}\right)\rightarrow 0.
	\label{c2eqn:sel}
	\end{align}
	\label{c2thm:sel}
\end{theorem}

Since coefficients that are too close to zero cannot be identified by any selection strategy, some threshold for the true nonzero coefficients is needed for detection. The requirement of a threshold is a fundamental limitation in high-dimensional setups. We make the following threshold, the so-called beta-min condition. The condition is made in view of the third assertion of Theorem~\ref{c2thm:thmmar}. The second assertion can also be used to make a similar threshold, but we only consider the given one below as it is generally weaker.

\begin{enumerate}[label=(C\arabic*),leftmargin=3.0\parindent,resume, itemsep=0.0cm]
	\item \label{c2con:betamin} The true parameter satisfies
	\begin{align*}
	\min_{\theta_{0,j}\ne 0} |\theta_{0,j}| >\frac{K_5 \sqrt{s_0\log p}}{\phi_2((K_4+1)s_0)\lVert X\rVert_\ast}.
	\end{align*}
\end{enumerate}

Since Theorem~\ref{c2thm:thetacon} implies that the posterior distribution of the support of $\theta$ includes that of the true support with probability tending to one, selection consistency is an easy consequence of Theorem~\ref{c2thm:sel} under the beta-min condition \ref{c2con:betamin}. Moreover, this improves the distributional approximation in \eqref{c2eqn:bvm} so that the posterior distribution can be approximated by a single component of the mixture; that is, the Bernstein-von Mises theorem holds for the parameter component $\theta_{S_0}$.
The arguments here are summarized in the following two corollaries, whose proofs are straightforward and thus are omitted.
\begin{corollary}[Selection consistency]
	Suppose that {\rm\ref{c2con:maxfrob}--\ref{c2con:trueparcon2}}, {\rm\ref{c2con:thm2mod1}},  {\rm\ref{c2con:sep}}, {\rm\ref{c2con:combo1r}--\ref{c2con:thm4asm2}}, and {\rm\ref{c2con:sepsp}--\ref{c2con:betamin}}    are satisfied for some orthogonal projection $H$.
	Then
	\begin{align}
	{\mathbb E}_0\Pi\left(\theta: S_\theta\ne S_0\,|\,Y^{(n)}\right)\rightarrow 0.
	\label{c2eqn:strsel}
	\end{align}
	\label{c2cor:strsel}
\end{corollary}
\begin{corollary}[Bernstein-von Mises]
	Suppose that {\rm\ref{c2con:maxfrob}--\ref{c2con:trueparcon2}}, {\rm\ref{c2con:thm2mod1}},  {\rm\ref{c2con:sep}}, {\rm\ref{c2con:combo1r}--\ref{c2con:thm4asm2}}, and {\rm\ref{c2con:sepsp}--\ref{c2con:betamin}}   are satisfied for some orthogonal projection $H$.
	Then
	\begin{align}
	\begin{split}
	{\mathbb E}_0 \bigg\lVert&\Pi(\theta\in\cdot\,|\,Y^{(n)})-\left( {\cal N}_{\hat\theta_{S_0},\tilde X_{S_0}^T(I-H)\tilde X_{S_0}}^S\otimes\delta_0^{S_0^c} \right) (\theta\in\cdot)\bigg\rVert_{\rm TV}\rightarrow 0.
	\end{split}
	\label{c2eqn:strbvm}
	\end{align}
	\label{c2cor:strbvm}
\end{corollary}

Corollary~\ref{c2cor:strbvm} enables us to quantify the remaining uncertainty of the parameter through the posterior distribution.
Specifically, we can construct credible sets for the individual components of $\theta_0$ as in \citet{castillo2015bayesian}.
It is easy to see that by the definition of $\hat\theta_{S_0}$, its $j$th component has a normal distribution, whose mean is the $j$th element of $\theta_{S_0}$ and variance is the $j$th diagonal element of $(\tilde X_{S_0}^T(I-H)\tilde X_{S_0})^{-1}$. 
Correct uncertainty quantification is thus guaranteed by the weak convergence.

\section{Applications}
\label{c2sec:app}

In this section, we apply the main results established in this study to the examples considered in Section~\ref{sec:example}. The main objective is to obtain nearly optimal posterior contraction rates and selection consistency via shape approximation to the posterior distribution with the Bernstein-von Mises phenomenon.

To use Corollary~\ref{c2cor:opt} for the optimal posterior contraction when $n\bar\epsilon_n^2=\log n$,  we simply assume that $s_0>0$ for all examples in this section, although Theorem~\ref{c2thm:thmmar} can also be applied under stronger conditions. The assumption $s_0>0$ is extremely mild rather than considering the ultra high-dimensional case, i.e., $\log n=o(\log p)$.
A large enough $A_4$ is also sufficient instead of the assumption $s_0>0$, but we do not pursue this direction as a specific threshold is not available.
We check the conditions of Theorem~\ref{c2thm:thmmar} only for more complicated models where $n\bar\epsilon_n^2>\log n$.

\subsection{Multiple response models with missing components}
We first apply the main results to Example~\ref{c2exm:missing}.
To recover posterior contraction of $\Sigma$ from the primitive results, it is necessary to assume that every entry of the response is jointly observed sufficiently many times. To be more specific, let $e_{ij}$ be 1 if the $j$th entry of $Y_i^{\rm aug}$ is observed and be zero otherwise. The contraction rate of the $(j,k)$th element of $\Sigma$ is directly determined by the order of $n^{-1}\sum_{i=i}^ne_{ij}e_{ik}$. The ideal case is when this quantity is bounded away from zero, that is, the entries are jointly observed at a rate proportional to $n$. Then the recovery is possible without any loss of information. If $n^{-1}\sum_{i=1}^ne_{ij}e_{ik}$ decays to zero, then the optimal recovery is not attainable, but consistent estimation may still be possible with slower rates. With an inverse Wishart prior on $\Sigma$, the following theorem studies the posterior asymptotic properties of the given model.

\begin{theorem}
	Assume that $s_0>0$, $1\lesssim\rho_{\min}(\Sigma_0)\le\rho_{\max}(\Sigma_0)\lesssim 1$, $\lVert \theta_0\rVert_\infty\lesssim \lambda^{-1} \log p$, and $\min_{j,k} n^{-1}\sum_{i=1}^n e_{ij}e_{ik}\gtrsim c_n^{-1}$ for some nondecreasing $c_n$ such that $c_n s_0 \log p=o(n)$. Then the following assertions hold.
	\begin{enumerate}[label={\rm (\alph*)}]
		\item	\label{c2thm:missingcon} The optimal posterior contraction rates for $\theta$ in \eqref{c2eqn:thmmar} are obtained.
		\item  \label{c2thm:missingcon2}
		The posterior contraction rate for $\Sigma$ is $\sqrt{c_n (s_0\log p)/n}$ with respect to the Frobenius norm. 
	\end{enumerate} 
	Assume further that $c_n(s_0 ^2\vee\log c_n)(s_0 \log p)^3=o(n)$ and  $\phi_1(Ds_0 )\gtrsim1$ for a sufficiently large $D$. Then the following assertions hold.
	\begin{enumerate}[label={\rm (\alph*)}, resume]
		\item	\label{c2thm:missingbvm} 
		For $H\in\mathbb{R}^{n_\ast\times n_\ast}$ the zero matrix, the distributional approximation in \eqref{c2eqn:bvm} holds.
		\item \label{c2thm:missingbvm2}  If $A_4>1$ and $s_0 \lesssim p^a$ for $a<A_4-1$, then the no-superset result in \eqref{c2eqn:sel} holds.
		\item \label{c2thm:missingbvm3}  Under the beta-min condition as well as the conditions for {\rm \ref{c2thm:missingbvm2}}, the selection consistency in \eqref{c2eqn:strsel} and the Bernstein-von Mises theorem in \eqref{c2eqn:strbvm} hold.
	\end{enumerate}
	\label{c2thm:exmmissing}
\end{theorem}


\subsection{Multivariate measurement error models}
We now consider Example~\ref{c2exm:mem}.
For convenience we write $Y^\ast=(Y_1^\ast,\dots Y_n^\ast)^T\in\mathbb{R}^n$, $W=(W_1^T,\dots,W_n^T)^T\in\mathbb{R}^{nq}$, and $X^\ast=(X_1^{\ast},\dots,X_n^{\ast})^T\in\mathbb{R}^{n\times p}$ in what follows.
In this subsection, we use the symbol $\otimes$ for the Kronecker product of matrices.
For priors of the nuisance parameters, normal prior distributions are assigned for the location parameters ($\alpha$, $\beta$, and $\mu$) and an inverse gamma and inverse Wishart prior are used for the scale parameters ($\sigma^2$ and $\Sigma$). The next theorem shows posterior asymptotic properties of the model.
In particular, specific forms of their mean and variance for shape approximation are provided considering the modeling structure.

\begin{theorem}
	Assume that $s_0>0$, $s_0\log p=o(n)$, $|\alpha_0|\vee\lVert\beta_0\rVert_\infty\vee\lVert\mu_0\rVert_\infty\lesssim 1$, $1\lesssim\sigma_0^2 \lesssim 1$, $1\lesssim\rho_{\min}(\Sigma_0)\le\rho_{\max}(\Sigma_0)\lesssim 1$, $\lVert \theta_0\rVert_\infty\lesssim \lambda^{-1} \log p$,
	and $\min\{\varsigma_{\min}([X_S^\ast,1_n]):s\le D s_0 \}\gtrsim 1$ for a sufficiently large $D$. Then the following assertions hold.
	\begin{enumerate}[label={\rm (\alph*)}]
		\item \label{c2thm:memcon}	The optimal posterior contraction rates for $\theta$ in \eqref{c2eqn:thmmar} are obtained.
		\item \label{c2thm:memcon2} The contraction rates for $\alpha$, $\beta$, $\mu$, and $\sigma^2$ are $\sqrt{(s_0\log p )/n}$  relative to the $\ell_2$-norms. The same rate is also obtained for $\Sigma$ with respect to the Frobenius norm.
	\end{enumerate}
	Assume further that $s_0 ^5\log^3 p=o(n)$ and $\phi_1(Ds_0 )\gtrsim1$ for a sufficiently large $D$. Then the following assertions hold.
	\begin{enumerate}[label={\rm (\alph*)}, resume]
		\item \label{c2thm:membvm} The distributional approximation in \eqref{c2eqn:bvm} holds with the mean vector 
		\begin{align*}
		\hat\theta_S=&(X_S^{\ast T} H^\ast X_S^\ast)^{-1}X_S^{\ast T}\Big\{H^\ast\Big[ \left(Y^\ast-(\alpha_0+\mu_0^T\beta_0)1_n\right)\\
		&\qquad\qquad\qquad\qquad-\left(I_n\otimes(\beta_0^T\Sigma_0(\Sigma_0+\Psi)^{-1})\right)\left(W-1_n\otimes\mu_0\right)\Big]\Big\}
		\end{align*}
		and the covariance matrix
		$(\sigma_0^2+\beta_0^T\Sigma_0(\Sigma_0+\Psi)^{-1}\Psi\beta_0) (X_S^{\ast T} H^\ast X_S^\ast)^{-1}$ for $H^\ast=I_n-n^{-1}1_n1_n^T$. 
		\item \label{c2thm:membvm2} If  $A_4>1$ and $s_0 \lesssim p^a$ for $a<A_4-1$, then the no-superset result in \eqref{c2eqn:sel} holds.
		\item \label{c2thm:membvm3} Under the beta-min condition as well as the conditions for {\rm \ref{c2thm:membvm2}}, the selection consistency in \eqref{c2eqn:strsel} and the Bernstein-von Mises theorem in \eqref{c2eqn:strbvm} hold.
	\end{enumerate}
	\label{c2thm:exmmem}
\end{theorem}
We note that the marginal law of $W_i$ is given by $W_i\sim {\rm N}(\mu,\Sigma+\Psi)$. This gives a hope that the rates for $\mu$ and $\Sigma$ may actually be improved up to the parametric rate $n^{-1/2}$ (possibly up to some logarithmic factors). However, other parameters are connected to the high-dimensional coefficients $\theta$, so such a parametric rate may not be obtained for them.

\subsection{Parametric correlation structure}
\label{c2sec:paramcor}
Next, our main results are applied to Example~\ref{c2exm:parcor}.
A correlation matrix $G_i(\alpha)$ should be chosen so that the conditions in the main theorems can be satisfied. Here we consider a compound-symmetric, a first order autoregressive, and a first order moving average correlation matrices: for $\alpha\in (b_1,b_2)$ with fixed boundaries $b_1$ and $b_2$ of the range, respectively, $\{G_i^{\rm CS}(\alpha)\}_{j,k}=\mathbbm{1}(j=k)+\alpha\mathbbm{1}(j\ne k)$, $\{G_i^{\rm AR}(\alpha)\}_{j,k}=\alpha^{|j-k|}$, and $\{G_i^{\rm MA}(\alpha)\}_{j,k}=\mathbbm{1}(j=k)+\alpha\mathbbm{1}(|j-k|=1)$.
The range is chosen so that the corresponding correlation matrix can be positive definite, i.e., $(b_1,b_2)=(0,1)$ for $G_i^{\rm CS}(\alpha)$, $(b_1,b_2)=(-1,1)$ for $G_i^{\rm AR}(\alpha)$, and  $(b_1,b_2)=(-1/2,1/2)$ for $G_i^{\rm MA}(\alpha)$.
Again, an inverse gamma prior is assigned to $\sigma^2$. For a prior on $\alpha$, we consider a density 
\begin{align*}
\Pi(d\alpha)\propto \exp\left\{-\frac{1}{(\alpha-b_1)^{c_1}(b_2-\alpha)^{c_2}}\right\},\quad \alpha\in(b_1,b_2),
\end{align*}
for some $c_1,c_2>0$ such that $\Pi(\alpha<t)\lesssim \exp(-(t-b_1)^{-c_1})$ for $t>b_1$ close to $b_1$ and $\Pi(\alpha>t)\lesssim \exp(-(b_2-t)^{-c_2})$ for $t<b_2$ close to $b_2$.

\begin{theorem}
	Assume that $s_0>0$, $s_0\log p=o(n)$, $\overline m n\asymp n_\ast$, $\lVert \theta_0\rVert_\infty\lesssim \lambda^{-1} \log p$, $\sigma_0^2\asymp 1$, $\alpha_0\in[b_1+\epsilon,b_2-\epsilon]$ for some fixed $\epsilon>0$.
	Suppose further that $ \overline m\lesssim 1$ for the compound-symmetric correlation matrix and $\log \overline m\lesssim \log p$ for the autoregressive and moving average correlation matrices.
	Then the following assertions hold.
	\begin{enumerate}[label={\rm (\alph*)}]
		\item  \label{c2thm:parcorcon} For any correlation matrix discussed above, the optimal posterior contraction rates for $\theta$ in \eqref{c2eqn:thmmar} are obtained.
		\item \label{c2thm:parcorcon2} For the autoregressive and moving average correlation matrices,  the posterior contraction rates for $\sigma^2$ and $\alpha$ are $\sqrt{(s_0\log p)/(\overline m n)}$ with respect to the $\ell_2$-norms. For the compound-symmetric correlation matrix, their contraction rates are $\sqrt{(s_0\log p)/ n}$  relative to the $\ell_2$-norm.
	\end{enumerate}
	Assume further that $s_0 ^5\log^3 p=o(n)$ and  $\phi_1(Ds_0 )\gtrsim1$ for a sufficiently large $D$. Then the following assertions hold.
	\begin{enumerate}[label={\rm (\alph*)}, resume]
		\item  \label{c2thm:parcorbvm} For $H\in\mathbb{R}^{n_\ast\times n_\ast}$ the zero matrix, the distributional approximation in \eqref{c2eqn:bvm} holds.
		\item \label{c2thm:parcorbvm2} If $A_4>1$ and $s_\star \lesssim p^a$ for $a<A_4-1$, then the no-superset result in \eqref{c2eqn:sel} holds.
		\item \label{c2thm:parcorbvm3} Under the beta-min condition as well as the conditions for {\rm \ref{c2thm:parcorbvm2}}, the selection consistency in \eqref{c2eqn:strsel} and the Bernstein-von Mises theorem in \eqref{c2eqn:strbvm} hold.
	\end{enumerate}
	\label{c2thm:exmparcor}
\end{theorem}

As for the prior for $\alpha$, the property that the tail probabilities decay to zero exponentially fast near both zero and one is crucial for the optimal posterior contraction rates.
It should be noted that many common probability distributions with compact supports may not be enough for this purpose (e.g., beta distributions).

The main difference between this example and those in the preceding subsections is that we consider possibly increasing $m_i$ here. Although we have the same form of contraction rates for $\theta$ as in previous examples, the implication is not the same due to a different order of $\lVert X\rVert_\ast$. For increasing $m_i$, it is expected to have $\lVert X\rVert_\ast\asymp \sqrt{n_\ast}$, which is commonly the case in regression settings. This is reduced to $\lVert X\rVert_\ast\asymp \sqrt{n}$ for the cases with fixed $m_i$, and hence increasing $m_i$ may help get faster rates.
While the increasing dimensionality of $m_i$ is often a benefit for contraction properties of $\theta$, this may or may not be the case for the nuisance parameters since it depends on the dimensionality of $\eta$. In the example in this subsection, the dimension of the nuisance parameters is fixed although $m_i$ can increase, which makes their posterior contraction rates faster than those with fixed $m_i$. However, this may not be true if $\eta$ is increasing dimensional. For example, see the example in Section~\ref{c2sec:graph}.

\subsection{Mixed effects models}
For the mixed effects models with sparse regression coefficients in Example~\ref{c2exm:mm}, we assume that the maximum of $\lVert Z_i\rVert_{\rm sp}$ is bounded, which is particularly mild if $\overline m$ is bounded. We also assume that $\sum_{i=1}^n\mathbbm{1}(m_i\ge q)\asymp n$ and $\min_i\{\varsigma_{\min}( Z_i) :m_i\ge q\}\gtrsim 1$, that is, $m_i$ is likely to be larger than $q$ with fixed probability and $Z_i$ is a full rank. 
These conditions are required for \ref{c2con:maxfrob} to hold.
We put an inverse Wishart prior on $\Psi$ as in other examples. The following theorem shows that the posterior asymptotic properties of the mixed effects models.

\begin{theorem}
	Assume that $s_0>0$, $s_0\log p=o(n)$, 
	$1\lesssim\rho_{\min}(\Psi_0)\le\rho_{\max}(\Psi_0)\lesssim 1$, $\lVert \theta_0\rVert_\infty\lesssim \lambda^{-1} \log p$, $\sum_{i=1}^n\mathbbm{1}(m_i\ge q)\asymp n$,  $\min_i\{\varsigma_{\min}( Z_i) :m_i\ge q\}\gtrsim 1$, and  $\max_i\lVert Z_i\rVert_{\rm sp}\lesssim 1$. Then the following assertions hold.
	\begin{enumerate}[label={\rm (\alph*)}]
		\item \label{c2thm:mmcon} The optimal posterior contraction rates for $\theta$ in \eqref{c2eqn:thmmar} are obtained.
		\item \label{c2thm:mmcon2} The posterior contraction rate for $\Psi$ is $\sqrt{(s_0\log p)/n}$ with respect to the Frobenius norm. 
	\end{enumerate}
	Assume further that $s_0 ^5\log^3 p=o(n)$ and $\phi_1(Ds_0)\gtrsim1$ for a sufficiently large $D$. Then the following assertions hold.
	\begin{enumerate}[label={\rm (\alph*)}, resume]
		\item \label{c2thm:mmbvm} For $H\in \mathbb{R}^{n_\ast\times n_\ast}$ the  zero matrix, the distributional approximation in \eqref{c2eqn:bvm} holds.
		\item \label{c2thm:mmbvm2} If  $A_4>1$ and $s_0 \lesssim p^a$ for $a<A_4-1$, then the no-superset result in \eqref{c2eqn:sel} holds.
		\item \label{c2thm:mmbvm3} Under the beta-min condition as well as the conditions for {\rm \ref{c2thm:mmbvm2}}, the selection consistency in \eqref{c2eqn:strsel} and the Bernstein-von Mises theorem in \eqref{c2eqn:strbvm} hold.
	\end{enumerate}
	\label{c2thm:exmmm}
\end{theorem}

Note that we assume that $\sigma^2$ is known, which is actually unnecessary at the modeling stage.
The assumption was made to find a sequence $a_n$ satisfying \ref{c2con:maxfrob} with ease. This can be relaxed only with stronger assumptions on $Z_i$. For example, if $q=1$ and $Z_i$ is an all-one vector, then the model is equivalent to that with a compound-symmetric correlation matrix in Section~\ref{c2sec:paramcor} with some reparameterization, in which $\sigma^2$ can be treated as unknown.

\subsection{Graphical structure with sparse precision matrices}
\label{c2sec:graph}
For the graphical structure models in Example~\ref{c2exm:graph}, we define an edge-inclusion indicator $\Upsilon=\{\upsilon_{jk}:1\le j\le k\le \overline m\}$ such that $\upsilon_{jk}=1$ if $\omega_{jk}\neq0$ and $\upsilon_{jk}=0$ otherwise, where $\omega_{jk}$ is the $(j,k)$th element of $\Omega$. 
We put a prior with a density $f_1$ on $(0,\infty)$ to the nonzero off-diagonal entries and a prior with a density $f_2$ on $\mathbb{R}$ to the diagonal entries of $\Omega$, such that 
the support is truncated to a matrix space with restricted eigenvalues and entries.
For the edge-inclusion indicator, we use a binomial prior with probability $\varpi$ when $|\Upsilon| \coloneqq \sum_{j,k}\upsilon_{jk}$ is given, and assign a prior to $|\Upsilon|$ such that $\log\Pi(|\Upsilon| \le \bar r)\lesssim -\bar r\log \bar r$.
The prior specification is summarized as
\begin{align*}
\Pi(\Omega|\Upsilon) &\propto \prod_{j,k:\upsilon_{jk}=1} f_1( \omega_{jk})\prod_{j=1}^{\overline m}f_2(\omega_{jj})\mathbbm{1}_{{\cal M}_0^+(L)}(\Omega),\\
\Pi(\Upsilon) &\propto  \varpi^{\bar r}(1-\varpi)^{\binom{\overline m}{2}-\bar r} \Pi(|\Upsilon| = \bar r), \quad \log\Pi(|\Upsilon| \le \bar r)\lesssim -\bar r\log \bar r,
\end{align*}
where ${\cal M}_0^+(L)$ is a collection of $\overline m\times\overline m$ positive definite matrices for a sufficiently large $L$, in which eigenvalues are between $[L^{-1},L]$ and entries are also bounded by $L$ in absolute value.

\begin{theorem}
	Let $s_\star =s_0\vee \bar s_\star$ for $\bar s_\star=(\overline m+d)(\log n)/\log p$.
	Assume that $s_0>0$, $s_0\log p=o(n)$, $\overline m\log n=o(n)$, $|\Upsilon_0|\le d$ for some $d$ such that $d \log n=o(n)$, $\Omega_0\in{\cal M}_0^+(cL)$ for some $0<c<1$, and $\lVert \theta_0\rVert_\infty\lesssim \lambda^{-1} \log p$.
	Then the following assertions hold.
	\begin{enumerate}[label={\rm (\alph*)}]
		\item \label{c2thm:graphcon} The posterior contraction rates for $\theta$ are given by \eqref{c2eqn:thm3}. If $\bar s_\star\lesssim 1$, the optimal rates in \eqref{c2eqn:thmmar} are obtained.
		\item \label{c2thm:graphcon2} The posterior contraction rate of $\Omega$ is $\sqrt{(s_0\log p \vee (\overline m+d)\log n )/n}$ with respect to the Frobenius norm. 
	\end{enumerate}
	If further $(\bar s_\star \vee \overline m^2)\bar s_\star \log p=o(n)$ and $\phi_1(D\bar s_\star )\gtrsim1$ for a sufficiently large $D$, then the following assertion holds.
	\begin{enumerate}[label={\rm (\alph*)}, resume]
		\item \label{c2thm:graphcon3} The optimal posterior contraction rates for $\theta$ in \eqref{c2eqn:thmmar} are obtained even if $\bar s_\star\rightarrow \infty$.
	\end{enumerate}
	Assume further that  $(s_\star \vee \overline m)^2(s_\star \log p)^3=o(n)$ and $\phi_1(Ds_\star )\gtrsim1$ for a sufficiently large $D$.	Then the following assertions hold.
	\begin{enumerate}[label={\rm (\alph*)}, resume]
		\item  \label{c2thm:graphbvm} For $H\in\mathbb{R}^{n_\ast\times n_\ast}$ the zero matrix, the distributional approximation in \eqref{c2eqn:bvm} holds.
		\item \label{c2thm:graphbvm2} If   $A_4>1$ and $s_\star \lesssim p^a$ for $a<A_4-1$, then the no-superset result in \eqref{c2eqn:sel} holds.
		\item \label{c2thm:graphbvm3} Under the beta-min condition as well as the conditions for {\rm \ref{c2thm:graphbvm2}}, the selection consistency in \eqref{c2eqn:strsel} and the Bernstein-von Mises theorem in \eqref{c2eqn:strbvm} hold.
	\end{enumerate}
	\label{c2thm:exmgraph}
\end{theorem}

Note that  increasing $\overline m$ is likely to improve the $\ell_2$-norm contraction rate for $\theta$ as we expect that $\lVert X\rVert_\ast\asymp \sqrt{\overline m n}$. In particular, the improvement is clearly the case if $d\lesssim \overline m$ and $\phi_2(Ds_\star )\gtrsim1$ for a sufficiently large $D$.
However, as pointed out in Section~\ref{c2sec:paramcor}, this is not the case for $\Omega$ as its dimension is also increasing.

If we assume that $\log n\lesssim \log \overline m$, then the term $\sqrt{(\overline m+d)(\log n )/n}$ arising from the sparse precision matrix $\Omega$ becomes $\sqrt{(\overline m+d)(\log \overline m)/n}$. The latter is comparable to the frequentist convergence rate of the graphical lasso in \citet{rothman2008sparse}. Therefore, our rate is deemed to be optimal considering the additional complication due to the mean term involving sparse regression coefficients.

\subsection{Nonparametric heteroskedastic regression models}
\label{c2sec:exhet}
Next, we use the main results for Example~\ref{c2exm:het}.
For a bounded, convex subset ${\cal X}\subset\mathbb{R}$, define the $\alpha$-H\"older class $\mathfrak{C}^\alpha ({\cal X})$ as the collection of functions $f:{\cal X}\rightarrow \mathbb{R}$ such that $\lVert f\rVert_{\mathfrak{C}^\alpha}<\infty$, where
\begin{align*}
\lVert f\rVert_{\mathfrak{C}^\alpha}=\max_{0\le k\le\lfloor\alpha\rfloor} \sup_{x\in{\cal X}}|f^{(k)}(x)|+\sup_{x,y\in{\cal X}:x\ne y}\frac{|f^{(\lfloor\alpha\rfloor)}(x)-f^{(\lfloor\alpha\rfloor)}(y)|}{|x-y|^{\alpha-\lfloor\alpha\rfloor}},
\end{align*}
with the $k$th derivative $f^{(k)}$ of $f$ and $\lfloor\alpha\rfloor$ the largest integer that is strictly smaller than $\alpha$.
Let the true function $v_0$ belong to $\mathfrak{C}^\alpha [0,1]$ with assumption that $v_0$ is strictly positive. While $\alpha>1/2$ suffices for the basic posterior contraction, we will see that the optimal posterior contraction for $\theta$ requires $\alpha>1$. The stronger condition $\alpha>2$ is even needed for the Bernstein-von Mises theorem and the selection consistency, but all these conditions are mild if the true function is sufficiently smooth.

We put a prior on $g$ through B-splines. The function is expressed as a linear combination of $J$-dimensional B-spline basis terms $B_J$ of order $q\ge\alpha$, i.e., $v_\beta(z)=\beta^T B_J(z)$, while an inverse Gaussian prior distribution is independently assigned to each entry of $\beta$. 
For any measurable function $f:[0,1]\mapsto \mathbb R$, we let $\lVert f\rVert_\infty=\sup_{z\in[0,1]}|f(z)|$ and $\lVert f\rVert_{2,n}=(n^{-1}\sum_{i=1}^n|f(z_i)|^2)^{1/2}$ denote the sup-norm and empirical $L_2$-norm, respectively.
To deploy the properties of B-splines, we assume that $z_i$ are sufficiently regularly distributed on $[0,1]$.

\begin{theorem}
	The true function $v_0$ is assumed to be strictly positive on $[0,1]$ and belong to $\mathfrak{C}^\alpha [0,1]$ with $\alpha>1/2$.
	We choose $J\asymp (n/\log n)^{1/(2\alpha+1)}$.
	Let $s_\star =s_0\vee \bar s_\star$ for $\bar s_\star=(\log n)^{2\alpha/(2\alpha+1)}n^{1/(2\alpha+1)}/\log p$ and assume that $s_0>0$, $Js_0\log p=o(n)$, and $\lVert \theta_0\rVert_\infty\lesssim \lambda^{-1} \log p$.
	 Then the following assertions hold.
	\begin{enumerate}[label={\rm (\alph*)}]
		\item	\label{c2thm:hetcon} The posterior contraction rates for $\theta$ are given by \eqref{c2eqn:thm3}. If $\bar s_\star\lesssim 1$, the optimal rates in \eqref{c2eqn:thmmar} are obtained.
		\item	\label{c2thm:hetcon2} The posterior contraction rate for $v$ is $\sqrt{(s_0\log p)/n}\vee(\log n/n)^{\alpha/(2\alpha+1)}$ with respect to the $\lVert\cdot\rVert_{2,n}$-norm.
	\end{enumerate}	
	If further $\alpha>1$ and $\phi_1(D\bar s_\star )\gtrsim1$ for a sufficiently large $D$, then the following assertion holds.
	\begin{enumerate}[label={\rm (\alph*)}, resume]
		\item \label{c2thm:hetcon3} The optimal posterior contraction rates for $\theta$ in \eqref{c2eqn:thmmar} are obtained even if $\bar s_\star\rightarrow \infty$.
	\end{enumerate}
	Assume further that $\alpha>2$, $J(s_\star ^2\vee J)(s_\star \log p)^3 =o(n)$ and $\phi_1(Ds_\star )\gtrsim1$ for a sufficiently large $D$. Then the following assertions hold.
	\begin{enumerate}[label={\rm (\alph*)}, resume]
		\item \label{c2thm:hetbvm} The distributional approximation in \eqref{c2eqn:bvm} holds with $H$ the $n\times n$ zero matrix. 
		\item  \label{c2thm:hetbvm2}  If $A_4>1$ and $s_\star \lesssim p^a$ for $a<A_4-1$, then the no-superset result in \eqref{c2eqn:sel} holds.
		\item \label{c2thm:hetbvm3}  Under the beta-min condition as well as the conditions for {\rm \ref{c2thm:hetbvm2}}, the selection consistency in \eqref{c2eqn:strsel} and the Bernstein-von Mises theorem in \eqref{c2eqn:strbvm} hold.
	\end{enumerate}
	\label{c2thm:exmhet}
\end{theorem}

An inverse Gaussian prior is used due to the property that its tail probabilities at both zero and infinity decay to zero exponentially fast.
The exponentially decaying tail probabilities in both directions are essential to obtain the optimal contraction rate. Note that standard choices such as gamma and inverse gamma distributions do not satisfy this property.

By investigating the proof, it can be seen that the condition $\alpha>1/2$ is required to satisfy  condition \ref{c2con:maxfrob} for posterior contraction, so this condition is not avoidable in applying the main theorems.
Unlike Theorem~\ref{c2thm:exmplm} below, assertion \ref{c2thm:hetcon3} does not require any further boundedness condition. This is because the restriction $\alpha>1$ makes the required bound tend to zero.
For the Bernstein-von Mises theorem and the selection consistency, it can be seen that $\alpha>2$ is necessary for the condition $J(s_\star ^2\vee J)(s_\star \log p)^3 =o(n)$ but not sufficient. Although the requirement $\alpha>2$ is implied by the latter condition, we specify this in the statement due to its importance.  We refer to the proof of Theorem~\ref{c2thm:exmhet} for more details.

\subsection{Partial linear models}
\label{c2sec:explm}
Lastly, we consider Example~\ref{c2exm:plm}. 
We assume that the true function $g_0$ belongs to $\mathfrak{C}^\alpha [0,1]$ for with $\alpha>0$. Any $\alpha>0$ suffices for the basic posterior contraction, but stronger restrictions are required for further assertions as in Theorem~\ref{c2thm:exmhet}.
We put a prior on $g$ through $J$-dimensional B-spline basis terms of order $q\ge a$, i.e., $g_\beta(z)=\beta^T B_J(z)$. 
With a given $J$, we define the design matrix $W_J=(B_J(z_1),\dots, B_J(z_n))^T\in\mathbb{R}^{n\times J}$.
The standard normal prior is independently assigned to each component of $\beta$ and an inverse gamma prior is assigned to $\sigma^2$.
Similar to Section~\ref{c2sec:exhet}, we assume that $z_i$ are sufficiently regularly distributed on $[0,1]$.

\begin{theorem} 
	The true function is assumed to satisfy $g_0\in\mathfrak{C}^\alpha [0,1]$ with $\alpha>0$.
	We choose $J\asymp (n/\log n)^{1/(2\bar\alpha+1)}$ for some $\bar\alpha\le\alpha$.
	Let $s_\star =s_0\vee \bar s_\star$ for $\bar s_\star=(\log n)^{2\bar\alpha/(2\bar\alpha+1)}n^{1/(2\bar\alpha+1)}/\log p$ and 
	assume that $s_0>0$, $s_0\log p=o(n)$, $\sigma_0^2\asymp1$, $\lVert \theta_0\rVert_\infty\lesssim \lambda^{-1} \log p$, and $\min\{\varsigma_{\min}([X_S,W_J]):s\le D s_\star \}\gtrsim 1$ for a sufficiently large $D$.
	Then the following assertions hold.
	\begin{enumerate}[label={\rm (\alph*)}]
		\item \label{c2thm:plmcon} The posterior contraction rates for $\theta$ are given by \eqref{c2eqn:thm3}. If $\bar s_\star\lesssim 1$, the optimal rates in \eqref{c2eqn:thmmar} are obtained.
		\item \label{c2thm:plmcon2} The contraction rates for $g$ and $\sigma^2$ are $\sqrt{(s_0\log p)/n}\vee(\log n/n)^{\bar\alpha/(2\bar\alpha+1)}$ with respect to the $\lVert\cdot\rVert_{2,n}$- and $\ell_2$-norms, respectively.
	\end{enumerate}	
	If further $1/2\le \bar\alpha<\alpha$, $(\log n)^{2(\alpha\wedge 2\bar\alpha)/(2\bar\alpha+1)}n^{(-2(\alpha\wedge 2\bar\alpha)+2\bar\alpha+1)/(2\bar\alpha+1)}=o(\log p)$, and $\phi_1(D\bar s_\star )\gtrsim1$ for a sufficiently large $D$, then the following assertion holds.
	\begin{enumerate}[label={\rm (\alph*)}, resume]
		\item \label{c2thm:plmcon3} The optimal posterior contraction rates for $\theta$ in \eqref{c2eqn:thmmar} are obtained even if $\bar s_\star\rightarrow \infty$.
	\end{enumerate}
	Assume that $1<\bar\alpha<\alpha-1/2$, $(s_\star^2\log p )(\log n)^{2\alpha/(2\bar\alpha+1)}n^{(2(\bar\alpha-\alpha)+1)/(2\bar\alpha+1)}=o(1)$, $s_\star^5\log^3 p=o(n)$, and $\phi_1(Ds_\star)\gtrsim1$ for a sufficiently large $D$. Then the following assertions hold.
	\begin{enumerate}[label={\rm (\alph*)}, resume]
		\item \label{c2thm:plmbvm} The distributional approximation in \eqref{c2eqn:bvm} holds for the projection matrix $H=W_J(W_J^T W_J)^{-1}W_J^T$.
		\item  \label{c2thm:plmbvm2}  If $A_4>1$ and $s_\star \lesssim p^a$ for $a<A_4-1$, then the no-superset result in \eqref{c2eqn:sel} holds.
		\item \label{c2thm:plmbvm3}  Under the beta-min condition as well as the conditions for {\rm \ref{c2thm:plmbvm2}}, the selection consistency in \eqref{c2eqn:strsel} and the Bernstein-von Mises theorem in \eqref{c2eqn:strbvm} hold.
	\end{enumerate}
	\label{c2thm:exmplm}
\end{theorem}

Here we elaborate more on the choices of the number $J$ of basis terms.
For assertions \ref{c2thm:plmcon}--\ref{c2thm:plmcon2}, $J$ can be chosen such that $\bar\alpha=\alpha$ which gives rise to the optimal rates for the nuisance parameters. This choice, however, does not satisfy \ref{c2con:thm4asm00} and \ref{c2con:thm4asm0}, and hence we need a better approximation for $\lVert(I-H)\tilde \xi_{\eta_0}\rVert_2$ with some $\bar\alpha<\alpha$ to strictly control the remaining bias. For example, if $\bar\alpha=\alpha$, the bondedness condition for \ref{c2thm:plmcon3} is reduced to $\bar s_\star=o(1)$, which gives the optimal contraction for $\theta$ by \ref{c2thm:plmcon}.
Therefore, to incorporate the case that $\bar s_\star\rightarrow \infty$, there is a need to consider some appropriate $\bar\alpha$ that is strictly smaller than $\alpha$.
For the Bernstein-von Mises theorem and the selection consistency, the required restriction becomes even stronger such that $\bar\alpha<\alpha-1/2$.

\appendix

\section{Proofs for the main results}
In this section, we provide proofs of the main theorems.
We first describe the additional notations used for the proofs. 
For a matrix $X$, we write $\rho_1(X)\ge\rho_2(X)\ge\cdots$ for the eigenvalues of $X$ in decreasing order. The notation $\Lambda_n(\theta,\eta)=\prod_{i=1}^n ({p_{\theta,\eta,i}}/{p_{0,i}})(Y_i)$ stands for the likelihood ratio of $p_{\theta,\eta}$ and  $p_0$.
Let $\mathbb{E}_{\theta,\eta}$ denote the expectation operator with the density $p_{\theta,\eta}$ and let $\mathbb{P}_0$ denote the probability operator with the true density.
For two densities $f$ and $g$, let
$K(f,g)=\int f\log (f/g)$ and $V(f,g)=\int f|\log(f/g)-K(f,g)|^2$ stand for the Kullback-Leibler divergence and variation, respectively.
Using some constants $\underline \rho_0,\overline\rho_0>0$, we rewrite \ref{c2con:trueparcon2} as $\underline \rho_0\le\min_i\rho_{\min}(\Delta_{\eta_0,i})\le \max_i\rho_{\max}(\Delta_{\eta_0,i})\le \overline\rho_0$ for clarity.

\subsection{Proof of Theorem~\ref{c2thm:dimen}}

We first state a lemma showing that the denominator of the posterior distribution is bounded below by a factor with probability tending to one, which will be used to prove the main theorems. 
\begin{lemma}
	Suppose that {\rm\ref{c2con:maxfrob}--\ref{c2con:trueparcon2}} are satisfied.
	Then there exists a constant $K_0$ such that
	\begin{align}
	\begin{split}
	{\mathbb P}_0\bigg(&\int_{\Theta\times{\cal H}}\Lambda_n(\theta,\eta)d\Pi(\theta,\eta)\ge\pi_p(s_0)e^{-K_0(s_0\log p+n\bar\epsilon_n^2)}\bigg)\rightarrow 1.
	\label{c2eqn:lmm1}
	\end{split}
	\end{align}
	\label{c2lmm:1}
\end{lemma}

\begin{proof}
	We define the Kullback-Leibler-type neighborhood ${\cal B}_n=\{(\theta,\eta)\in\Theta\times{\cal H}: \sum_{i=1}^n K(p_{0,i},p_{\theta,\eta,i})\le C_1 n \bar\epsilon_n^2, \sum_{i=1}^n V(p_{0,i},p_{\theta,\eta,i})\le C_1 n\bar\epsilon_n^2\}$ for a sufficiently large $C_1$. Then Lemma~10 of \citet{ghosal2007convergence} implies that for any $C>0$,
	\begin{align}
	{\mathbb P}_0\left(\int_{{\cal B}_n}\Lambda_n(\theta,\eta)d\Pi(\theta,\eta)\le e^{-(1+C)C_1 n\bar\epsilon_n^2}\Pi({\cal B}_n)\right)\le \frac{1}{C^2C_1n\bar\epsilon_n^2}.
	\label{c2eqn:lgv}
	\end{align}
	Hence, it suffices to show that $\Pi({\cal B}_n)$ is bounded below as in the lemma. 
	By Lemma~\ref{c2lmm:kullback}, the Kullback-Leibler divergence and variation of the $i$th observation are given by
	\begin{align*}
	K(p_{0,i},p_{\theta,\eta,i})&=\frac{1}{2}\bigg\{-\sum_{k=1}^{m_i} \log \rho_{i,k}^\ast-\sum_{k=1}^{m_i} (1-\rho_{i,k}^\ast)\\
	&\qquad\quad+\lVert\Delta_{\eta,i}^{-1/2}
	(X_i(\theta-\theta_0)+\xi_{\eta,i}-\xi_{\eta_0,i})\rVert_2^2\bigg\},\\	V(p_{0,i},p_{\theta,\eta,i})&=
	\frac{1}{2}\sum_{k=1}^{m_i}(1-\rho_{i,k}^\ast)^2+\lVert\Delta_{\eta_0,i}^{1/2}\Delta_{\eta,i}^{-1}
	(X_i(\theta-\theta_0)+\xi_{\eta,i}-\xi_{\eta_0,i})\rVert_2^2,
	\end{align*}
	where $\rho_{i,k}^\ast,~k=1,\dots,m_i,$ are the eigenvalues of $\Delta_{\eta_0,i}^{1/2}\Delta_{\eta,i}^{-1}\Delta_{\eta_0,i}^{1/2}$.
	For ${\cal I}_{n,\delta}=\{1\le i\le n : \sum_{k=1}^{m_i}(1-\rho_{i,k}^\ast)^2\ge\delta\}$ with small $\delta>0$ and $|{\cal I}_{n,\delta}|$ the cardinality of ${\cal I}_{n,\delta}$, we see that on ${\cal B}_n$,
	\begin{align}
	\begin{split}
	a_n\bar\epsilon_n^2&\gtrsim \frac{a_n}{n}\sum_{i=1}^n \sum_{k=1}^{m_i}(1-\rho_{i,k}^\ast)^2\ge\frac{a_n\delta|{\cal I}_{n,\delta}|}{n}+\frac{a_n}{n}\sum_{i\notin{\cal I}_{n,\delta}}\sum_{k=1}^{m_i}(1-\rho_{i,k}^\ast)^2 .
	\label{c2eqn:klbound}
	\end{split}
	\end{align}
Since every $i\notin\mathcal I_{n,\delta}$ satisfies $ \sum_{k=1}^{m_i}(1-\rho_{i,k}^\ast)^2 < \delta$ for small $\delta>0$, observe that
	\begin{align*}
\sum_{i\notin{\cal I}_{n,\delta}}\sum_{k=1}^{m_i}(1-\rho_{i,k}^\ast)^2&\gtrsim	\sum_{i\notin{\cal I}_{n,\delta}}\sum_{k=1}^{m_i}(1-1/\rho_{i,k}^\ast)^2\ge \frac{1}{\overline\rho_0^2}\sum_{i\notin{\cal I}_{n,\delta}}\lVert\Delta_{\eta,i}-\Delta_{\eta_0,i}\rVert_{\rm F}^2,
	\end{align*}
	where the first inequality follows by the relation $|1-x|\asymp|1-x^{-1}|$ as $x\rightarrow 1$ and the second inequality holds by (i) of Lemma~\ref{c2lmm:2} in Appendix. Since $a_n|{\cal I}_{n,\delta}|/n\lesssim a_n\bar\epsilon_n^2$ by \eqref{c2eqn:klbound}, it follows using \eqref{c2eqn:maxfrob2} that for some constants $C_2,C_3>0$,
	\begin{align*}
	\frac{a_n}{n}\sum_{i\notin{\cal I}_{n,\delta}}\lVert\Delta_{\eta,i}-\Delta_{\eta_0,i}\rVert_{\rm F}^2&\ge a_n d_{B,n}^2(\eta,\eta_0)-\frac{ a_n|{\cal I}_{n,\delta}|}{n}\max_{1\le i\le n}\lVert\Delta_{\eta,i}-\Delta_{\eta_0,i}\rVert_{\rm F}^2\\
	&\ge (C_2-C_3 a_n\bar\epsilon_n^2) \max_{1\le i\le n}\lVert\Delta_{\eta,i}-\Delta_{\eta_0,i}\rVert_{\rm F}^2-e_n.
	\end{align*}
	Combining this with \eqref{c2eqn:klbound}, we conclude that $a_n\bar\epsilon_n^2+e_n\gtrsim\max_i\lVert\Delta_{\eta,i}-\Delta_{\eta_0,i}\rVert_{\rm F}^2$ on ${\cal B}_n$, which implies that $\max_{i,k} |1-\rho_{i,k}^\ast|$ is small for all sufficiently large $n$, by (i) of Lemma~\ref{c2lmm:2} and the inequality $|1-x|\asymp|1-x^{-1}|$ as $x\rightarrow 1$.
	Hence, $\log \rho_{i,k}^\ast$ can be expanded in the powers of $(1-\rho_{i,k}^\ast)$ to get $-\log \rho_{i,k}^\ast-(1-\rho_{i,k}^\ast)\sim(1-\rho_{i,k}^\ast)^2/2$ for every $i$ and $k$.
	Furthermore, since $\max_{i,k} |1-\rho_{i,k}^\ast|$ is sufficiently small,
	we obtain that
	$\sum_{k=1}^{m_i}(1-\rho_{i,k}^\ast)^2\lesssim\sum_{k=1}^{m_i}(1-1/\rho_{i,k}^\ast)^2\lesssim\lVert\Delta_{\eta,i}-\Delta_{\eta_0,i}\rVert_{\rm F}^2$ by (i) of Lemma~\ref{c2lmm:2}, and that  $\lVert\Delta_{\eta,i}^{-1}\rVert_{\rm sp}\lesssim\lVert\Delta_{\eta_0,i}^{1/2}\Delta_{\eta,i}^{-1}\Delta_{\eta_0,i}^{1/2}\rVert_{\rm sp}\lesssim 1$ by the restriction on the eigenvalues of $\Delta_{\eta_0,i}$.
	Combining these results, it follows that on ${\cal B}_n$, both $n^{-1}\sum_{i=1}^n K(p_{0,i},p_{\theta,\eta,i})$ and $n^{-1}\sum_{i=1}^n V(p_{0,i},p_{\theta,\eta,i})$ are bounded above by a constant multiple of $n^{-1}\lVert X(\theta-\theta_0)\rVert_2^2+d_n^2(\eta,\eta_0)$. Hence, $C_1$ can be chosen sufficiently large such that
	\begin{align}
	\begin{split}
	\Pi({\cal B}_n) & \ge \Pi\left\{(\theta,\eta)\in\Theta\times{\cal H}:n^{-1}\lVert X\rVert_\ast^2\lVert \theta-\theta_0\rVert_1^2+d_n^2(\eta,\eta_0)\le 2\bar\epsilon_n^2\right\}\\
	&\ge\Pi\left\{\theta\in\Theta:n^{-1}\lVert X\rVert_\ast^2\lVert \theta-\theta_0\rVert_1^2\le \bar\epsilon_n^2\right\} \Pi\left\{\eta\in{\cal H}:d_n^2(\eta,\eta_0)\le \bar\epsilon_n^2\right\},
	\end{split}
	\label{c2eqn:prconsp}
	\end{align}
	by the inequality $\lVert X\theta\rVert_2\le\sum_{j=1}^p\lvert\theta_j\rvert\lVert X_{\cdot j}\rVert_2\le\lVert X\rVert_\ast\lVert \theta\rVert_1$. The logarithm of the second term on the rightmost side is bounded below by a constant multiple of $-n\bar\epsilon_n^2$ by \ref{c2con:thm1prcon}. To find the lower bound for the first term, we shall first work with the case $s_0\ge1$, and then show that the same lower bound is obtained even when $s_0=0$.

	Now, assume that $s_0\ge1$ and let $\Theta_{0,n}=\{\theta_{S_0}\in\mathbb{R}^{s_0}:n^{-1/2}\lVert X\rVert_{\ast}\lVert\theta_{S_0}-\theta_{0,S_0}\rVert_1\le\epsilon\}$ for $\epsilon>0$ to be chosen later. Then
	\begin{align}
	\begin{split}
	&\Pi\{\theta\in\Theta:n^{-1/2}\lVert X\rVert_{\ast}\lVert\theta-\theta_0\rVert_1\le\epsilon\}\\
	&\quad\ge\frac{\pi_p(s_0)}{\binom p {s_0}}\int_{\Theta_{0,n}}g_{S_0}(\theta_{S_0})d\theta_{S_0}\\
	&\quad\ge\frac{\pi_p(s_0)}{\binom p {s_0}}e^{-\lambda\lVert\theta_0\rVert_1}\int_{\Theta_{0,n}}g_{S_0}(\theta_{S_0}-\theta_{0,S_0})d\theta_{S_0}
	\end{split}
	\label{c2eqn:2}
	\end{align}
	by the inequality $g_{S_0}(\theta_{S_0})\ge e^{-\lambda\lVert\theta_0\rVert_1}g_{S_0}(\theta_{S_0}-\theta_{0,S_0})$. Using the relation (6.2) of \citet{castillo2015bayesian} and the assumption on the prior in \eqref{c2eqn:lambda}, the integral on the rightmost side satisfies
	\begin{align}
	\begin{split}
	\int_{\Theta_{0,n}}g_{S_0}(\theta_{S_0}-\theta_{0,S_0})d\theta_{S_0}&\ge e^{-\lambda\epsilon\sqrt{n}/\lVert X\rVert_{\ast}}\frac{(\lambda\epsilon\sqrt{n}/\lVert X\rVert_{\ast})^{s_0}}{s_0!}\\
	&\ge e^{-L_3 \epsilon }\frac{(\epsilon\sqrt{n}/L_1p^{L_2})^{s_0}}{s_0!},
	\end{split}
	\label{c2eqn:22}
	\end{align}
	for $s_0>0$, and thus the rightmost side of \eqref{c2eqn:2} is bounded below by
	\begin{align*}
	\pi_p(s_0)(\epsilon\sqrt{n})^{s_0} \exp\left\{-\lambda\lVert\theta_0\rVert_1-L_3\epsilon -(L_1+1)s_0\log p-s_0\log L_1\right\},
	\end{align*}
	by the inequality $\binom p {s_0}s_0!\le p^{s_0}$. Choosing $\epsilon=\bar\epsilon_n$, the first term on the rightmost side of \eqref{c2eqn:prconsp} satisfies
	\begin{align*}
	&\Pi\left\{\theta\in\Theta:n^{-1}\lVert X\rVert_\ast^2\lVert \theta-\theta_0\rVert_1^2\le \bar\epsilon_n^2\right\}\\
	&\quad\ge	\pi_p(s_0)(n\bar\epsilon_n^2)^{s_0/2} \exp\left\{-\lambda\lVert\theta_0\rVert_1-L_3\bar\epsilon_n -(L_1+1)s_0\log p-s_0\log L_1\right\}.
	\end{align*}
	Note that $n\bar\epsilon_n^2>1$ and $s_0+\bar\epsilon_n +s_0\log p\lesssim s_0\log p$ if $s_0>0$, and thus the last display implies that there exists a constant $C_4>0$ such that
	\begin{align*}
	\Pi({\cal B}_n)\ge \pi_p(s_0)\exp\left\{-C_4(\lambda\lVert\theta_0\rVert_1+s_0\log p+n\bar\epsilon_n^2)\right\}.
	\end{align*}
	If $s_0=0$, the first term of~\eqref{c2eqn:prconsp} is clearly bounded below by $\pi_p(0)$, so that the same lower bound for $\Pi({\cal B}_n)$ in the last display is also obtained since we have $\lambda\lVert\theta_0\rVert_1+s_0\log p=0$. Finally, the lemma follows from \eqref{c2eqn:lgv}.
\end{proof}

\begin{proof}[Proof of Theorem~\ref{c2thm:dimen}]
	For the set ${\cal B}=\{(\theta,\eta):s_\theta > \bar s\}$ with any integer $\bar s\ge s_0$, we see that $\Pi({\cal B})$ is equal to
	\begin{align*}
	\sum_{s=\bar s+1}^p \pi_p(s)\le\pi_p(s_0)\sum_{s=\bar s+1}^p \left(\frac{A_2}{p^{A_4}}\right)^{s-s_0}\le \pi_p(s_0)\left(\frac{A_2}{p^{A_4}}\right)^{\bar s+1-s_0}\sum_{j=0}^\infty\left(\frac{A_2}{p^{A_4}}\right)^j.
	\end{align*}
	Let ${\cal E}_n$ be the event in \eqref{c2eqn:lmm1}. Since $\Lambda_n(\theta,\eta)$ is nonnegative, by Fubini's theorem and Lemma~\ref{c2lmm:1},
	\begin{align}
	\begin{split}
	{\mathbb E}_0\Pi({\cal B}\,|\,Y^{(n)})\mathbbm{1}_{{\cal E}_n}&={\mathbb E}_0\left[\frac{\int_{{\cal B}} \Lambda_n(\theta,\eta)d\Pi(\theta,\eta)}{\int \Lambda_n(\theta,\eta)d\Pi(\theta,\eta)}\mathbbm{1}_{{\cal E}_n}\right]\\
	&\le\pi_p(s_0)^{-1}\exp\{C_1 (s_0\log p+n\bar\epsilon_n^2)\}\Pi({\cal B})\\
	&\lesssim\exp\left\{(\bar s+1-s_0)(\log A_2-A_4\log p)+2C_1 s_\star \log p\right\},
	\end{split}
	\label{c2eqn:thm1prf}
	\end{align}
	for some constant $C_1$ and sufficiently large $p$. For a sufficiently large constant $C_2$,
	choose the largest integer that is smaller than $C_2 s_\star $ for $\bar s$. Replacing $\bar s+1$ by $C_2 s_\star $  in the last display, it is easy to see that the rightmost side goes to zero.
	The proof is complete since  ${\mathbb P}_0({\cal E}_n^c)\rightarrow 0$  by Lemma~\ref{c2lmm:1}.
\end{proof}

\subsection{Proof of Theorems~\ref{c2thm:helcon}--\ref{c2thm:thetacon} and Corollary~\ref{c2cor:opt}}

The following lemma shows that a small piece of the alternative centered at any $(\theta_1,\eta_1)\in\Theta\times{\cal H}$ are locally testable with exponentially small errors, provided that the center is sufficiently separated from the truth with respect to the average R\'enyi divergence. Theorem~\ref{c2thm:helcon} for posterior contraction relative to the average R\'enyi divergence will then be proved by showing that the number of those pieces is controlled by the target rate. We write $p_1$ for the density with $(\theta_1,\eta_1)$,
and $\mathbb{E}_1$ and $\mathbb{P}_1$ for the expectation and probability with $p_1$, respectively.

\begin{lemma} 
	For a given sequence $\gamma_n'>0$, a sequence $a_n$ satisfying {\rm\ref{c2con:maxfrob}}, given $(\theta_1,\eta_1)\in\Theta\times{\cal H}$ such that $R_n(p_0,p_1)\ge\delta_n^2$ with $\delta_n=o(\sqrt{\overline m})$, define
	\begin{align}
	\begin{split}
	{\cal F}_{1,n}=\bigg\{(\theta,\eta)\in\Theta\times{\cal H} \,:\,  \frac{1}{n}\sum_{i=1}^n\lVert X_i(\theta-\theta_1)+\xi_{\eta,i}-\xi_{\eta_1,i}\rVert_2^2\le \frac{\delta_n^2}{16\gamma_n'},\\  d_{B,n}(\eta,\eta_1)\le \frac{\delta_n^2}{2\overline m \gamma_n'\sqrt{a_n}}, \, \max_{1\le i\le n}\lVert\Delta_{\eta,i}^{-1}\rVert_{\rm sp}\le \gamma_n'\bigg\}.
	\label{c2eqn:smallpiece}
	\end{split}
	\end{align}
	Then under \ref{c2con:maxfrob}, there exists a test $\bar\varphi_n$ such that
	\begin{align*}
	\mathbb{E}_0 \bar\varphi_n \le e^{-n\delta_n^2},\qquad \sup_{(\theta,\eta)\in {\cal F}_{1,n} }\mathbb{E}_{\theta,\eta} (1-\bar\varphi_n)\le e^{-n\delta_n^2/16}.
	\end{align*}
	\label{c2lmm:test}
\end{lemma}

\begin{proof}
	For given $(\theta_1,\eta_1)\in\Theta\times{\cal H}$ such that $R_n(p_0,p_1)\ge\delta_n^2$, consider the most powerful test $\bar\varphi_n=\mathbbm{1}_{\{\Lambda_n(\theta_1,\eta_1)\ge 1\}}$ given by the Neyman-Pearson lemma. It is then easy to see that
	\begin{align}
	\begin{split}
	\mathbb{E}_0 \bar\varphi_n &= \mathbb{P}_0 \left(\sqrt{\Lambda_n(\theta_1,\eta_1)}\ge 1\right)\le \int \sqrt{p_0 p_1}\le e^{-n\delta_n^2},\\
	\mathbb{E}_1 (1-\bar\varphi_n)& =\mathbb{P}_1 \left(\sqrt{\Lambda_n(\theta_1,\eta_1)}\le 1\right)\le \int \sqrt{p_0 p_1}\le e^{-n\delta_n^2}.
	\end{split}	
	\label{c2eqn:consttest}
	\end{align}
	The first inequality of the lemma is a direct consequence of the first line of the preceding display. For the second inequality of the lemma, 
	note that by the Cauchy-Schwarz inequality, we have
	\begin{align*}
	\left\{\mathbb{E}_{\theta,\eta} (1-\bar\varphi_n)\right\}^2 \le\mathbb{E}_1 (1-\bar\varphi_n)\;\mathbb{E}_1 (({p_{\theta,\eta}}/{p_1})(Y^{(n)}))^2.
	\end{align*}
	Thus, by the second line of \eqref{c2eqn:consttest}, it suffices to show $\mathbb{E}_1 (({p_{\theta,\eta}}/{p_1})(Y^{(n)}))^2\le e^{7n\delta_n^2/8}$ for every $(\theta,\eta)\in {\cal F}_{1,n}$.
	Defining $\Delta_{\eta,i}^\ast=\Delta_{\eta,i}^{-1/2}\Delta_{\eta_1,i}\Delta_{\eta,i}^{-1/2}$, observe that
	\begin{align*}
	\max_{1\le i\le n}\lVert\Delta_{\eta,i}^\ast-I\rVert_{\rm sp}&\le\max_{1\le i\le n}\lVert\Delta_{\eta,i}^{-1}\rVert_{\rm sp}\lVert\Delta_{\eta,i}-\Delta_{\eta_1,i}\rVert_{\rm sp} \\
	&\le \max_{1\le i\le n}\lVert\Delta_{\eta,i}^{-1}\rVert_{\rm sp}\sqrt{a_n} d_{B,n}(\eta,\eta_1)\le \frac{\delta_n^2}{2\overline m},
	\end{align*}
	on the set ${\cal F}_{1,n}$, where the second inequality is due to \ref{c2con:maxfrob}.
	Since the leftmost side of the display is further bounded below by $\max_i|\rho_k(\Delta_{\eta,i}^\ast)-1|$ for every $k\le m_i$, we have that
	\begin{align}
	1-\frac{\delta_n^2}{2\overline m} \le \min_{1\le i\le n}\rho_{\min}(\Delta_{\eta,i}^\ast)\le \max_{1\le i\le n}\rho_{\max}(\Delta_{\eta,i}^\ast)\le 1+\frac{\delta_n^2}{2\overline m}.
	\label{c2eqn:eigrest}
	\end{align}
	Since $\delta_n^2/\overline m\rightarrow 0$ and $\rho_k(2\Delta_{\eta,i}^\ast-I)=2\rho_k(\Delta_{\eta,i}^\ast)-1$ for every $k\le m_i$, \eqref{c2eqn:eigrest} implies that $2\Delta_{\eta,i}^\ast-I$ is nonsingular for every $i\le n$, and hence on ${\cal F}_{1,n}$, it can be  shown that $\mathbb{E}_1 (({p_{\theta,\eta}}/{p_1})(Y^{(n)}))^2$ can be written as being equal to
	\begin{align}
	\begin{split}
&\prod_{i=1}^{n}\left\{\det(\Delta_{\eta,i}^\ast)^{1/2}\det(2I-{\Delta_{\eta,i}^{\ast-1}})^{-1/2}\right\}\\
	&\times \exp\Bigg\{\sum_{i=1}^n \lVert (2\Delta_{\eta,i}^\ast-I)^{-1/2} \Delta_{\eta,i}^{-1/2} (X_i(\theta-\theta_1)+\xi_{\eta,i}-\xi_{\eta_1,i})\rVert_2^2\Bigg\}.
	\label{c2eqn:expectest}
	\end{split}
	\end{align}
	To bound this, note that $	\det(\Delta_{\eta,i}^\ast)^{1/2} \det(2I-{\Delta_{\eta,i}^{\ast-1}})^{-1/2}$ is equal to
	\begin{align}
	\begin{split}
	\prod_{k=1}^{m_i}\left\{\frac{\rho_k (\Delta_{\eta,i}^\ast)}{2-\rho_k^{-1}(\Delta_{\eta,i}^\ast)}\right\}^{1/2} \le\left(\frac{1-\delta_n^4/4\overline m^2}{1-\delta_n^2/\overline m}\right)^{m_i/2} \le \left(1+\frac{3\delta_n^2}{2\overline m}\right)^{m_i/2} \le e^{3\delta_n^2/4},
	\end{split}
	\label{c2eqn:expectest1}
	\end{align}
	where the first inequality holds by \eqref{c2eqn:eigrest}, the second inequality holds by the inequality $(1-x^2)/(1-2x)\le 1+3x$ for small $x>0$, and the last inequality holds by the inequality $x+1\le e^x$. Now, for every $(\theta,\eta)\in {\cal F}_{1,n}$, observe that the exponent in \eqref{c2eqn:expectest} is bounded above by
	\begin{align*}
	\max_{1\le i \le n}\lVert (2\Delta_{\eta,i}^\ast-I)^{-1}\rVert_{\rm sp} \max_{1\le i \le n}\lVert\Delta_{\eta,i}^{-1} \rVert_{\rm sp} \sum_{i=1}^n \lVert X_i(\theta-\theta_1)+\xi_{\eta,i}-\xi_{\eta_1,i}\rVert_2^2\le \frac{n\delta_n^2}{8},
	\end{align*}
	since $\max_i\lVert (2\Delta_{\eta,i}^\ast-I)^{-1}\rVert_{\rm sp}\le 2$ for large $n$.  
	Combined with \eqref{c2eqn:expectest} and \eqref{c2eqn:expectest1}, the display completes the proof.
\end{proof}

\begin{proof}[Proof of Theorem~\ref{c2thm:helcon}]
	Let $\Theta_n=\left\{\theta\in\Theta:s_\theta\le K_1 s_\star \right\}$ and $R_n^\star(\theta,\eta)=R_n(p_{\theta,\eta},p_0)$. Then for every $\epsilon>0$,
	\begin{align}
	\begin{split}
	&{\mathbb E}_0\Pi\left((\theta,\eta)\in\Theta\times{\cal H} : \sqrt{R_n^\star(\theta,\eta)}>\epsilon \,|\,Y^{(n)}\right) \\
	&\quad\le {\mathbb E}_0\Pi\left((\theta,\eta)\in\Theta_n\times{\cal H} : \sqrt{R_n^\star(\theta,\eta)}>\epsilon \,|\,Y^{(n)}\right)+{\mathbb E}_0\Pi\left(\Theta_n^c\,|\,Y^{(n)}\right),
	\end{split}
	\label{c2eqn:cons1}
	\end{align}
	where the second term on the right hand side goes to zero by Theorem~\ref{c2thm:dimen}. Hence, it suffices to show that the first term goes to zero for $\epsilon>0$ chosen to be the threshold in the theorem.
	Now, let $\Theta_n^\ast=\{\theta\in\Theta:s_\theta\le K_1 s_\star , \lVert\theta\rVert_\infty\le p^{L_2+2}/\lVert X\rVert_\ast \}$ and define ${\cal F}_{1,n}$ as in \eqref{c2eqn:smallpiece} with $\gamma_n'=\gamma_n$ and $\delta_n=\epsilon_n$. Then Lemma~\ref{c2lmm:test} implies that small pieces of the alternative densities can be tested with exponentially small errors as long as the center is $\epsilon_n$-separated from the true parameter values relative to the average R\'enyi divergence. To complete the proof, we shall show that the minimal number $N_n^\ast$ of those small pieces that are needed to cover $\Theta_n^\ast\times{\cal H}_n$ is controlled appropriately in terms of $\epsilon_n$, and that the prior mass of $\Theta_n\setminus\Theta_n^\ast$ and ${\cal H}\setminus{\cal H}_n$ decreases fast enough to balance the denominator of the posterior distribution. (For more discussion on a construction of a test using metric entropies, see Section D.2 and Section D.3 of \citet{ghosal2017fundamentals}.)

	Note that for every $\theta,\theta'\in\Theta$ and $\eta,\eta'\in{\cal H}$,
	\begin{align*}
	\frac{1}{n}\sum_{i=1}^n\lVert X_i(\theta-\theta')+\xi_{\eta,i}-\xi_{\eta',i}\rVert_2^2\le 2\left\{ \frac{p^2}{n}\lVert X\rVert_\ast^2\lVert\theta-\theta'\rVert_\infty^2+d_{A,n}^2(\eta,\eta')\right\},
	\end{align*}
	by the inequality $\lVert X(\theta-\theta')\rVert_2\le \lVert X\rVert_\ast\lVert\theta-\theta'\rVert_1\le p\lVert X\rVert_\ast\lVert\theta-\theta'\rVert_\infty$ and the Cauchy-Schwarz inequality.
	Since $a_n<n$  and $\epsilon_n^2>n^{-1}$, 
	it is easy to see that we have ${\cal F}_{1,n}\supset {\cal F}_{1,n}'$ for
	\begin{align*}
	{\cal F}_{1,n}'=\bigg\{(\theta,\eta)\in \Theta\times{\cal H}: \, &\frac{p^2}{n}\lVert X\rVert_\ast^2\lVert\theta-\theta_1\rVert_\infty^2+d_n^2(\eta,\eta_1)\le\frac{1}{32\overline m^2\gamma_n^2n^3}, \\
	&\max_{1\le i\le n}\lVert\Delta_{\eta,i}^{-1}\rVert_{\rm sp}\le \gamma_n\bigg\},
	\end{align*}
	with the same $(\theta_1,\eta_1)$ used to define ${\cal F}_{1,n}$. Hence, $\log N_n^\ast$ is bounded above by
	\begin{align}
	\log N\left(\frac{1}{6 \overline m\gamma_n n p\lVert X\rVert_\ast},\Theta_n^\ast,\lVert\cdot\rVert_\infty\right)+\log N\left(\frac{1}{6\overline m\gamma_n n^{3/2}},{\cal H}_{n},d_n\right).
	\label{c2eqn:covnum}
	\end{align}
	Note that for any small $\delta>0$,
	\begin{align*}
	N(\delta,\Theta_n^\ast,\lVert\cdot\rVert_\infty)\le\binom{p}{\lfloor K_1s_\star \rfloor }\left(\frac{3p^{L_2+2}}{\delta\lVert X\rVert_{\ast}}\right)^{\lfloor K_1s_\star \rfloor}\le \left(\frac{3p^{L_2+3}}{\delta\lVert X\rVert_{\ast}}\right)^{ K_1s_\star },
	\end{align*}
	and thus we obtain
	\begin{align*}
	\log N\left(\frac{1}{6 \overline m\gamma_n n p\lVert X\rVert_\ast},\Theta_n^\ast,\lVert\cdot\rVert_\infty\right)&\lesssim  s_\star (\log \overline m+\log \gamma_n+\log p)\lesssim n\epsilon_n^2.
	\end{align*}
	Using the last display and the entropy condition \eqref{c2eqn:thm2c2}, the right hand side of~\eqref{c2eqn:covnum} is bounded above by a constant multiple of
	$n\epsilon_n^2$.
	Hence, by Lemma D.3 of \citet{ghosal2017fundamentals}, for every $\epsilon>\epsilon_n$,
	there exists a test $\varphi_n$ such that for some $C_1>0$,
	${\mathbb E}_0\varphi_n \le 2 \exp(C_1n \epsilon_n^2-n\epsilon^2)$
	and ${\mathbb E}_{\theta,\eta}(1-\varphi_n) \le \exp(-n\epsilon^2/16)$ for every $(\theta,\eta)\in\Theta_n^\ast\times{\cal H}_{n}$ such that $\sqrt{R_n^\star(\theta,\eta)}>\epsilon$.
	Note that under condition \eqref{c2eqn:prdim} on the prior distribution, we have $-\log\pi_p(s_0)\lesssim s_0\log p-\log \pi_p(0)\lesssim s_\star\log p$ since $\pi_p(0)$ is bounded away from zero.
	Hence, for ${\cal E}_n$ the event in~\eqref{c2eqn:lmm1} and some constant $C_2>0$, the first term on the right hand side of \eqref{c2eqn:cons1} is bounded by
	\begin{align*}
	&{\mathbb E}_0\Pi\left((\theta,\eta)\in\Theta_n\times{\cal H} : \sqrt{R_n^\star(\theta,\eta)}>\epsilon \,|\,Y^{(n)}\right)\mathbbm{1}_{{\cal E}_n}(1-\varphi_n)+{\mathbb E}_0(\varphi_n+\mathbbm{1}_{{\cal E}_n^c})\\
	&~~\le\Bigg\{\!\sup_{ (\theta,\eta)\in\Theta_n^\ast\times{\cal H}_n :R_n^\star(\theta,\eta)>\epsilon^2 }\!{\mathbb E}_{\theta,\eta}(1-\varphi_n)+\Pi(\Theta_n\!\setminus\!\Theta_n^\ast)+\Pi({\cal H}\!\setminus\!{\cal H}_{n})\Bigg\}e^{C_2 s_\star\log p}\\
	&\qquad+{\mathbb E}_0\varphi_n+{\mathbb P}_0{\cal E}_n^c,
	\end{align*}
	where the term ${\mathbb P}_0{\cal E}_n^c$ converges to zero by Lemma \ref{c2lmm:1}.
	Choosing $\epsilon=C_3\epsilon_n$ for a sufficiently large $C_3$, we have
	\begin{align*}
	{\mathbb E}_0\varphi_n\rightarrow 0,\quad\sup_{(\theta,\eta)\in\Theta_n^\ast\times{\cal H}_n:R_n^\star(\theta,\eta)>\epsilon^2}{\mathbb E}_{\theta,\eta}(1-\varphi_n)e^{C_2 s_\star\log p}\rightarrow 0.
	\end{align*}
	Furthermore, $\Pi({\cal H}\setminus{\cal H}_{n})e^{C_2 s_\star\log p}$ goes to zero by  condition \eqref{c2eqn:thm2c3}.
	Now, to show that $\Pi(\Theta_n\setminus\Theta_n^\ast)$ goes to zero exponentially fast, observe that
	\begin{align*}
	\Pi(\Theta_n\setminus\Theta_n^\ast)&=\Pi\left\{\theta\in\Theta:s_\theta\le K_1 s_\star , \lVert\theta\rVert_\infty > p^{L_2+2}/\lVert X\rVert_\ast \right\}\\
	&= \sum_{S:s\le K_1 s_\star }\frac{\pi_p(s)}{\binom{p}{s}}\int_{\{\theta_S:\lVert\theta_S\rVert_\infty > p^{L_2+2}/\lVert X\rVert_\ast\}} g_S(\theta_S)d\theta_S\\
	&\le\sum_{S:s\le K_1 s_\star }\frac{(A_2 p^{-A_4})^s}{\binom{p}{s}}\int_{\{\theta_S:\lVert\theta_S\rVert_\infty > p^{L_2+2}/\lVert X\rVert_\ast\}} g_S(\theta_S)d\theta_S.
	\end{align*}
	by the inequality  $\pi_p(s)\le (A_2 p^{-A_4})^s \pi_p(0)$ for every $S$. Since the tail probability of the Laplace distribution is given by $\int_{|x|>t}2^{-1}\lambda e^{-\lambda|x|}dx=\exp(-\lambda t)$ for every $t>0$, the rightmost side of the last display is bounded above by a constant multiple of
	\begin{align*}
\sum_{s=1}^{K_1 s_\star } s e^{-\lambda p^{L_2+2}/\lVert X\rVert_\ast} \left(\frac{A_2}{p^{A_4}}\right)^s \lesssim s_\star  e^{-\lambda p^{L_2+2}/\lVert X\rVert_\ast}.
	\end{align*}
	Since $\lambda p^{L_2+2}/\lVert X\rVert_\ast\gtrsim p^2$ by \eqref{c2eqn:lambda}, the right hand side is bounded by $e^{-C_4p^2}$ for some $C_4>0$, and thus $\Pi(\Theta_n\setminus\Theta_n^\ast)e^{C_2 s_\star\log p}$ goes to zero since $s_\star\log p=o(p^2)$. Finally, we conclude that the left hand side of \eqref{c2eqn:cons1} goes to zero  with $\epsilon=C_3\epsilon_n$.
\end{proof}

\begin{proof}[Proof of Theorem~\ref{c2thm:thetacon}]
	By Theorem~\ref{c2thm:helcon}, we obtain the contraction rate of the posterior distribution with respect to the average R\'enyi divergence $R_n(p_{\theta,\eta},p_0)$ between $p_{\theta,\eta}$ and $p_0$ given by
	\begin{align*}
	R_n(p_{\theta,\eta},p_0)=& -\frac{1}{n}\sum_{i=1}^n\log\left\{ \frac{(\det\Delta_{\eta,i})^{1/4}(\det\Delta_{\eta_0,i})^{1/4}}{\det((\Delta_{\eta,i}+\Delta_{\eta_0,i})/2)^{1/2}}\right\}\\
	&+\frac{1}{4n}\sum_{i=1}^n \lVert (\Delta_{\eta,i}+\Delta_{\eta_0,i})^{-1}(X_i(\theta-\theta_0)+\xi_{\eta,i}-\xi_{\eta_0,i}) \rVert_2^2.
	\end{align*}
	Define
	\begin{align}
	g^2(\Delta_{\eta,i},\Delta_{\eta_0,i})=1-\frac{(\det\Delta_{\eta,i})^{1/4}(\det\Delta_{\eta_0,i})^{1/4}}{\det((\Delta_{\eta,i}+\Delta_{\eta_0,i})/2)^{1/2}}.
	\label{c2eqn:defg}
	\end{align}
	Then Theorem~\ref{c2thm:helcon} implies that by the last display,
	\begin{align}
	\epsilon_n^2\gtrsim-\frac{1}{n}\sum_{i=1}^n\log(1-g^2(\Delta_{\eta,i},\Delta_{\eta_0,i}))\ge\frac{1}{n}\sum_{i=1}^n g^2(\Delta_{\eta,i},\Delta_{\eta_0,i}),
	\label{c2eqn:renrec}
	\end{align}
	where the second inequality holds by the inequality $\log x\le x-1$. 
	Note that by combining (i) and (ii) of Lemma~\ref{c2lmm:2} in Appendix, we obtain
	$g^2(\Delta_{\eta,i},\Delta_{\eta_0,i})\gtrsim\lVert \Delta_{\eta,i}-\Delta_{\eta_0,i} \rVert_{\rm F}^2$ if the left hand side is small.
	Thus, using the same approach in the proof of Lemma~\ref{c2lmm:1}, 
	\eqref{c2eqn:renrec} is further bounded below by
	\begin{align}
	\begin{split}
	&C_1 d_{B,n}^2(\eta,\eta_0)- C_2\epsilon_n^2\max_{1\le i\le n}\lVert\Delta_{\eta,i}-\Delta_{\eta_0,i}\rVert_{\rm F}^2\\
	&\quad \ge  (C_1-C_3a_n\epsilon_n^2 )d_{B,n}^2(\eta,\eta_0)- C_3 e_n\epsilon_n^2,
	\label{c2eqn:renrec2}
	\end{split}
	\end{align}
	for some constants $C_1,C_2,C_3>0$.
	Since $C_1-C_3a_n\epsilon_n^2$ is bounded away from zero and $e_n$ is decreasing,
	\eqref{c2eqn:renrec} and \eqref{c2eqn:renrec2} imply that $\epsilon_n \gtrsim d_{B,n}(\eta,\eta_0)$.
	Now, it is easy to see that by \eqref{c2eqn:maxfrob2},
	\begin{align*}
	\max_{1\le i\le n}\lVert\Delta_{\eta,i}+\Delta_{\eta_0,i}\rVert_{\rm sp}^2&\le 2\max_{1\le i\le n}\lVert\Delta_{\eta,i}-\Delta_{\eta_0,i}\rVert_{\rm sp}^2+8\max_{1\le i\le n}\lVert\Delta_{\eta_0,i}\rVert_{\rm sp}^2 \\
	&\lesssim e_n+a_nd_{B,n}^2(\eta,\eta_0)+1,
	\end{align*}
	which is bounded since $e_n+a_n \epsilon_n^2 =o(1)$. Hence, we see that for $\eta_\ast$ satisfying \ref{c2con:sep}, $n^{-1}\lVert X(\theta-\theta_0)\rVert_2^2+d_{A,n}^2(\eta,\eta_0)$ is bounded by a constant multiple of
	\begin{align*}
	&\frac{1}{n}\lVert X(\theta-\theta_0)\rVert_2^2+d_{A,n}^2(\eta,\eta_\ast)+d_{A,n}^2(\eta_\ast,\eta_0)\\
	&\quad\lesssim\frac{1}{n}\sum_{i=1}^n \lVert X_i(\theta-\theta_0)+\xi_{\eta,i}-\xi_{\eta_\ast,i} \rVert_2^2+d_{A,n}^2(\eta_\ast,\eta_0)\\
	&\quad\lesssim\frac{1}{n}\sum_{i=1}^n \lVert (\Delta_{\eta,i}+\Delta_{\eta_0,i})^{-1}(X_i(\theta-\theta_0)+\xi_{\eta,i}-\xi_{\eta_0,i}) \rVert_2^2+d_{A,n}^2(\eta_\ast,\eta_0).
	\end{align*}
	The display implies that $\lVert X(\theta-\theta_0)\rVert_2^2+n d_{A,n}^2(\eta,\eta_0)\lesssim n\epsilon_n^2$ by Theorem~\ref{c2thm:helcon} and \ref{c2con:sep}.
	Combining the results verifies the third and fourth assertions of the theorem. 
	For the remainder, observe that
	$s_{\theta-\theta_0}\le s_\theta+s_0\le K_1s_\star +s_0\lesssim s_\star $ for $\theta$  such that $s_\theta\le K_1 s_\star $. Therefore by Theorem~\ref{c2thm:dimen}, the first and the second assertions readily follow from the definitions of $\phi_1$ and $\phi_2$.
\end{proof}

\begin{proof}[Proof of Corollary~\ref{c2cor:opt}]
	We first verify the assertion~\ref{c2cor:opt-a}. If $s_0>0$ the assertion is trivial. If $s_0=0$, the condition $n\bar\epsilon_n^2/\log p\rightarrow 0$ implies that  $s_\star\rightarrow 0$, and  hence Theorem~\ref{c2thm:dimen} holds with $s_\star=0$. Since this means that $\theta=\theta_0=0$ if $s_0=0$, we can plug in $s_0$ for $s_\star$ in Theorem~\ref{c2thm:thetacon}.
	
	Similarly, the assertion~\ref{c2cor:opt-b} trivially holds if $s_0>0$ and we only need to verify the case $s_0=0$.
	By reading the proof of Theorem~\ref{c2thm:dimen}, one can see that \eqref{c2eqn:thm1prf} goes to zero for large enough $A_4$ if $s_0=0$. This completes the proof.
\end{proof}

\subsection{Proof of Theorem~\ref{c2thm:thmmar}}
To prove Theorem~\ref{c2thm:thmmar}, we first provide preliminary results. Some of these will also be used to prove Theorems~\ref{c2thm:bvm}--\ref{c2thm:sel}.

\begin{lemma}
	Suppose that {\rm\ref{c2con:maxfrob}}, {\rm\ref{c2con:thm1prcon}}, {\rm\ref{c2con:combo1}}, {\rm\ref{c2con:thm4asm00}} and {\rm\ref{c2con:thm4asm20}} are satisfied for some orthogonal projection $H$. Then, for $\Lambda_n^\ast(\theta,\eta)=(p_{\theta,\eta}/p_{\theta_0,{\tilde\eta_n(\theta,\eta)}})(Y^{(n)})$ and $\Lambda_n^\star(\theta)$ in \eqref{c2eqn:applambda} with the corresponding $H$, 
	there exists a positive sequence $\delta_n\rightarrow 0$ such that for any $\theta$ with $s_\theta\le K_1\bar s_\star$, 
	\begin{align}
	\begin{split}
	\mathbb P_0\Bigg(&\sup_{\eta\in\widetilde{\mathcal H}_n} |\log\Lambda_n^\ast(\theta,\eta)-\log\Lambda_n^\star(\theta)|\\
	&\quad\le \delta_n \left\{\lVert X(\theta-\theta_0)\rVert_2\sqrt{(s_\theta+s_0)\log p}+\lVert X(\theta-\theta_0)\rVert_2^2\right\}\Bigg)\rightarrow 1.
	\end{split}
	\label{eqn:lblbss}
	\end{align}
	\label{lmm:n1}
\end{lemma}
\begin{proof}
	If $s_\theta=s_0=0$, the left hand side in the probability operator is zero, and the assertion trivially holds. We thus only consider the case $s_\theta+s_0>0$ below.
	
	By Markov's inequality, it suffices to show that there exists a positive sequence $\delta_n'=o(\delta_n)$ such that
	\begin{align}
	\begin{split}
	&\mathbb E_0\sup_{\eta\in\widetilde{\mathcal H}_n} |\log\Lambda_n^\ast(\theta,\eta)-\log\Lambda_n^\star(\theta)|\\
	&\quad\le  \delta_n' \left\{\lVert X(\theta-\theta_0)\rVert_2\sqrt{(s_\theta+s_0)\log p}+\lVert X(\theta-\theta_0)\rVert_2^2\right\}.
	\end{split}
	\label{eqn:expbo11}
	\end{align}
	Let  $\Delta_\eta^\star\in\mathbb{R}^{n_\ast\times n_\ast}$ be the block-diagonal matrix formed by stacking $\Delta_{\eta_0,i}^{1/2}\Delta_{\eta,i}^{-1}\Delta_{\eta_0,i}^{1/2}$, $i=1,\dots,n$,
	and observe that
	\begin{align*}
	\log\Lambda_n^\ast(\theta,\eta)
	=&-\frac{1}{2}\lVert {\Delta_\eta^\star}^{1/2} (I-H)\tilde X(\theta-\theta_0)\rVert_2^2\\
	&+(\theta-\theta_0)^T \tilde X^T(I-H) \Delta_\eta^\star \{U-(\tilde\xi_\eta-\tilde\xi_{\eta_0})-H\tilde X(\theta-\theta_0)\}.
	\end{align*}
	The left hand side of \eqref{eqn:expbo11} is thus bounded by the sum of the following terms:
	\begin{align}
	\label{eqn:qq11}
	\sup_{\eta\in\widetilde{\cal H}_n}\big\lvert (\theta-\theta_0)^T\tilde X^T (I-H)(I-\Delta_\eta^\star) (I-H)\tilde X(\theta-\theta_0)\big\rvert&,\\
	\label{eqn:qq12}
	\sup_{\eta\in\widetilde{\cal H}_n}\big\lvert (\theta-\theta_0)^T\tilde X^T(I-H)\Delta_\eta^\star (\tilde\xi_\eta-\tilde\xi_{\eta_0}+H\tilde X(\theta-\theta_0))\big\rvert&,\\
	\label{eqn:qq13}
	{\mathbb E}_0\sup_{\eta\in\widetilde{\cal H}_n}\big\lvert (\theta-\theta_0)^T\tilde X^T(I-H)(I-\Delta_\eta^\star) U\big\rvert&.
	\end{align}
	
	First, observe that \eqref{eqn:qq11} is bounded above by a constant multiple of
	\begin{align}
	\begin{split}
	&\sup_{\eta\in\widetilde{\cal H}_n}\lVert I-\Delta_\eta^\star\rVert_{\rm sp}\lVert \tilde X(\theta-\theta_0)\rVert_2^2 \lesssim \lVert X(\theta-\theta_0)\rVert_2^2\sup_{\eta\in\widetilde{\cal H}_n}\max_{1\le i\le n}\lVert \Delta_{\eta,i}^{-1}-\Delta_{\eta_0,i}^{-1} \rVert_{\rm F}.
	\label{c2eqn:bv111}
	\end{split}
	\end{align}
	Using (i) of Lemma~\ref{c2lmm:2} and the inequality $|1-x|\asymp|1-x^{-1}|$ as $x\rightarrow 1$, we obtain that for $\rho_{i,k}^\ast=\rho_k(\Delta_{\eta_0,i}^{1/2}\Delta_{\eta,i}^{-1}\Delta_{\eta_0,i}^{1/2})$,
	\begin{align}
	\lVert \Delta_{\eta,i}^{-1}-\Delta_{\eta_0,i}^{-1}\rVert_{\rm F}^2&\lesssim\sum_{k=1}^{m_i}\left(1-\rho_{i,k}^\ast\right)^2
	\lesssim\sum_{k=1}^{m_i}\left(1-1/\rho_{i,k}^\ast\right)^2\lesssim\lVert \Delta_{\eta,i}-\Delta_{\eta_0,i}\rVert_{\rm F}^2,
	\label{c2eqn:invrev1}
	\end{align}
	provided that the rightmost side is sufficiently small. Because $\max_i\lVert \Delta_{\eta,i}-\Delta_{\eta_0,i} \rVert_{\rm F}^2\le e_n+a_n  d_{B,n}^2(\eta,\eta_0)\lesssim e_n+a_n \bar\epsilon_n^2$ on $\widetilde{\cal H}_n$, \eqref{c2eqn:invrev1} holds. This implies that for all sufficiently large $n$, the right hand side of \eqref{c2eqn:bv111} is bounded above by a constant multiple of
	\begin{align*}
	&\lVert X(\theta-\theta_0)\rVert_2^2\sup_{\eta\in\widetilde{\cal H}_n}\sqrt{ e_n+a_n d_{B,n}^2(\eta,\eta_0)}\lesssim \lVert X(\theta-\theta_0)\rVert_2^2\sqrt{e_n+a_n \bar\epsilon_n^2},
	\end{align*}
	where $e_n+a_n \bar\epsilon_n^2=o(1)$ due to \ref{c2con:maxfrob} and \ref{c2con:thm1prcon}.
	
	Next, \eqref{eqn:qq12} is equal to
	\begin{align*}
	\sup_{\eta\in\widetilde{\cal H}_n}\Big\lvert (\theta-\theta_0)^T\tilde X^T(I-H)\left\{(\tilde\xi_{\eta}-\tilde\xi_{\eta_0})-(I-\Delta_\eta^\star) (\tilde\xi_\eta-\tilde\xi_{\eta_0}+H\tilde X(\theta-\theta_0))\right\}\Big\rvert.
	\end{align*}
	By the triangle inequality, the display is bounded by a constant multiple of
	\begin{align}
	\begin{split}
	&\lVert X(\theta-\theta_0)\rVert_2\sup_{\eta\in\widetilde{\cal H}_n}\lVert (I-H)(\tilde\xi_{\eta}-\tilde\xi_{\eta_0})\rVert_2\\
	&+\sup_{\eta\in\widetilde{\cal H}_n}\Big\{\lVert X(\theta-\theta_0)\rVert_2^2+\lVert X(\theta-\theta_0)\rVert_2 \sqrt{n}d_{A,n}(\eta,\eta_0)\Big\}\max_{1\le i\le n}\lVert\Delta_{\eta,i}^{-1}-\Delta_{\eta_0,i}^{-1}\rVert_{\rm sp}.
	\end{split}
	\label{eqn:bvm111}
	\end{align}
	Using the same approach used in \eqref{c2eqn:invrev1}, the second term is further bounded above by a constant multiple of 
	\begin{align*}
	&\lVert X(\theta-\theta_0)\rVert_2^2\sqrt{e_n+a_n\bar\epsilon_n^2} + \lVert X(\theta-\theta_0)\rVert_2 \sqrt{n\bar\epsilon_n^2(e_n+a_n\bar\epsilon_n^2)}.
	\end{align*} 
	Therefore, by \ref{c2con:thm4asm00} and \ref{c2con:thm4asm20}, \eqref{eqn:bvm111} is bounded by
	$
	\delta_n'\{\lVert X(\theta-\theta_0)\rVert_2\sqrt{(s_0\vee 1)\log p}+\lVert X(\theta-\theta_0)\rVert_2^2\}
	$
	for some $\delta_n'\rightarrow 0$. This is not more than the right hand side of \eqref{eqn:expbo11} if $s_\theta+s_0>0$.

	Note also that \eqref{eqn:qq13} is bounded by
	\begin{align*}
	&\left\lVert \theta-\theta_0\right\rVert_1 {\mathbb E}_0\sup_{\eta\in\widetilde{\cal H}_n}\lVert \tilde X^T(I-H)(I-\Delta_\eta^\star) U\rVert_\infty\\
	&\quad\le \frac{\sqrt{s_\theta+s_0}\lVert X(\theta-\theta_0)\rVert_2}{\phi_1(s_\theta+s_0)\lVert X \rVert_\ast}{\mathbb E}_0\sup_{\eta\in\widetilde{\cal H}_n}\lVert \tilde X^T(I-H)(I-\Delta_\eta^\star) U\rVert_\infty.
	\end{align*}
	We have that $\phi_1(s_\theta+s_0)\ge \phi_1(K_1\bar s_\star+s_0)\gtrsim 1$ by condition \ref{c2con:combo1}. 
	By Lemma~\ref{lmm:supgp} below,  one can see that
	\begin{align}
	\begin{split}
	&{\mathbb E}_0\sup_{\eta\in\widetilde{\cal H}_n}\lVert \tilde X^T(I-H)(I-\Delta_\eta^\star) U\rVert_\infty\\
	&\quad\lesssim \lVert X \rVert_\ast\sqrt{\log p}\left\{ \sqrt{e_n+a_n\bar\epsilon_n^2}+ \sqrt{a_n}\int_0^{C_3\bar\epsilon_n}\sqrt{\log N(\delta,\widetilde {\mathcal H}_n,d_{B,n})} d\delta\right\},
	\label{eqn:empbo}
	\end{split}
	\end{align}
	for some $C_3>0$. The term in the braces goes to zero by \ref{c2con:thm4asm20}.
	Combining the bounds, we easily see that there exists $\delta_n'\rightarrow 0$ satisfying \eqref{eqn:expbo11}. The assertion holds by choosing $\delta_n=\sqrt{\delta_n'}$.
\end{proof}

\begin{lemma}
	Consider a neighborhood $\mathcal H_n^\ast=\{\eta\in\mathcal H: d_{B,n}(\eta,\eta_0)\le \zeta_n \}$	with any given $\zeta_n=o(a_n^{-1/2})$ for $a_n$ satisfying \ref{c2con:maxfrob}.
	Then, for any orthogonal projection $P$ and  a sufficiently large $C>0$, we have that under \ref{c2con:maxfrob},
	\begin{align*}
	&{\mathbb E}_0\sup_{\eta\in{\cal H}_n^\ast}\lVert \tilde X^T P (I-\Delta_\eta^\star) U\rVert_\infty \\
	&\quad\lesssim \lVert X\rVert_\ast\sqrt{\log p}\left\{\sqrt{e_n+a_n\zeta_n^2}+\sqrt{a_n} \int_0^{C\zeta_n}\sqrt{\log N\left(\delta,{\cal H}_n^\ast,d_{B,n}\right)}d\delta\right\},
	\end{align*}
	where $\Delta_\eta^\star\in\mathbb{R}^{n_\ast\times n_\ast}$ is the block-diagonal matrix formed by stacking the matrices $\Delta_{\eta_0,i}^{1/2}\Delta_{\eta,i}^{-1}\Delta_{\eta_0,i}^{1/2}$, $i=1,\dots,n$.
	\label{lmm:supgp}
\end{lemma}

\begin{proof}
	Let $W_{\eta,j}=\tilde X_{\cdot j}^T P(I-\Delta_\eta^\star) U$ for $\tilde X_{\cdot j}\in\mathbb{R}^{n_\ast}$ the $j$th column of $\tilde X$. Then, by Lemma 2.2.2 of \citet{van1996weak} applied with $\psi(x)=e^{x^2}-1$,
	the expectation in the lemma is equal to
	\begin{align}
	\begin{split}
	{\mathbb E}_0\max_{1\le j\le p}\sup_{\eta\in{\cal H}_n^\ast} |W_{\eta,j}|&\le\bigg\lVert\max_{1\le j\le p}\sup_{\eta\in{\cal H}_n^\ast} |W_{\eta,j}|\bigg\rVert_{\psi} \lesssim\sqrt{\log p}\max_{1\le j\le p}\bigg\lVert\sup_{\eta\in{\cal H}_n^\ast} |W_{\eta,j}|\bigg\rVert_{\psi},
	\end{split}
	\label{c2eqn:maxbou}
	\end{align}
	where $\lVert\cdot\rVert_{\psi}$ is the Orlicz norm for $\psi$.
	For any $\eta_1,\eta_2\in{\cal H}_n^\ast$, define the standard deviation pseudo-metric between $W_{\eta_1,j}$ and $W_{\eta_2,j}$ as
	\begin{align*}
	d_{\sigma,j}(\eta_1,\eta_2)
	& \coloneqq \sqrt{{\rm Var}(W_{\eta_1,j}-W_{\eta_2,j})}=\lVert (\Delta_{\eta_1}^\star-\Delta_{\eta_2}^\star)P\tilde X_{\cdot j}\rVert_2.
	\end{align*}
	Using the tail bound for normal distributions and Lemma 2.2.1 of \citet{van1996weak}, we see that $\lVert W_{\eta_1,j}-W_{\eta_2,j} \rVert_\psi\lesssim d_{\sigma,j}(\eta_1,\eta_2)$ for every $\eta_1,\eta_2\in {\cal H}_n^\ast$.
	We shall show that ${\cal H}_n^\ast$ is a separable pseudo-metric space with $d_{\sigma,j}$ for every $j\le p$.
	Then, under the true model ${\mathbb P}_0$, we see that  $\{W_{\eta,j}:\eta\in{\cal H}_n^\ast\}$ is a separable Gaussian process for $d_{\sigma,j}$.
	Hence, by Corollary 2.2.5 of \citet{van1996weak}, for any fixed $\eta'\in{\cal H}_n^\ast$, 
	\begin{align}
	\bigg\lVert\sup_{\eta\in{\cal H}_n^\ast} |W_{\eta,j}|\bigg\rVert_{\psi}\lesssim\lVert W_{\eta',j}\rVert_{\psi}+\int_0^{{\rm diam}_j({\cal H}_n^\ast)}\sqrt{\log N(\epsilon/2,{\cal H}_n^\ast,d_{\sigma,j})}d\epsilon,
	\label{c2eqn:gaumax}
	\end{align}
	where ${\rm diam}_j({\cal H}_n^\ast)=\sup\{d_{\sigma,j}(\eta_1,\eta_2):{\eta_1,\eta_2\in{\cal H}_n^\ast}\}$.
	It is clear that $W_{\eta',j}$ possesses a normal distribution with mean zero and variance $\lVert (I-\Delta_{\eta'}^\star)P\tilde X_{\cdot j}\rVert_2^2$. 
	
	Using Lemma 2.2.1 of \citet{van1996weak} again, we see that
	\begin{align}
	\begin{split}
	\lVert W_{\eta',j}\rVert_{\psi}&\lesssim\lVert (I-\Delta_{\eta'}^\star)P\tilde X_{\cdot j}\rVert_2 \\
	&\lesssim \max_{1\le i\le n}\lVert\Delta_{\eta',i}^{-1}-\Delta_{\eta_0,i}^{-1}\rVert_2\lVert X_{\cdot j}\rVert_2 \\
	&\lesssim \lVert X\rVert_\ast\sqrt{e_n+a_n\zeta_n^2},
	\end{split}
	\label{c2eqn:orbou}
	\end{align}
	for every $\eta'\in {\cal H}_n^\ast$. Here the last inequality holds by using
	\eqref{c2eqn:invrev1} and the fact that
	$\max_i\lVert \Delta_{\eta,i}-\Delta_{\eta_0,i} \rVert_{\rm F}^2\le e_n+a_n  d_{B,n}^2(\eta,\eta_0)\lesssim e_n+a_n \zeta_n^2=o(1)$ on ${\cal H}_n^\ast$, under \ref{c2con:maxfrob}.
	
	Next, to further bound the second term in \eqref{c2eqn:gaumax}, note that for every $\eta_1,\eta_2\in{\cal H}_n^\ast$,
	\begin{align*}
	a_n\zeta_n^2&\gtrsim\sum_{k=1}^2 2a_n d_{B,n}^2(\eta_k,\eta_0) \ge a_n d_{B,n}^2(\eta_1,\eta_2)\ge \max_{1\le i \le n}\lVert \Delta_{\eta_1,i}-\Delta_{\eta_2,i} \rVert_{\rm F}^2,
	\end{align*}
	which is further bounded below by
	\begin{align*}
	\min_{1\le i \le n} \rho_{\min}^2(\Delta_{\eta_2,i}) \max_{1\le i \le n} \sum_{k=1}^{m_i} \left\{1-1/\rho_k(\Delta_{\eta_2,i}^{1/2}\Delta_{\eta_1,i}^{-1}\Delta_{\eta_2,i}^{1/2})\right\}^2,
	\end{align*}
	using (i) of Lemma~\ref{c2lmm:2}.
	In the last display, we see that $\min_i \rho_{\min}(\Delta_{\eta_2,i})$ is bounded away from zero since
	\begin{align*}
	\max_{1\le i\le n}\lVert \Delta_{\eta_2,i}^{-1}\rVert_{\rm sp}&\le\max_{1\le i\le n}\lVert \Delta_{\eta_2,i}^{-1}-\Delta_{\eta_0,i}^{-1}\rVert_{\rm sp}+\max_{1\le i\le n}\lVert \Delta_{\eta_0,i}^{-1}\rVert_{\rm sp} \lesssim \sqrt{e_n+a_n\zeta_n^2}+1,
	\end{align*} 
	and hence every eigenvalue $\rho_k(\Delta_{\eta_2,i}^{1/2}\Delta_{\eta_1,i}^{-1}\Delta_{\eta_2,i}^{1/2})$  is bounded below and above by a multiple of its reciprocal, as $a_n\zeta_n^2\rightarrow0$.
	This implies that $a_n\zeta_n^2$ is further bounded below by a constant multiple of
	\begin{align*}
	&\max_{1\le i \le n} \sum_{k=1}^{m_i} \left\{1-\rho_k(\Delta_{\eta_2,i}^{1/2}\Delta_{\eta_1,i}^{-1}\Delta_{\eta_2,i}^{1/2})\right\}^2\\
	&\quad\ge\min_{1\le i \le n} \rho_{\min}^2(\Delta_{\eta_2,i})\max_{1\le i \le n}\lVert \Delta_{\eta_1,i}^{-1}-\Delta_{\eta_2,i}^{-1} \rVert_{\rm F}^2.
	\end{align*}
	By the definition of $d_{\sigma,j}$ and the preceding displays, we thus obtain
	\begin{align}
	\begin{split}
	d_{\sigma,j}(\eta_1,\eta_2)&\le\lVert \Delta_{\eta_1}^\star-\Delta_{\eta_2}^\star\rVert_{\rm sp}\lVert\tilde X_{\cdot j}\rVert_2 \\
	&\lesssim \lVert X_{\cdot j}\rVert_2\max_{1\le i\le n}\lVert \Delta_{\eta_1,i}^{-1}-\Delta_{\eta_2,i}^{-1}\rVert_{\rm sp} \\
	&\lesssim \lVert X_{\cdot j}\rVert_2 \sqrt{a_n}d_{B,n}(\eta_1,\eta_2),
	\end{split}
	\label{c2eqn:ddbound}
	\end{align}
	for every $\eta_1,\eta_2\in{\cal H}_n^\ast$.
	Hence, using that ${\rm diam}_j({\cal H}_n^\ast)\lesssim\lVert X_{\cdot j}\rVert_2\zeta_n\sqrt{a_n}$,
	we can bound the second term in \eqref{c2eqn:gaumax} above by a constant multiple of
	\begin{align*}
	\int_0^{C_1\lVert X_{\cdot j}\rVert_2\zeta_n\sqrt{a_n}}\sqrt{\log N\left({\epsilon}/{C_2\lVert X_{\cdot j}\rVert_2 \sqrt{a_n}},{\cal H}_n^\ast,d_{B,n}\right)}d\epsilon,
	\end{align*}
	for some $C_1,C_2>0$. This can be further bounded by replacing $\lVert X_{\cdot j}\rVert_2$ in the display by $\lVert X\rVert_\ast$.
	Then, using \eqref{c2eqn:maxbou}, \eqref{c2eqn:gaumax}, and \eqref{c2eqn:orbou}, and by the substitution $\delta={\epsilon}/{(C_2\lVert X\rVert_\ast \sqrt{a_n})}$ for the last display, we bound \eqref{c2eqn:maxbou} above by a constant multiple of
	\begin{align*}
	\lVert X\rVert_\ast\sqrt{\log p}\left\{\sqrt{e_n+a_n\zeta_n^2}+\sqrt{a_n} \int_0^{C_3\zeta_n}\sqrt{\log N\left(\delta,{\cal H}_n^\ast,d_{B,n}\right)}d\delta\right\},
	\end{align*}
	for some $C_3>0$.

	To complete the proof, it remains to show
	that ${\cal H}_n^\ast$ is a separable pseudo-metric space with $d_{\sigma,j}$ for every $j\le p$. By \eqref{c2eqn:ddbound}, we see that $d_{\sigma,j}(\eta_1,\eta_2)\lesssim\lVert X\rVert_\ast \sqrt{a_n}d_{B,n}(\eta_1,\eta_2)$	for every $\eta_1,\eta_2\in{\cal H}_n^\ast$. This implies that ${\cal H}_n^\ast$ is separable with $d_{\sigma,j}$ since $\mathcal H$ is separable with $d_{B,n}$.
\end{proof}

\begin{lemma}
	For any orthogonal projection $P$,
	\begin{align*}
	{\mathbb P}_0\left(\lVert\tilde X^T P U\rVert_\infty>2\underline \rho_0^{-1/2}\sqrt{\log p}\lVert X\rVert_\ast\right)&\le\frac{2}{p}.
	\end{align*}
	\label{lmm:gausian}
\end{lemma}
\begin{proof}
	Note first that $\tilde X_{\cdot j}^T P U$ has a normal distribution with mean zero and variance $\lVert P\tilde X_{\cdot j}\rVert_2^2$, and hence we have 
	\begin{align*}
	{\mathbb P}_0\left(\lVert\tilde X^T P U\rVert_\infty>t\max_{1\le j\le p}\lVert P \tilde X_{\cdot j}\rVert_2\right)\le 2pe^{-t^2/2},\quad t>0,
	\end{align*}
	by the tail probabilities of normal distributions. By choosing $t=2\sqrt{\log p}$ and using the inequality $\lVert P\tilde X_{\cdot j}\rVert_2\le\lVert\tilde X_{\cdot j}\rVert_2\le \underline \rho_0^{-1/2}\lVert X\rVert_\ast$ for every $j\le p$, we verify the assertion.
\end{proof}

\begin{lemma} 
	If {\rm\ref{c2con:combo1}} and {\rm\ref{c2con:thm4asm20}} are satisfied and $s_0\log p\lesssim n\bar\epsilon_n^2$, there exists a constant $K_0'>0$ such that
	\begin{align}
	\mathbb P_0\left(\inf_{\eta\in\widetilde{\mathcal H}_n} \int \frac{ p_{\theta,\eta} }{ p_{\theta_0,\eta} } (Y^{(n)}) d\Pi(\theta)\ge e^{-K_0'(1+s_0\log p)} \right)\rightarrow 1.
	\label{eqn:lbinfrat}
	\end{align}
	\label{lmm:n2}
\end{lemma}

\begin{proof}
	Let $\Theta_n^\ast=\{\theta\in\Theta: s_\theta = s_0 ,\lVert X(\theta-\theta_0)\rVert_2^2\le 1\}$. Restricting the integral to this set, the left hand side of the inequality in \eqref{eqn:lbinfrat} is bounded below by
	\begin{align}
	\begin{split}
	\inf_{\eta\in\widetilde{\mathcal H}_n} \int_{\Theta_n^\ast} \frac{ p_{\theta,\eta} }{ p_{\theta_0,\eta} } (Y^{(n)}) 	 d\Pi(\theta)&\ge\int_{\Theta_n^\ast}\inf_{\eta\in\widetilde{\mathcal H}_n}  \frac{ p_{\theta,\eta} }{ p_{\theta_0,\eta} } (Y^{(n)}) 	 d\Pi(\theta)\\
	&=\int_{\Theta_n^\ast}\exp\left(\inf_{\eta\in\widetilde{\mathcal H}_n} \log  \frac{ p_{\theta,\eta} }{ p_{\theta_0,\eta} } (Y^{(n)}) \right)	 d\Pi(\theta).
	\end{split}
	\label{eqn:lbinf0}
	\end{align}
	The exponent is equal to
	\begin{align}
	\begin{split}
	&\inf_{\eta\in\widetilde{\mathcal H}_n} \left\{(\theta-\theta_0)^T \tilde X^T \Delta_\eta^\star (U-\tilde\xi_\eta+\tilde\xi_{\eta_0})-\frac{1}{2} \lVert \Delta_\eta^{\star 1/2} \tilde X(\theta-\theta_0)\rVert_2^2\right\}\\
	&\quad\gtrsim -\lVert\theta-\theta_0\rVert_1 \sup_{\eta\in\widetilde{\mathcal H}_n}\lVert\tilde X^T \Delta_\eta^\star U\rVert_\infty \\
	&\qquad-\lVert X(\theta-\theta_0)\rVert_2\sup_{\eta\in\widetilde{\mathcal H}_n}\lVert\tilde\xi_\eta-\tilde\xi_{\eta_0}\rVert_2- \lVert X(\theta-\theta_0)\rVert_2^2,
	\label{eqn:lbinf}
	\end{split}
	\end{align}
	since $\lVert \Delta_\eta^\star \rVert_{\rm sp}\lesssim 1$ on $\widetilde{\mathcal H}_n$.
	We first consider the case $s_0>0$.
	Observe that
	$\sup_{\eta\in\widetilde{\mathcal H}_n}\lVert\tilde X^T \Delta_\eta^\star U\rVert_\infty \le \lVert\tilde X^T  U\rVert_\infty+\sup_{\eta\in\widetilde{\mathcal H}_n}\lVert\tilde X^T (I-\Delta_\eta^\star) U\rVert_\infty$, where the first term is bounded by a constant multiple of $\lVert X\rVert_\ast\sqrt{\log p}$ with $\mathbb P_0$-probability tending to one, due to Lemma~\ref{lmm:gausian}. By Lemma~\ref{lmm:supgp} applied with $P=I$ together with \ref{c2con:thm4asm20}, the expected value of the second term is bounded by $\delta_n\lVert X\rVert_\ast\sqrt{\log p}$ for some $\delta_n\rightarrow0$. Hence, for any $M_n\rightarrow \infty$, 
	\begin{align*}
	\mathbb P_0\left(\sup_{\eta\in\widetilde{\mathcal H}_n}\lVert\tilde X^T (I-\Delta_\eta^\star) U\rVert_\infty \le M_n\delta_n\lVert X\rVert_\ast\sqrt{\log p} \right)\rightarrow 1.
	\end{align*}
	Consequently, taking a sufficiently slowly increasing $M_n$ for the above, \eqref{eqn:lbinf} is bounded below by a constant multiple of 
	\begin{align*}
	-\lVert X\rVert_\ast\lVert\theta-\theta_0\rVert_1 \sqrt{\log p} - \lVert X(\theta-\theta_0)\rVert_2^2,
	\end{align*}
	with $\mathbb P_0$-probability tending to one. Note that $\lVert X\rVert_\ast\lVert\theta-\theta_0\rVert_1 \le \sqrt{s_\theta+s_0}\lVert X(\theta-\theta_0)\rVert_2/ \phi_1(s_\theta+s_0)$ and $\phi_1(s_\theta+s_0)= \phi_1(2s_0)\gtrsim 1$ on $\Theta_n^\ast$ by \ref{c2con:combo1}, if $s_0\log p\lesssim n\bar\epsilon_n^2$.
	The last display is thus bounded below by $-C_1s_0\log p$ for some $C_1>0$, uniformly over $\theta\in \Theta_n^\ast$.
	Consequently, with $\mathbb P_0$-probability tending to one, \eqref{eqn:lbinf0} is bounded below by
	\begin{align*}
	e^{-C_1 s_0\log p}\Pi(\Theta_n^\ast)\ge  \pi_p(s_0)e^{-C_2 s_0\log p},
	\end{align*}
	for some $C_2>0$, where the inequality holds by \eqref{c2eqn:2} and \eqref{c2eqn:22} since $\lambda\lVert \theta_0 \rVert_1\le s_0\log p$ by \ref{c2con:trueparcon1}. Since $-\log \pi_p(s_0)\lesssim s_0\log p$ if $s_0>0$, the display is further bounded below as in the assertion.
	
	If $s_0=0$, \eqref{eqn:lbinf} is equal to zero on $\Theta_n^\ast$, as this is a singleton set $\{\theta:\theta=0\}$. This means that \eqref{eqn:lbinf0} is bounded below by $\pi_p(0)$, which is also bounded away from zero. This leads to the desired assertion.
\end{proof}

\begin{proof}[Proof of Theorem~\ref{c2thm:thmmar}]
	The idea of our proof is similar in part to that of Theorem 3.5 in \citet{chae2016bayesian}.
	We only need to verify the first and fourth assertions. The second and third assertions then follow from the definitions of $\phi_1$ and $\phi_2$.
	Note also that we only need to consider the case $s_0\log p\lesssim n\bar\epsilon_n^2$, as the assertions follow from Theorems~\ref{c2thm:dimen} and \ref{c2thm:thetacon} if $s_0\log p\gtrsim n\bar\epsilon_n^2$.
	
	Let $\mathcal B_n=\{ \theta\in \Theta : s_\theta>K_4 s_0\}\cup \{ \theta\in \Theta : \lVert X(\theta-\theta_0)\rVert_2^2 > K_5s_0\log p \}$.
	Also define $\widetilde{\mathcal H}_n'$ as $\widetilde{\mathcal H}_n$ but using a constant $\tilde M_2'\le \tilde M_2$ such that $\widetilde{\mathcal H}_n'\subset \widetilde{\mathcal H}_n$.
	Then, by Theorem~\ref{c2thm:thetacon}, we have that
	\begin{align*}
	\mathbb E_0\Pi(\theta\in\mathcal B_n| Y^{(n)} )&\le \mathbb E_0\Pi(\theta\in\mathcal B_n\cap \widetilde\Theta_n,\eta\in\widetilde{\mathcal H}_n'| Y^{(n)} ) + o(1) \\
	&\le \mathbb E_0\Pi(\theta\in\mathcal B_n\cap \widetilde\Theta_n,\eta\in\widetilde{\mathcal H}_n'|Y^{(n)},\eta\in\widetilde{\mathcal H}_n ) + o(1).
	\end{align*}
	Let $\Omega$ be the event that is an intersection of the events in \eqref{eqn:lblbss}, \eqref{eqn:lbinfrat}, and the event $\{\lVert\tilde X^T (I-H) U\rVert_\infty\le2\underline \rho_0^{-1/2}\sqrt{\log p}\lVert X\rVert_\ast\}$ whose probability goes to zero by Lemma~\ref{lmm:gausian}.
	Since $\mathbb P_0(\Omega^c)\rightarrow0$, it suffices to show that 
	\begin{align}
	\begin{split}
	&\mathbb E_0\Pi(\theta\in\mathcal B_n\cap \widetilde\Theta_n,\eta\in\widetilde{\mathcal H}_n'|Y^{(n)},\eta\in\widetilde{\mathcal H}_n ) \mathbbm 1_{\Omega} \\
&\quad=\mathbb E_0\frac{\int_{\widetilde\Theta_n\cap {\mathcal B}_n} \int_{\widetilde{\mathcal H}_n'} p_{\theta,\eta} (Y^{(n)}) d\Pi(\eta) d\Pi(\theta)}{\int \int_{\widetilde{\mathcal H}_n} p_{\theta,\eta} (Y^{(n)}) d\Pi(\eta) d\Pi(\theta)}
	\end{split}
	\label{eqn:ratio}
	\end{align} 
	tends to zero.
	Observe that by Fubini's theorem, the denominator of the ratio is equal to
	\begin{align*}
	&\int_{\widetilde{\mathcal H}_n}\int \frac{ p_{\theta,\eta} }{ p_{\theta_0,\eta} } (Y^{(n)}) d\Pi(\theta) p_{\theta_0,\eta} d\Pi(\eta)\\
	 &\quad\ge \left\{\inf_{\eta\in\widetilde{\mathcal H}_n} \int \frac{ p_{\theta,\eta} }{ p_{\theta_0,\eta} } (Y^{(n)}) d\Pi(\theta) \right\} \int_{\widetilde{\mathcal H}_n} p_{\theta_0,\eta} (Y^{(n)}) d\Pi(\eta).
	\end{align*}
	By Lemma~\ref{lmm:n2}, the term in the braces on the right hand side is further bounded below by $e^{-K_0'(1+s_0\log p)}$ on the event $\Omega$.
	Note also that the numerator of the ratio in \eqref{eqn:ratio} is equal to
	\begin{align*}
	&\int_{\widetilde\Theta_n\cap {\mathcal B}_n}  \int_{\widetilde{\mathcal H}_n'}  \Lambda_n^\ast(\theta,\eta)  p_{\theta_0,\tilde\eta_n(\theta,\eta)} (Y^{(n)}) d\Pi(\eta) d\Pi(\theta)\\
	&\quad\le \left\{\int_{\widetilde\Theta_n\cap {\mathcal B}_n}  \Lambda_n^\star(\theta) \sup_{\eta\in\widetilde{\mathcal H}_n'}\frac{\Lambda_n^\ast(\theta,\eta)}{\Lambda_n^\star(\theta)}  d\Pi(\theta)\right\}  \sup_{\theta\in\widetilde\Theta_n\cap {\mathcal B}_n}\int_{\widetilde{\mathcal H}_n'}   p_{\theta_0,\tilde\eta_n(\theta,\eta)} (Y^{(n)}) d\Pi(\eta).
	\end{align*}
	Combining the bounds, on the event $\Omega$,
	the ratio in \eqref{eqn:ratio} is bounded by
	\begin{align*}
&	e^{K_0'(1+s_0\log p)}\sup_{\theta\in\widetilde\Theta_n\cap {\mathcal B}_n}\frac{\int_{\widetilde{\mathcal H}_n'}   p_{\theta_0,\tilde\eta_n(\theta,\eta)} (Y^{(n)}) d\Pi(\eta)}{\int_{\widetilde{\mathcal H}_n} p_{\theta_0,\eta} (Y^{(n)}) d\Pi(\eta)}\\
&\quad\times \int_{\widetilde\Theta_n\cap {\mathcal B}_n}   \Lambda_n^\star(\theta) \sup_{\eta\in\widetilde{\mathcal H}_n'}\frac{\Lambda_n^\ast(\theta,\eta)}{\Lambda_n^\star(\theta)}  d\Pi(\theta).
	\end{align*}
	At the end of this proof, we will verify that
	\begin{align}
	\sup_{\theta\in\widetilde\Theta_n\cap {\mathcal B}_n}\frac{\int_{\widetilde{\mathcal H}_n'}   p_{\theta_0,\tilde\eta_n(\theta,\eta)} (Y^{(n)}) d\Pi(\eta)}{\int_{\widetilde{\mathcal H}_n} p_{\theta_0,\eta} (Y^{(n)}) d\Pi(\eta)}\lesssim 1,
	\label{eqn:conq1}
	\end{align}
	with $\mathbb P_0$-probability tending to one.
	Assuming that this is true for now and letting $\Omega^\ast$ be the event satisfying \eqref{eqn:conq1}, we see that
	\eqref{eqn:ratio} is bounded by
	\begin{align*}
	e^{K_0'(1+s_0\log p)} \mathbb E_0\int_{\widetilde\Theta_n\cap {\mathcal B}_n}   \Lambda_n^\star(\theta) \sup_{\eta\in\widetilde{\mathcal H}_n'}\frac{\Lambda_n^\ast(\theta,\eta)}{\Lambda_n^\star(\theta)}  d\Pi(\theta)\mathbbm 1_{\Omega\cap \Omega^\ast} + o(1).
	\end{align*}
	To show that this tends to zero, for $\delta_n$ in Lemma~\ref{lmm:n1}, define
	$\mathcal B_{1,n}=\{ \theta\in \widetilde\Theta_n : s_\theta>K_4 s_0 , \lVert X(\theta-\theta_0)\rVert_2^2\le \delta_n^{-1/2}(s_\theta+s_0)\log p\}$,
	$\mathcal B_{2,n}=\{ \theta\in \widetilde\Theta_n : s_\theta>K_4 s_0 , \lVert X(\theta-\theta_0)\rVert_2^2 > \delta_n^{-1/2}(s_\theta+s_0)\log p\}$, and
	$\mathcal B_{3,n}=\{ \theta\in \widetilde\Theta_n : s_\theta\le K_4 s_0 , \lVert X(\theta-\theta_0)\rVert_2^2 > K_5 s_0\log p\}$ such that $\widetilde\Theta_n\cap {\mathcal B}_n = \cup_{k=1}^3 \mathcal B_{k,n}$.
	Below we will show that
	\begin{align*}
	A(\mathcal B_{k,n}) &\coloneqq  e^{K_0'(1+s_0\log p)} \\
	 &\quad\times\mathbb E_0\int_{{\mathcal B}_{k,n}}    \Lambda_n^\star(\theta) \sup_{\eta\in\mathcal H_n}\frac{\Lambda_n^\ast(\theta,\eta)}{\Lambda_n^\star(\theta)}  d\Pi(\theta)\mathbbm 1_{\Omega\cap \Omega^\ast}\rightarrow 0, \quad k=1,2,3.
	\end{align*}
	
	Since $\mathbb E_0 \Lambda_n^\star(\theta)=1$ by the moment generating function of normal distributions, we obtain that
	\begin{align*}
	A(\mathcal B_{1,n})&\le \mathbb E_0\int_{ {\mathcal B}_{1,n}} \Lambda_n^\star(\theta) e^{K_0' (1+s_0\log p)+2\delta_n^{1/2} (s_\theta+s_0)\log p} d\Pi(\theta)\\
	&\le\pi_p(0)\sum_{s>K_4 s_0}e^{K_0' (1+s_0\log p)+2\delta_n^{1/2} (s+s_0)\log p} \left(\frac{A_2}{p^{A_4}}\right)^{s-s_0}.
	\end{align*}
	If $s_0=0$, the rightmost side goes to zero for any $K_4>0$. If $s_0>0$, it still goes to zero  for $K_4$ that is much larger than $K_0'$. 
	
	Note also that by conditions \ref{c2con:trueparcon2}, \ref{c2con:combo1} and \ref{c2con:thm4asm00},  we have that for some $C_1,C_2>0$ and any $\theta$,
	\begin{align}
	\begin{split}
	\log\Lambda_n^\star(\theta)&=-\frac{1}{2}\lVert(I-H)\tilde X(\theta-\theta_0)\rVert_2^2+(\theta-\theta_0)^T \tilde X^T(I-H) U\\
	& \le -C_1\lVert X(\theta-\theta_0)\rVert_2^2 + \lVert\theta-\theta_0\rVert_1 \lVert \tilde X^T(I-H) U\rVert_\infty\\
	& \le -C_1\lVert X(\theta-\theta_0)\rVert_2^2 + C_2\lVert X(\theta-\theta_0)\rVert_2 \sqrt{(s_\theta+s_0)\log p},
	\end{split}
	\label{eqn:probbo2}
	\end{align}
	on the event $\Omega$. Hence by \eqref{eqn:lblbss} and \eqref{eqn:probbo2}, for every $\theta\in\mathcal B_{2,n}$,
	\begin{align*}
	\log\left\{\Lambda_n^\star(\theta) \sup_{\eta\in\mathcal H_n}\frac{\Lambda_n^\ast(\theta,\eta)}{\Lambda_n^\star(\theta)}\right\}
	&\le(C_2\delta_n^{1/4}+\delta_n+\delta_n^{5/4}-C_1 )\lVert X(\theta-\theta_0)\rVert_2^2\le 0,
	\end{align*}
	on the event $\Omega$. 
	Therefore,
	\begin{align*}
	A(\mathcal B_{2,n})
	&\le e^{K_0'(1+s_0\log p)}\int_{{\mathcal B}_{2,n}} d\Pi(\theta)+o(1)\\
	&\le\pi_p(0) e^{K_0'(1+s_0\log p)}\sum_{s>K_4 s_0}\left(\frac{A_2}{p^{A_4}}\right)^{s-s_0}+o(1).
	\end{align*}
	This tends to zero if $K_4$ is sufficiently large.
	
	If $s_0=0$, $\mathcal B_{3,n}$ is the empty set as it implies $\theta=\theta_0=0$.  Hence it suffices to consider the case that $s_0>0$ below.
	By \eqref{eqn:lblbss} and \eqref{eqn:probbo2} again, there exists a constant $C_3>0$ such that for every $\theta\in\mathcal B_{3,n}$,
	\begin{align*}
	&\log\left\{\Lambda_n^\star(\theta) \sup_{\eta\in\mathcal H_n}\frac{\Lambda_n^\ast(\theta,\eta)}{\Lambda_n^\star(\theta)}\right\}\\
	&\quad\le-C_1\lVert X(\theta-\theta_0)\rVert_2^2 + \left\{C_2\sqrt{\frac{K_4+1}{K_5}}+\delta_n\left(1+\frac{1}{\sqrt{K_5}}\right)\right\}\lVert X(\theta-\theta_0)\rVert_2^2 \\
	&\quad\le-C_3\lVert X(\theta-\theta_0)\rVert_2^2,
	\end{align*}
	on the event $\Omega$, where the last inequality holds by choosing $K_5$ much larger than $K_4$.
	Therefore, 
	\begin{align*}
	A(\mathcal B_{3,n})&\le e^{K_0'(1+s_0\log p)}\int_{{\mathcal B}_{3,n}}  e^{-C_3\lVert X(\theta-\theta_0)\rVert_2^2}d\Pi(\theta)\\
	&\le   e^{K_0'(1+s_0\log p)-C_3K_5s_0\log p},
	\end{align*}
	which tends to zero for $K_5$ that is much larger than $K_0'$, if $s_0>0$.
	
	It only remains to show \eqref{eqn:conq1}. Since the map $\eta\mapsto\tilde\eta_n(\theta,\eta)$ is bijective for every fixed $\theta$, for the set defined by $\tilde\eta_n(\theta,\widetilde{\cal H}_n')=\{\tilde\eta_n(\theta,\eta):\eta\in\widetilde{\cal H}_n'\}$ with given $\theta\in\widetilde\Theta_n$, we see that 
	\begin{align}
	\int_{\widetilde{\cal H}_n'}p_{\theta_0,\tilde\eta_n(\theta,\eta)}(Y^{(n)})d \Pi(\eta)&=\int_{\tilde\eta_n(\theta,\widetilde{\cal H}_n')}p_{\theta_0,\eta}(Y^{(n)})d \Pi_{n,\theta}(\eta),
	\label{c2eqn:intsub1}
	\end{align}
	by the substitution in the integral.
	Writing $\Delta_0^\ast$ the block diagonal matrix formed by stacking $\Delta_{\eta_0,i}^{1/2}$, $i=1,\dots,n$, it can be seen that
	\begin{align*}
\tilde\eta_n(\theta,\widetilde{\cal H}_n')=\bigg\{\eta\in{\cal H} : \sqrt{ \lVert\Delta_0^\ast(\tilde\xi_{\eta}-\tilde\xi_0-H\tilde X(\theta-\theta_0))\rVert_2^2 + d_{B,n}^2(\eta,\eta_0)}\le \tilde M_2' \bar\epsilon_n\bigg\}.
\end{align*}
	Hence, we see that $\tilde M_2$ can be chosen sufficiently larger than $\tilde M_2'$ such that $\tilde\eta_n(\theta,\widetilde{\cal H}_n')\subset \widetilde{\cal H}_n$ for every $\theta\in\widetilde\Theta_n$ as we have $\sqrt{n}d_{A,n}(\eta,\eta_0)\lesssim\lVert \tilde\xi_{\eta}-\tilde\xi_{\eta_0}- H\tilde X(\theta-\theta_0)\rVert_2+\lVert X(\theta-\theta_0)\rVert_2$.
	Therefore, \eqref{c2eqn:intsub1} is bounded by
	\begin{align*}
	\int_{\widetilde{\cal H}_n }p_{\theta_0,\eta}(Y^{(n)})\exp\left(\left|\log \frac{d \Pi_{n,\theta}(\eta)}{d \Pi(\eta) }\right|\right)d \Pi(\eta) \lesssim \int_{\widetilde{\cal H}_n}p_{\theta_0,\eta}(Y^{(n)})d \Pi(\eta) ,
	\end{align*}
	by \ref{c2con:thm4asm10}, since $d\Pi(\eta)=d\Pi_{n,\theta_0}(\eta)$. This verifies \eqref{eqn:conq1} and thus the proof is complete.
\end{proof}

\subsection{Proof of Theorems~\ref{c2thm:bvm}--\ref{c2thm:sel}}

To prove the shape approximation in Theorem~\ref{c2thm:bvm} and the selection results in Theorem~\ref{c2thm:sel}, we first obtain two lemmas. The first shows that the remainder of the approximation goes to zero in $\mathbb{P}_0$- probability, which is a stronger version of Lemma~\ref{lmm:n1}. The second implies that with a point mass prior  for $\theta$ at $\theta_0$, we also obtain a rate which is not worse than that in Theorem~\ref{c2thm:thetacon}.
\begin{lemma} 
	Suppose that {\rm\ref{c2con:maxfrob}}, {\rm\ref{c2con:trueparcon2}}, \rm{\ref{c2con:thm4asm0}}, and {\rm\ref{c2con:thm4asm2}} are satisfied for some orthogonal projection $H$. Then, for   $\Lambda_n^\ast(\theta,\eta)=(p_{\theta,\eta}/p_{\theta_0,{\tilde\eta_n(\theta,\eta)}})(Y^{(n)})$ and $\Lambda_n^\star(\theta)$ in \eqref{c2eqn:applambda} with the corresponding $H$, we have that
	\begin{align*}
	{\mathbb E}_0\sup_{\theta\in\widehat\Theta_n}\sup_{\eta\in\widehat{\cal H}_n}\left\lvert\log\Lambda_n^\ast(\theta,\eta)-\log\Lambda_n^\star(\theta)\right\rvert\rightarrow 0.
	\end{align*}
	\label{c2lmm:nor}
\end{lemma}

\begin{proof}
	Similar to the proof of Lemma~\ref{lmm:n1}, it suffices to show the following three assertions:
	\begin{align}
	\sup_{\theta\in\widehat\Theta_n}\sup_{\eta\in\widehat{\cal H}_n}\big\lvert (\theta-\theta_0)^T\tilde X^T (I-H)(I-\Delta_\eta^\star) (I-H)\tilde X(\theta-\theta_0)\big\rvert&\rightarrow 0,	\label{c2eqn:bvmthmpr1}\\
	\sup_{\theta\in\widehat\Theta_n}\sup_{\eta\in\widehat{\cal H}_n}\big\lvert (\theta-\theta_0)^T\tilde X^T(I-H)\Delta_\eta^\star (\tilde\xi_\eta-\tilde\xi_{\eta_0}+H\tilde X(\theta-\theta_0))\big\rvert&\rightarrow 0,	\label{c2eqn:bvmthmpr2}\\
	{\mathbb E}_0 \sup_{\theta\in\widehat\Theta_n}\sup_{\eta\in\widehat{\cal H}_n}\big\lvert (\theta-\theta_0)^T\tilde X^T(I-H)(I-\Delta_\eta^\star) U\big\rvert&\rightarrow 0.
	\label{c2eqn:bvmthmpr3}
	\end{align}
	First, note that the left side of \eqref{c2eqn:bvmthmpr1} is bounded above by a constant multiple of
	\begin{align}
	\begin{split}
	&\sup_{\theta\in\widehat\Theta_n}\sup_{\eta\in\widehat{\cal H}_n}\lVert I-\Delta_\eta^\star\rVert_{\rm sp}\lVert \tilde X(\theta-\theta_0)\rVert_2^2 \\
	&\quad\lesssim \sup_{\theta\in\widehat\Theta_n}\lVert X(\theta-\theta_0)\rVert_2^2\sup_{\eta\in\widehat{\cal H}_n}\max_{1\le i\le n}\lVert  \Delta_{\eta,i}-\Delta_{\eta_0,i}\rVert_{\rm F},
	\label{c2eqn:bv11}
	\end{split}
	\end{align}
	where the inequality holds by \eqref{c2eqn:invrev1} and the fact that $\max_i\lVert \Delta_{\eta,i}-\Delta_{\eta_0,i} \rVert_{\rm F}^2\le e_n+a_n  d_{B,n}^2(\eta,\eta_0)\lesssim e_n+a_n (s_\star\log p) /n =o(1)$ on $\widehat{\cal H}_n$.
	We see that \eqref{c2eqn:bv11} is bounded above by a constant multiple of
	\begin{align*}
	&\sup_{\theta\in\widehat\Theta_n}\lVert X\rVert_\ast\lVert\theta-\theta_0\rVert_1^2\sup_{\eta\in\widehat{\cal H}_n}\sqrt{ e_n+a_n d_{B,n}^2(\eta,\eta_0)}\lesssim s_\star^2 \log p \sqrt{e_n+\frac{a_ns_\star \log p}{n}},
	\end{align*}
	which goes to zero by \ref{c2con:thm4asm2}.
	
	Next, similar to \eqref{eqn:bvm111}, the left side of \eqref{c2eqn:bvmthmpr2} is bounded by
	\begin{align*}
	&\sup_{\theta\in\widehat\Theta_n}\lVert X(\theta-\theta_0)\rVert_2\sup_{\eta\in\widehat{\cal H}_n}\lVert (I-H)(\tilde\xi_{\eta}-\tilde\xi_{\eta_0})\rVert_2\\
	&+\sup_{\theta\in\widehat\Theta_n}\sup_{\eta\in\widehat{\cal H}_n}\Big\{\Big(\lVert X(\theta-\theta_0)\rVert_2^2+\lVert X(\theta-\theta_0)\rVert_2 \sqrt{n}d_{A,n}(\eta,\eta_0)\Big) \\
	&\qquad\qquad\qquad\times\max_{1\le i\le n}\lVert\Delta_{\eta,i}^{-1}-\Delta_{\eta_0,i}^{-1}\rVert_{\rm sp}\Big\}.
	\end{align*}
	Using the same approach used in \eqref{c2eqn:invrev1}, the display is further bounded above by a constant multiple of 
	\begin{align*}
	s_\star \sqrt{\log p}\sup_{\eta\in\widehat{\cal H}_n}\lVert (I-H)(\tilde\xi_{\eta}-\tilde\xi_{\eta_0})\rVert_2+s_\star^2\log p\sqrt{e_n+\frac{a_n s_\star\log p}{n}},
	\end{align*} which goes to zero by \ref{c2con:thm4asm0} and  \ref{c2con:thm4asm2}.
	
	Now, using Lemma~\ref{lmm:supgp}, note that \eqref{c2eqn:bvmthmpr3} is bounded above by
	\begin{align*}
	\begin{split}
	&\sup_{\theta\in\widehat\Theta_n}\left\lVert \theta-\theta_0\right\rVert_1{\mathbb E}_0\sup_{\eta\in\widehat{\cal H}_n}\lVert \tilde X^T(I-H)(I-\Delta_\eta^\star) U\rVert_\infty\\
	&\quad\lesssim s_\star \log p\Bigg\{ \sqrt{e_n+\frac{a_n s_\star \log p}{n}}\\
	&\qquad\qquad\qquad + \sqrt{a_n}\int_0^{C_1\sqrt{(s_\star \log p)/n}}\sqrt{\log N(\delta,\widehat {\mathcal H}_n,d_{B,n})} d\delta\Bigg\},
	\end{split}
	\end{align*}
	for some $C_1>0$. This tends to zero by \ref{c2con:thm4asm2}.
	
\end{proof}

\begin{lemma}
	Suppose that {\rm\ref{c2con:maxfrob}--\ref{c2con:trueparcon2}}, {\rm\ref{c2con:thm2mod1}}, and {\rm\ref{c2con:sep}} are satisfied.
	Then
	there exists a constant $K_6>0$ such that
	\begin{align*}
	\mathbb{E}_0\Pi^{\theta_0}\left(d_n(\eta,\eta_0)>K_6\bar\epsilon_n\,\big|\, Y^{(n)}\right)\rightarrow 0,
	\end{align*}
	where $\Pi^{\theta_0}(\cdot\,|\,Y^{(n)})$ is the posterior distribution induced by the point mass prior for $\theta$ at $\theta_0$, i.e., $\delta_{\theta_0}(\theta)$, in place of the prior in \eqref{c2eqn:lambda}.
	\label{c2lmm:truecon}
\end{lemma}

\begin{proof}
	Since the prior for $\theta$ is the point mass at $\theta_0$, we can reduce to a low dimensional model $Y_i^\ast\coloneqq Y_i-X_i\theta_0=\xi_{\eta,i}+\varepsilon_i$, $i=1,\dots,n$. Then the lemma can be easily verified using the main results on posterior contraction in Section~\ref{c2sec:pcr}.
	The denominator of the posterior distribution with the Dirac prior at $\theta_0$ is bounded as in Lemma~\ref{c2lmm:1}, which can be shown using \eqref{c2eqn:lgv} for the prior concentration condition \ref{c2con:thm1prcon} and the expressions for the Kullback-Leibler divergence $K(p_{0,i},p_{\theta_0,\eta,i})$ and variation $V(p_{0,i},p_{\theta_0,\eta,i})$ with the true value $\theta_0$. 
	For a local test relative to the average R\'enyi divergence, Lemma~\ref{c2lmm:test} applied with ${\cal F}_{1,n}$, modified so that it can be involved only with a given $\eta_1$ such that $R_n(p_0,p_{\theta_0,\eta_1})\ge \bar\epsilon_n^2$, implies that a small piece of the alternative is tested with exponentially small errors.
	Hence, by \ref{c2con:thm2mod1}, we obtain the  contraction rate $\bar\epsilon_n^2$ relative to $R_n(p_0,p_{\theta_0,\eta})$ for $\Pi^{\theta_0}(\cdot\,|\,Y^{(n)})$, as in the proof of Theorem~\ref{c2thm:helcon}. The lemma is then obtained by recovering the contraction rate  of $\eta$ with respect to $d_n$ using the approach in the proof of Theorem~\ref{c2thm:thetacon}. 
\end{proof}

\begin{proof}[Proof of Theorem~\ref{c2thm:bvm}]
	Our proof is based on the proof of Theorem~6 in \citet{castillo2015bayesian}, but is more involved due to $\eta$.
	We use the fact that for any probability measure $Q$ and its renormalized restriction $Q_{\cal A}(\cdot)=Q(\cdot\cap {\cal A})/Q({\cal A})$ to a set ${\cal A}$, we have $\lVert Q-Q_{\cal A} \rVert_{\rm TV}\le 2Q({\cal A}^c)$.
	First, using a sufficiently large constant $\hat M_2'$ that is smaller than $\hat M_2$, define $\widehat{\cal H}_n'$ as $\widehat{\cal H}_n$ in \eqref{c2eqn:subsets} such that $\widehat{\cal H}_n'\subset\widehat{\cal H}_n$.
	Let $\widetilde\Pi((\theta,\eta)\in\cdot)$ be the prior distribution restricted and renormalized on $\widehat\Theta_n\times\widehat{\cal H}_n'$ and $\widetilde\Pi((\theta,\eta)\in\cdot\,|\,Y^{(n)})$ be the corresponding posterior distribution. Also, $\widetilde\Pi^\infty(\theta\in\cdot\,|\,Y^{(n)})$ is the restricted and renormalized version of $\Pi^\infty(\theta\in\cdot\,|\,Y^{(n)})$ to the set $\widehat\Theta_n$.
	Then the left hand side of the theorem is bounded above by
	\begin{align}
	\begin{split}
	&\left\lVert\Pi(\theta\in\cdot\,|\,Y^{(n)})-\widetilde\Pi(\theta\in\cdot\,|\,Y^{(n)})\right\rVert_{\rm TV}+\left\lVert\widetilde\Pi(\theta\in\cdot\,|\,Y^{(n)})-\widetilde\Pi^\infty(\theta\in\cdot\,|\,Y^{(n)})\right\rVert_{\rm TV}\\
	&+\left\lVert\Pi^\infty(\theta\in\cdot\,|\,Y^{(n)})-\widetilde\Pi^\infty(\theta\in\cdot\,|\,Y^{(n)})\right\rVert_{\rm TV},
	\label{c2eqn:thm5tv}
	\end{split}
	\end{align}
	where the first summand goes to zero in ${\mathbb P}_0$-probability since
	$\Pi((\theta,\eta)\in\widehat\Theta_n\times\widehat{\cal H}_n'\,|\,Y^{(n)})\rightarrow 1$ in ${\mathbb P}_0$-probability by Theorem \ref{c2thm:dimen} and Theorem \ref{c2thm:thetacon}.
	
	To show that the second summand goes to zero in ${\mathbb P}_0$-probability, note that for every measurable ${\cal B}\subset\mathbb{R}^p$, we obtain
	\begin{align*}
	\widetilde\Pi(\theta\in {\cal B}\,|\,Y^{(n)})&\propto\int_{{\cal B}\cap\widehat\Theta_n}\int_{\widehat{\cal H}_n'}{p_{\theta,\eta}}(Y^{(n)})  \: e^{-\lambda\lVert\theta\rVert_1} d\Pi(\eta) dV(\theta)\\
	&=\int_{{\cal B}\cap\widehat\Theta_n}\int_{\widehat{\cal H}_n'}\Lambda_n^\ast(\theta,\eta) \: e^{-\lambda\lVert\theta\rVert_1} p_{\theta_0,\tilde\eta_n(\theta,\eta)}(Y^{(n)})\,d \Pi(\eta)dV(\theta),
	\\
	\widetilde\Pi^\infty(\theta\in {\cal B}\,|\,Y^{(n)})&\propto\int_{{\cal B}\cap\widehat\Theta_n}\Lambda_n^\star(\theta)  dV(\theta)\\
	&\propto\int_{{\cal B}\cap\widehat\Theta_n}\Lambda_n^\star(\theta) \: e^{-\lambda\lVert\theta_0\rVert_1} \int_{{\cal H}}p_{\theta_0,\eta}(Y^{(n)})d\Pi(\eta)dV(\theta),
	\end{align*} 
	where $dV(\theta)=\sum_{S: s\le K_1 s_\star }{\pi_p(s)}{\binom p {s}}^{-1}(\lambda/2)^{s} d\{\mathcal L(\theta_S)\otimes\delta_0(\theta_{S^c})\}$.
	In the last line, the factor $e^{-\lambda\lVert\theta_0\rVert_1} \int_{{\cal H}}p_{\theta_0,\eta}(Y^{(n)})d\Pi(\eta)$ cancels out in the normalizing constant, but is inserted for the sake of comparison.
	For any sequences of measures $\{\mu_S\}$ and $\{\nu_S\}$, if $\nu_S$ is absolutely continuous with respect to $\mu_S$ with the Radon-Nikodym derivative $d{\nu_S}/d{\mu_S}$, then it can be easily verified that
	\begin{align*}
	\left\lVert\frac{\sum_S\mu_S}{\lVert\sum_S\mu_S\rVert_{\rm TV}}-\frac{\sum_S\nu_S}{\lVert\sum_S\nu_S\rVert_{\rm TV}}\right\rVert_{\rm TV}&\le \frac{2\sum_S\lVert\mu_S-\nu_S\rVert_{\rm TV}}{\lVert\sum_S\mu_S\rVert_{\rm TV}}
	\le 2\sup_S\left\lVert 1-\frac{d{\nu_S}}{d{\mu_S}}\right\rVert_\infty.
	\end{align*}
	Hence, for $C_n=\int_{{\cal H}}p_{\theta_0,\eta}(Y^{(n)})d\Pi(\eta)$, we see that the second summand of \eqref{c2eqn:thm5tv} is bounded by
	\begin{align*}
	2 \sup_{\theta\in\widehat\Theta_n}\left\lvert 1- \frac{1}{C_n}\int_{\widehat{\cal H}_n'}\frac{ \Lambda_n^\ast(\theta,\eta)e^{-\lambda\lVert\theta\rVert_1}}{\Lambda_n^\star(\theta)e^{-\lambda\lVert\theta_0\rVert_1}}p_{\theta_0,\tilde\eta_n(\theta,\eta)}(Y^{(n)})d\Pi(\eta)\right\rvert.
	\end{align*}
	Using the fact that $|\lambda(\lVert\theta\rVert_1-\lVert\theta_0\rVert_1)|\le\lambda\lVert\theta-\theta_0\rVert_1\lesssim \lambda s_\star \sqrt{\log p}/\lVert X\rVert_\ast\rightarrow 0$ on $\widehat\Theta_n$ and that $\sup\{|1-\Lambda_n^\ast(\theta,\eta)/\Lambda_n^\star(\theta)|:\theta\in\widehat\Theta_n,\eta\in\widehat{\cal H}_n'\}$ goes to zero in ${\mathbb P}_0$-probability by Lemma~\ref{c2lmm:nor}, the last display is further bounded by
	\begin{align}
	2\sup_{\theta\in\widehat\Theta_n}\left\lvert 1- \left\{1+o(1)+o_{{\mathbb P}_0}(1)\right\}\frac{1}{C_n}\int_{\widehat{\cal H}_n'}p_{\theta_0,\tilde\eta_n(\theta,\eta)}(Y^{(n)})d\Pi(\eta)\right\rvert.
	\label{c2eqn:thm5scub2}
	\end{align}
	Now, note that the map $\eta\mapsto\tilde\eta_n(\theta,\eta)$ is bijective for every fixed $\theta\in\widehat\Theta_n$. Thus for the set defined by $\tilde\eta_n(\theta,\widehat{\cal H}_n')=\{\tilde\eta_n(\theta,\eta):\eta\in\widehat{\cal H}_n'\}$ with given $\theta\in\widehat\Theta_n$, we see that 
	\begin{align}
	\int_{\widehat{\cal H}_n'}p_{\theta_0,\tilde\eta_n(\theta,\eta)}(Y^{(n)})d \Pi(\eta)&=\int_{\tilde\eta_n(\theta,\widehat{\cal H}_n')}p_{\theta_0,\eta}(Y^{(n)})d \Pi_{n,\theta}(\eta),
	\label{c2eqn:intsub}
	\end{align}
	by the substitution in the integral. Similar to the proof of Theorem~\ref{c2thm:thmmar}, observe that 
	\begin{align*}
	\tilde\eta_n(\theta,\widehat{\cal H}_n')=\Big\{\eta\in{\cal H} : & \, \lVert\Delta_0^\ast(\tilde\xi_{\eta}-\tilde\xi_0-H\tilde X(\theta-\theta_0))\rVert_2\le \hat M_2' s_\star \sqrt{(\log p)/n}, \\
	& ~ d_{B,n}(\eta,\eta_0)\le \hat M_2'\sqrt{(s_\star\log p)/n} \Big\}.
		\end{align*}
	Hence, we see that $\hat M_2$ can be chosen sufficiently large such that $\tilde\eta_n(\theta,\widehat{\cal H}_n')\subset \widehat{\cal H}_n$ for every $\theta\in\widehat\Theta_n$ as we have $\sqrt{n}d_{A,n}(\eta,\eta_0)\lesssim\lVert \tilde\xi_{\eta}-\tilde\xi_{\eta_0}- H\tilde X(\theta-\theta_0)\rVert_2+\lVert X\rVert_\ast\lVert\theta-\theta_0\rVert_1$.
	Therefore, since $d\Pi(\eta)=d\Pi_{n,\theta_0}(\eta)$, one can see that \eqref{c2eqn:intsub} is written as
	\begin{align*}\{1+o(1)\}\int_{\tilde\eta_n(\theta,\widehat{\cal H}_n')}p_{\theta_0,\eta}(Y^{(n)})d \Pi(\eta),
	\end{align*}
	by \ref{c2con:thm4asm1}, and hence \eqref{c2eqn:thm5scub2} is equal to
	\begin{align}
	2\sup_{\theta\in\widehat\Theta_n}\left\lvert 1- \left\{1+o_{{\mathbb P}_0}(1)\right\}\frac{\int_{\tilde\eta_n(\theta,\widehat{\cal H}_n')}p_{\theta_0,\eta}(Y^{(n)})d \Pi(\eta)}{\int_{{\cal H}}p_{\theta_0,\eta}(Y^{(n)})d\Pi(\eta)}\right\rvert.
	\label{c2eqn:thm5scub3}
	\end{align}
	Now, observe that we also have the inequality of the other direction: $\lVert \tilde \xi_{\eta}-\tilde \xi_{\eta_0}-H \tilde X(\theta-\theta_0)\rVert_2\lesssim \sqrt{n}d_{A,n}(\eta,\eta_0)+\lVert X\rVert_\ast\lVert\theta-\theta_0\rVert_1$. This means that $\hat M_2'$ can be chosen sufficiently large such that $\{\eta\in {\cal H} : d_n(\eta,\eta_0)\le K_6 \bar\epsilon_n\}\subset\tilde\eta_n(\theta,\widehat{\cal H}_n')$ for every $\theta\in\widehat\Theta_n$. 
	Hence, with appropriately chosen constants, we obtain
	\begin{align*}
	\inf_{\theta\in\widehat\Theta_n}\frac{\int_{\tilde\eta_n(\theta,\widehat{\cal H}_n')}p_{\theta_0,\eta}(Y^{(n)})d \Pi(\eta)}{\int_{{\cal H}}p_{\theta_0,\eta}(Y^{(n)})d \Pi(\eta)}&=\inf_{\theta\in\widehat\Theta_n}\Pi^{\theta_0}\left(\eta\in\tilde\eta_n(\theta,\widehat{\cal H}_n') \,|\, Y^{(n)}\right)\\
	&\ge\Pi^{\theta_0}\left(d_n(\eta,\eta_0)\le K_6\bar\epsilon_n \,\big|\, Y^{(n)}\right).
	\end{align*}
	The rightmost term goes to one with probability tending to one by Lemma~\ref{c2lmm:truecon}. This implies that \eqref{c2eqn:thm5scub3} goes to zero in $\mathbb{P}_0$-probability, completing the proof for the second part of \eqref{c2eqn:thm5tv}.

	Next, we show that $\Pi^\infty(\theta\in\widehat\Theta_n\,|\,Y^{(n)})$ goes to one in ${\mathbb P}_0$-probability to verify that the last summand in \eqref{c2eqn:thm5tv} goes to zero in ${\mathbb P}_0$-probability. Observe  that $\Pi^\infty(\theta\in\widehat\Theta_n^c\,|\,Y^{(n)})$ is equal to
	\begin{align}
	\frac{\int_{\widehat\Theta_n^c} \exp\left\{{-\frac{1}{2}\lVert (I-H)\tilde X(\theta-\theta_{0})\rVert_2^2+U^T(I-H)\tilde X(\theta-\theta_{0})}\right\}dV(\theta)}{\int_{\mathbb{R}^p} \exp\left\{{-\frac{1}{2}\lVert (I-H)\tilde X(\theta-\theta_{0})\rVert_2^2+U^T(I-H)\tilde X(\theta-\theta_{0})}\right\} dV(\theta)}.
	\label{c2eqn:out1}
	\end{align}
	Clearly, the denominator is bounded below by
	\begin{align}
	\begin{split}
	\frac{\pi_p(s_0)}{\binom p {s_0}}\left(\frac{\lambda}{2}\right)^{s_0}\int_{\mathbb{R}^{s_0}} \exp\bigg\{&-\frac{1}{2}\lVert (I-H)\tilde X_{S_0}(\theta_{S_0}-\theta_{0,{S_0}})\rVert_2^2
	\\
	&+U^T(I-H)\tilde X_{S_0}(\theta_{S_0}-\theta_{0,{S_0}})\bigg\}d\theta_{S_0}.
	\end{split}
	\label{c2eqn:thm5den}
	\end{align}
	Since the measure $Q$ defined by $Q(d\theta_{S_0})=\exp\{-(1/2)\lVert(I-H)\tilde X_{S_0}(\theta_{S_0}-\theta_{0,S_0})\rVert_2^2\}$
	is symmetric about $\theta_{0,S_0}$, the mean of $(\theta_{S_0}-\theta_{0,S_0})$ with respect to the normalized probability measure $\widetilde Q=Q/Q(\mathbb{R}^{s_0})$ is zero.
	Note also that $\Gamma_S=\tilde X_{S}^T(I-H)\tilde X_{S}$ is nonsingular for every $S$ such that $s\le K_1s_\star $ by \ref{c2con:thm4asm0}. Thus, by Jensen's inequality, \eqref{c2eqn:thm5den} is bounded below by
	\begin{align*}
	&\frac{\pi_p(s_0)}{\binom p {s_0}}\left(\frac{\lambda}{2}\right)^{s_0}\int_{\mathbb{R}^{s_0}} \exp\left\{{-\frac{1}{2}\lVert (I-H)\tilde X_{S_0}(\theta_{S_0}-\theta_{0,{S_0}})\rVert_2^2}\right\}d\theta_{S_0} \\
	&\quad=\frac{\pi_p(s_0)}{\binom p {s_0}}\left(\frac{\lambda}{2}\right)^{s_0}\frac{(2\pi)^{s_0/2}}{\det(\Gamma_{S_0})^{1/2}}.
	\end{align*}
	Applying the arithmetic-geometric mean inequality to the eigenvalues, we obtain $\det(\Gamma_{S_0})\le ({\rm tr}(\Gamma_{S_0})/s_0)^{s_0}\le\lVert(I-H)\tilde X_{S_0}\rVert_\ast^{2s_0}\le\underline \rho_0^{-s_0}\lVert X\rVert_\ast^{2s_0}$, and hence $\det(\Gamma_{S_0})^{1/2}/\lambda^{s_0}\le \underline \rho_0^{-s_0/2}(L_1 p^{L_2})^{s_0}$ by \eqref{c2eqn:lambda}. Furthermore, we have $\pi_p(s_0)\gtrsim A_1^{s_0}p^{-A_3s_0}$ by \eqref{c2eqn:prdim} and $\binom{p}{s_0}\le p^{s_0}$. Hence, the preceding display is further bounded below by a constant multiple of
	\begin{align}
	p^{-(1+L_2+A_3)s_0}\left(\frac{A_1\sqrt{\underline \rho_0\pi}}{L_1\sqrt{2}}\right)^{s_0}.
	\label{c2eqn:thm4lb2}
	\end{align}
	To bound the numerator of \eqref{c2eqn:out1}, let $D_n=2\underline \rho_0^{-1/2}\sqrt{\log p}\lVert X\rVert_\ast$ and ${\cal U}_n=\{\lVert\tilde X^T (I-H)U\rVert_\infty\le D_n\}$. Then
	it suffices to show that \eqref{c2eqn:out1} goes to zero in ${\mathbb P}_0$-probability on the set ${\cal U}_n$ as $\mathbb P_0({\cal U}_n^c)\rightarrow 0$ by Lemma~\ref{lmm:gausian}.
	Note that on the set ${\cal U}_n$ we have
	\begin{align*}
	U^T(I-H)\tilde X(\theta-\theta_{0}) &\le
	D_n\lVert\theta-\theta_0\rVert_1 \\
	&\le D_n\frac{2\sqrt{\overline\rho_0}\lVert \tilde X(\theta-\theta_0)\rVert_2 |S_{\theta-\theta_0}|^{1/2}}{\lVert X\rVert_\ast \phi_1(|S_{\theta-\theta_0}|)}-D_n\lVert\theta-\theta_0\rVert_1.
	\end{align*}
	Using that $\lVert u\rVert_2\lesssim\lVert(I-H)u\rVert_2$ for every $u\in{\rm span}(\tilde X_S)$ with $s\le K_1 s_\star $  by \ref{c2con:thm4asm0}, the preceding display is, for some constant $C_1>0$, further bounded above by
	\begin{align*}
	&D_n \frac{2\sqrt{\overline\rho_0}C_1\lVert (I-H)\tilde X(\theta-\theta_0)\rVert_2 |S_{\theta-\theta_0}|^{1/2}}{\lVert X\rVert_\ast \phi_1(|S_{\theta-\theta_0}|)}-D_n\lVert\theta-\theta_0\rVert_1\\
	&\quad\le \frac{1}{2}\lVert (I-H)\tilde X(\theta-\theta_0)\rVert_2^2+\frac{2\overline\rho_0 C_1^2D_n^2|S_{\theta-\theta_0}|}{\lVert X\rVert_\ast^2 \phi_1(|S_{\theta-\theta_0}|)^2}-D_n\lVert\theta-\theta_0\rVert_1,
	\end{align*}
	by the Cauchy-Schwarz inequality. We have $s_{\theta-\theta_0}\le K_1s_\star +s_0$ on the support of the measure $V$. Hence, on the event ${\cal U}_n$, the numerator of \eqref{c2eqn:out1} is bounded above by
	\begin{align*}
	&\exp\left\{\frac{2\overline\rho_0C_1^2D_n^2(K_1s_\star +s_0)}{ \lVert X\rVert_\ast^2\phi_1(K_1s_\star +s_0)^2}-\frac{\hat M_1 D_ns_\star \sqrt{ \log p}}{2\lVert X\rVert_\ast}\right\} \\
	&\times\sum_{S: s\le K_1 s_\star }\frac{\pi_p(s)}{\binom p {s}}\int\left(\frac{\lambda}{2}\right)^{s}e^{-({D_n}/{2})\lVert\theta_S-\theta_{0,S}\rVert_1}d\theta_S \\
	&\quad\le \exp\left\{\frac{8\overline\rho_0 C_1^2(K_1+1)s_\star \log p}{ \underline \rho_0\phi_1(K_1s_\star +s_0)^2}-\frac{\hat M_1 s_\star  \log p}{\sqrt{\underline \rho_0}}\right\}\sum_{s=0}^p\pi_p(s) \left(L_3\sqrt\frac{{\underline \rho_0}}{n}\right)^s,
	\end{align*}
	since $D_n/2\ge \lambda\sqrt{n}/(L_3\sqrt{\underline \rho_0})$.
	Note that we have 
	\begin{align*}
	\sum_{s=0}^p\pi_p(s) \left(L_3\sqrt\frac{{\underline \rho_0}}{n}\right)^s\lesssim\sum_{s=0}^p \left(\frac{A_2 L_3 }{p^{A_4}} \sqrt\frac{{\underline \rho_0}}{n} \right)^s\lesssim 1,
	\end{align*}
	by \eqref{c2eqn:prdim} and that $\phi_1(K_1s_\star +s_0)$ in the denominators is bounded away from zero by the assumption.
	Thus, the last display combined with \eqref{c2eqn:thm4lb2} shows that \eqref{c2eqn:out1} goes to zero on the event ${\cal U}_n$, provided that $\hat M_1$ is chosen sufficiently large.
	
	Finally we conclude that \eqref{c2eqn:thm5tv} goes to zero in $\mathbb{P}_0$-probability. Since the total variation metric is bounded by 2, the convergence in mean holds as in the assertion.
\end{proof}

\begin{proof}[Proof of Theorem~\ref{c2thm:sel}]
	Our proof follows the proof of Theorem~4 in \citet{castillo2015bayesian}.
	Since ${\mathbb E}_0\lVert\Pi(\theta\in\cdot\,|\,Y^{(n)})-\Pi^\infty(\theta\in\cdot\,|\,Y^{(n)})\rVert_{\rm TV}$ tends to zero by Theorem~\ref{c2thm:bvm}, it suffices to show that ${\mathbb E}_0\Pi^\infty(\theta:S_\theta\in{\cal S}_n \,|\,Y^{(n)})\rightarrow 0$ for ${\cal S}_n=\{S:s\le K_1s_\star , S\supset S_0, S\ne S_0\}$. For the orthogonal projection defined by $\tilde H_S=(I-H)\tilde X_S\Gamma_S^{-1}\tilde X_S^T(I-H)$ with $\Gamma_S=\tilde X_S^T(I-H)\tilde X_S$, we see that $\Pi^\infty(\theta:S_\theta\in {\cal S}_n\,|\,Y^{(n)})$ is bounded by
	\begin{align*}
	\sum_{s=s_0+1}^{K_1s_\star }\frac{\pi_p(s)\binom{p}{s_0}\binom{p-s_0}{s-s_0}}{\pi_p(s_0)\binom{p}{s}}\left(\frac{\lambda\sqrt{\pi}}{\sqrt{2}}\right)^{s-s_0}\max_{S\in{\cal S}_n:|S|=s}\left[\frac{\det(\Gamma_{S_0})^{1/2}}{\det(\Gamma_{S})^{1/2}}e^{\lVert (\tilde H_S-\tilde H_{S_0}) U\rVert_2^2/2}\right],
	\end{align*}
	by \eqref{c2eqn:norapp},
	since $(\tilde H_S-\tilde H_{S_0})\tilde X \theta_0=(\tilde H_S-\tilde H_{S_0})(I-H)\tilde X_{S_0}\theta_{0,S_0}=0$ for every $S\in{\cal S}_n$ due to $S_0\subset S$ on ${\cal S}_n$.
	Note that $\rho_k(\Gamma_{S_0})\le\rho_k(\Gamma_S)$ for $k=1,\dots,s_0$, because $\Gamma_{S_0}$ is a principal submatrix of $\Gamma_S$. Hence, $\det(\Gamma_{S_0})$ is equal to
	\begin{align}
	\begin{split}
\prod_{k=1}^{s_0}\rho_k(\Gamma_{S_0})\le\prod_{k=1}^{s_0}\rho_k(\Gamma_{S})\le\frac{\det(\Gamma_S)}{ \rho_{\min}(\Gamma_S)^{s-s_0}}\le\frac{\det(\Gamma_S)}{(C_1\overline\rho_0^{-1/2}\phi_2(s)\lVert X\rVert_\ast)^{2(s-s_0)}},
	\end{split}
	\label{c2eqn:eigub}
	\end{align}
	for some $C_1>0$.
	The last inequality holds since by \ref{c2con:thm4asm0}, there exists a constant $C_1>0$ such that $C_1^2\lVert v\rVert_2^2\le \lVert (I-H) v\rVert_2^2$ for every $v\in{\rm span}(\tilde X_S)$ with $s\le K_1 s_\star $, and hence we have that by the definition of $\phi_2$,
	\begin{align*}
	\rho_{\min}(\Gamma_S)=\inf_{u\in\mathbb{R}^s,u\ne0}\frac{\lVert(I-H)\tilde X_S u\rVert_2^2}{\lVert u\rVert_2^2}\ge \frac{C_1^2\phi_2(s)^2\lVert X\rVert_\ast^2}{\overline\rho_0}.
	\end{align*}
	Now, we shall show that for any fixed $b>2$,
	\begin{align}
	{\mathbb P}_0\left(\lVert (\tilde H_S-\tilde H_{S_0}) U\rVert_2^2 \le b(s-s_0)\log p ,\text{~for every $S\in{\cal S}_n$}\right)\rightarrow 1.
	\label{c2eqn:chiub}
	\end{align}
	Note that $\lVert (\tilde H_S-\tilde H_{S_0}) U\rVert_2^2$ has a chi-squared distribution with degree of freedom $s-s_0$. Therefore, by Lemma~5 of \citet{castillo2015bayesian}, there exists a constant $C_2$ such that for every $b>2$ and given $s\ge s_0+1$,
	\begin{align*}
	{\mathbb P}_0\left(\max_{S\in{\cal S}_n:|S|=s}\lVert (\tilde H_S-\tilde H_{S_0}) U\rVert_2^2> b\log N_s \right)
	&\le  \left(\frac{1}{N_s}\right)^{(b-2)/4} e^{C_2(s-s_0)},
	\end{align*}
	where $N_s=\binom{p-s_0}{s-s_0}$ is the cardinality of the set $\{S\in{\cal S}_n:|S|=s\}$. Since $N_s\le (p-s_0)^{s-s_0}\le p^{s-s_0}$, for ${\cal T}_n$ the event in the relation \eqref{c2eqn:chiub}, it follows that
	\begin{align*}
	{\mathbb P}_0({\cal T}_n^c) \le \sum_{s=s_0+1}^{K_1s_\star } \left(\frac{1}{N_s}\right)^{(b-2)/4} e^{C_2(s-s_0)}.
	\end{align*}
	This goes to zero as $p\rightarrow\infty$, since for $s\le K_1s_\star $,
	\begin{align*}
	N_s\ge \frac{(p-s)^{s-s_0}}{(s-s_0)!}\ge \frac{(p-K_1s_\star )^{s-s_0}}{(s-s_0)^{s-s_0}}\ge \left(\frac{p-K_1s_\star }{K_1s_\star }\right)^{s-s_0},
	\end{align*}
	and $s_\star /p=o(1)$. To complete the proof, it remains to show that $\Pi^\infty(\theta:S_\theta\in {\cal S}_n\,|\,Y^{(n)})$ goes to zero on the set ${\cal T}_n$.
	Combining \eqref{c2eqn:eigub} and \eqref{c2eqn:chiub}, we see that $\Pi^\infty(\theta:S_\theta\in {\cal S}_n\,|\,Y^{(n)})\mathbbm{1}_{{\cal T}_n}$ is bounded by
	\begin{align*}
	&\sum_{s=s_0+1}^{K_1s_\star }\frac{\pi_p(s)\binom{p}{s_0}\binom{p-s_0}{s-s_0}}{\pi_p(s_0)\binom{p}{s}}\left(\frac{\lambda\sqrt{\pi}}{\sqrt{2}}\right)^{s-s_0}\left(\frac{\sqrt{\overline\rho_0 p^b}}{C_1\phi_2(s)\lVert X\rVert_\ast}\right)^{s-s_0}\\
	&\quad\le\sum_{s=s_0+1}^{K_1s_\star }\left(\frac{A_2}{p^{A_4}}\right)^{s-s_0}\binom{s}{s_0} \left(\frac{L_3}{C_1\phi_1(K_1s_\star ) }\sqrt{\frac{K_1 s_\star \pi\overline\rho_0 p^b }{2n}}\right)^{s-s_0},
	\end{align*}
	which holds by the inequalities ${\pi_p(s)}/{\pi_p(s_0)}\le (A_2p^{-A_4})^{s-s_0}$ and $\binom{p}{s_0}\binom{p-s_0}{s-s_0}/\binom{p}{s}=\binom{s}{s_0}$. Note that for $s\le K_1s_\star $, we have that $\binom{s}{s_0}=\binom{s}{s-s_0}\le (K_1s_\star )^{s-s_0}\le (K_1C_2p^a)^{s-s_0}$ for some $C_2>0$. Hence, the preceding display goes to zero provided that $a-A_4+b/2<0$ since $s_\star =o(n)$. This condition can be translated to $a<A_4-1$ by choosing $b$ arbitrarily close to 2.
\end{proof}

\section{Proofs for the applications}

\subsection{Proof of Theorem~\ref{c2thm:exmmissing}}
We first verify the conditions for Theorem~\ref{c2thm:thetacon} to prove assertions \ref{c2thm:missingcon} and \ref{c2thm:missingcon2}. 
\begin{itemize}[leftmargin=0.4cm, itemsep=0.0cm]
	\item {\it Verification of \rm\ref{c2con:maxfrob}}: 
	Let $\bar\sigma_{jk}$ be the $(j,k)$th element of $\Sigma-\Sigma_0$. Observe that $d_n^2(\Sigma,\Sigma_0)$ is equal to
	\begin{align}
	\begin{split}
\frac{1}{n}\sum_{i=1}^n\lVert E_i^T(\Sigma-\Sigma_0) E_i\rVert_{\rm F}^2
	&= \frac{1}{n}\sum_{j=1}^{\overline m}\sum_{k=1}^{\overline m}\left[\bar\sigma_{jk}^2 \sum_{i=1}^n e_{ij}e_{ik}\right]\gtrsim\frac{1}{c_n}\lVert\Sigma-\Sigma_0\rVert_{\rm F}^2.
	\end{split}
	\label{c2eqn:examem2}
	\end{align}
	Hence, we see that $c_n$ has the same role as $a_n$. We also have $e_n=0$ as the true $\Sigma_0$ belongs to the support of the prior.
	\item {\it Verification of \rm\ref{c2con:thm1prcon}}: Note that 
	\begin{align}
	d_n^2(\Sigma_1,\Sigma_2)=\frac{1}{n}\sum_{i=1}^n\lVert E_i^T(\Sigma_1-\Sigma_2)E_i\rVert_{\rm F}^2\le\lVert\Sigma_1-\Sigma_2\rVert_{\rm F}^2,
	\label{c2eqn:examem1}
	\end{align}
	for every $\Sigma_1,\Sigma_2\in{\cal H}$. Hence
	we obtain that for every $\bar\epsilon_n>n^{-1/2}$,
	\begin{align*}
	\log \Pi(d_n(\Sigma,\Sigma_0)\le\bar\epsilon_n)\ge\log \Pi(\lVert \Sigma-\Sigma_0\rVert_{\rm F}\le\bar\epsilon_n)\gtrsim \log \bar\epsilon_n\gtrsim -\log n,
	\end{align*}
	since $1\lesssim\rho_{\min}(\Sigma_0)\le\rho_{\max}(\Sigma_0)\lesssim 1$.
	This leads us to choose $\bar\epsilon_n=\sqrt{(\log n )/ n}$ for \ref{c2con:thm1prcon} to be satisfied.
	\item {\it Verification of \rm\ref{c2con:trueparcon1}}: The assumption $\lVert\theta_0\rVert_\infty\lesssim \lambda^{-1}\log p$ given in the theorem directly satisfies \ref{c2con:trueparcon1}.
	\item {\it Verification of \rm\ref{c2con:trueparcon2}}: We have the  inequalities $\rho_{\min}(\Sigma_0)\le\rho_{\min}(E_i^T\Sigma_0 E_i)\le\rho_{\max}(E_i^T\Sigma_0 E_i)\le \rho_{\max}(\Sigma_0)$ for every $i\le n$ as $E_i^T\Sigma_0 E_i$ is a principal submatrix of $\Sigma_0$. Hence \ref{c2con:trueparcon2} is directly satisfied by the assumption on $\Sigma_0$.
	\item {\it Verification of \rm\ref{c2con:thm2mod1}}:  For a sufficiently large $M>0$ and $s_\star =s_0\vee(\log n /\log p)$, choose ${\cal H}_{n}=\{\Sigma:n^{-M}\le\rho_{\min}(\Sigma)\le\rho_{\max}(\Sigma)\le e^{Ms_\star \log p}\}$. Since $E_i^T\Sigma E_i$ is a principal submatrix of $\Sigma$, we have $\rho_{\min}(E_i^T\Sigma E_i)\ge\rho_{\min}(\Sigma)\ge n^{-M}$ for every $i\le n$ and $\Sigma\in{\cal H}_{n}$. Hence the minimum eigenvalue condition \eqref{c2eqn:thm2c1} is satisfied with $\log \gamma_n\asymp\log  n$.
	Also, the entropy relative to $d_n$ is given by
	\begin{align*}
	&\log N\left(\frac{1}{6\overline m n^{M+3/2}},{\cal H}_n,d_n\right)\\
	&\quad \le\log N\left(\frac{1}{6\overline m n^{M+3/2}},\left\{\Sigma:\lVert\Sigma\rVert_{\rm F}\le \sqrt{\overline m}e^{Ms_\star \log p}\right\},\lVert\cdot\rVert_{\rm F}\right)\\
	&\quad \lesssim \log n+{s_\star \log p}.
	\end{align*}
	The entropy condition in \eqref{c2eqn:thm2c2} is thus satisfied if we choose $\epsilon_n=\sqrt{(s_\star \log p)/n}$. 
	To verify the sieve condition \eqref{c2eqn:thm2c3}, note that for some positive constants $b_1$, $b_2$, $b_3$, $b_4$ and $b_5$, an inverse Wishart distribution satisfies 
	\begin{align}
	\begin{split}
	\Pi(\Sigma:\rho_{\min}(\Sigma)<n^{-M})&\le b_1 e^{-b_2 n^{b_3M}},\\
	\Pi(\Sigma:\rho_{\max}(\Sigma)> e^{Ms_\star \log p})&\le b_4 e^{-b_5 Ms_\star \log p};
	\label{c2eqn:iwtail}
	\end{split}
	\end{align}
	see, for example, Lemma 9.16 of \citet{ghosal2017fundamentals}. 
The sieve condition \eqref{c2eqn:thm2c3} is met provided that $M$ is chosen sufficiently large. Note that the condition $a_n\epsilon_n^2\rightarrow 0$ is satisfied by the assumption $c_ns_\star \log p =o(n)$. 
	\item {\it Verification of \rm\ref{c2con:sep}}: The separability condition is trivially satisfied in this example as there is no nuisance mean part.
\end{itemize}
Therefore, the contraction properties in Theorem~\ref{c2thm:thetacon} are obtained with $s_\star =s_0\vee(\log n/\log p)$, but $s_\star$ is replaced by $s_0$ since $s_0>0$ and $\log n\lesssim \log p$. The contraction rate for $\Sigma$ with respect to the Frobenius norm follows from \eqref{c2eqn:examem2}.
The optimal posterior contraction directly follows from Corollary~\ref{c2cor:opt}. Assertions \ref{c2thm:missingcon} and \ref{c2thm:missingcon2} are thus proved.

Next,
we verify conditions \ref{c2con:thm4asm0}--\ref{c2con:thm4asm2} and \ref{c2con:sepsp} to apply Theorems~\ref{c2thm:bvm}--\ref{c2thm:sel} and Corollaries~\ref{c2cor:strsel}--\ref{c2cor:strbvm}.
\begin{itemize}[leftmargin=0.4cm, itemsep=0.0cm]
	\item {\it Verification of \rm\ref{c2con:thm4asm0}--\ref{c2con:thm4asm1}}: These conditions are trivially satisfied with the zero matrix $H$ as there is no nuisance mean part. 
	\item {\it Verification of \rm\ref{c2con:thm4asm2}}: 	Since the entropy in \ref{c2con:thm4asm2} is bounded above by a constant multiple of $\log N(\delta,\{\Sigma:\lVert\Sigma-\Sigma_0\rVert_{\rm F}\le \hat M_2 \sqrt{c_n}\epsilon_n\},\lVert\cdot\rVert_{\rm F})\lesssim0\vee\log(3\hat M_2\sqrt{c_n}\epsilon_n/\delta)$ using \eqref{c2eqn:examem2} and \eqref{c2eqn:examem1}, the term in \ref{c2con:thm4asm2} is bounded by a multiple of $( s_\star \vee\sqrt{\log c_n})\sqrt{c_n(s_\star \log p)^3/n}$ by Remark~\ref{rmk:3}. This term tends to zero as  $s_\star$ can be replaced by $s_0$.
	\item {\it Verification of \rm\ref{c2con:sepsp}}: Note that $d_{B,n}(\Sigma_1,\Sigma_2)\le \lVert \Sigma_1-\Sigma_2\rVert_{\rm F}$ for every $\Sigma_1,\Sigma_2$ by \eqref{c2eqn:examem1}, and hence it suffices to show that $\cal H$ is a separable metric space with the Frobenius norm. Since the support of the prior for $\Sigma$ is Euclidean, separability with the Frobenius norm is trivial.
\end{itemize}
Hence, under \ref{c2con:combo1r}, Theorem \ref{c2thm:bvm} can be applied to obtain the distributional approximation in \eqref{c2eqn:bvm} with the zero matrix $H$. Under \ref{c2con:combo1r} and \ref{c2con:selpri}, Theorem \ref{c2thm:sel} implies the no-superset result in \eqref{c2eqn:sel}. If the beta-min condition \ref{c2con:betamin} is also met, the strong results in Corollary \ref{c2cor:strsel} and Corollary \ref{c2cor:strbvm} hold. These establish \ref{c2thm:missingbvm}--\ref{c2thm:missingbvm3}.

\subsection{Proof of Theorem~\ref{c2thm:exmmem}}
We first verify the conditions for Theorem~\ref{c2thm:thetacon} for \ref{c2thm:memcon} and \ref{c2thm:memcon2}. 
\begin{itemize}[leftmargin=0.4cm, itemsep=0.0cm]
	\item {\it Verification of \rm\ref{c2con:maxfrob}}: Since $\Delta_{\eta,i}$ is the same for every $i\le n$ and the true parameters belong to the support of the prior, we see that $a_n=1$ and $e_n=0$ satisfy \ref{c2con:maxfrob}.
	\item {\it Verification of \rm\ref{c2con:thm1prcon}}: Observe that for every $\eta_1,\eta_2\in{\cal H}$,
	\begin{align}
	\begin{split}
	\lVert\xi_{\eta_1}-\xi_{\eta_2}\rVert_2^2&=\lvert (\alpha_1-\alpha_2)+(\mu_1^T\beta_1-\mu_2^T\beta_2)
	\rvert^2+\lVert\mu_1-\mu_2\rVert_2^2\\
	&\lesssim\lvert\alpha_1-\alpha_2\rvert^2+\lVert\mu_1\rVert_2^2\lVert\beta_1-\beta_2\rVert_2^2+(\lVert\beta_2\rVert_2^2+1)\lVert\mu_1-\mu_2\rVert_2^2,\\
	\lVert\Delta_{\eta_1}-\Delta_{\eta_2}\rVert_{\rm F}^2&=|(\beta_1^T\Sigma_1\beta_1-\beta_2^T\Sigma_2\beta_2)+(\sigma_1^2-\sigma_2^2)|^2\\
	&\quad+2\lVert\Sigma_1\beta_1-\Sigma_2\beta_2\rVert_2^2+\lVert\Sigma_1-\Sigma_2\rVert_{\rm F}^2\\
	&\lesssim (\lVert\beta_1\rVert_2^2+1)^2\lVert\Sigma_1-\Sigma_2\rVert_{\rm F}^2+|\sigma_1^2-\sigma_2^2|^2\\
	&\quad+(\lVert\beta_1\rVert_2^2+\lVert\beta_2\rVert_2^2+1)\lVert\Sigma_2\rVert_{\rm F}^2\lVert\beta_1-\beta_2\rVert_2^2.
	\label{c2eqn:membo}
	\end{split}
	\end{align}
	Since $\lVert\beta_0\rVert_2$, $\lvert\sigma_0^2\rvert$, and $\lVert\Sigma_0\rVert_{\rm F}$ are bounded, it follows from the last display that there exists a constant $C_1$ such that $\lvert\alpha-\alpha_0\rvert+\lVert\beta-\beta_0\rVert_2+\lVert\mu-\mu_0\rVert_2+\lvert\sigma^2-\sigma_0^2\rvert+\lVert\Sigma-\Sigma_0\rVert_{\rm F}\le C_1\bar\epsilon_n$ implies $d_n(\eta,\eta_0)\le\bar\epsilon_n$ for any small $\bar\epsilon_n$. 
	This shows that \ref{c2con:thm1prcon} is satisfied as long as we choose $\bar\epsilon_n=\sqrt{\log n/n}$, as we have $|\alpha_0|\vee\lVert\beta_0\rVert_\infty\vee\lVert\mu_0\rVert_\infty\lesssim 1$, $\sigma_0^2\asymp1$, and $1\lesssim\rho_{\min}(\Sigma_0)\le\rho_{\max}(\Sigma_0)\lesssim 1$.
	\item {\it Verification of \rm\ref{c2con:trueparcon1}}: The assumption $\lVert\theta_0\rVert_\infty\lesssim \lambda^{-1}\log p$ given in the theorem directly satisfies \ref{c2con:trueparcon1}.
	\item {\it Verification of \rm\ref{c2con:trueparcon2}}: 	Since $\Delta_\eta$ can be written as the sum of two positive definite matrices as
	\begin{align*}
	\Delta_\eta=\begin{pmatrix}
	\beta^T\Sigma\beta & \beta^T\Sigma \\ \Sigma\beta & \Sigma
	\end{pmatrix} + \begin{pmatrix}
	\sigma^2 & 0 \\ 0 & \Psi
	\end{pmatrix},
	\end{align*}
	condition \ref{c2con:trueparcon2} is satisfied as we obtain $\sigma_0^2\wedge\rho_{\min}(\Psi)\le \rho_{\min}(\Delta_{\eta_0})\le \rho_{\max}(\Delta_{\eta_0})\le\lVert\Delta_{\eta_0} \rVert_{\rm F}$  by Weyl's inequality.
	\item {\it Verification of \rm\ref{c2con:thm2mod1}}:  For a sufficiently large $M$ and $s_\star =s_0\vee (\log n/\log p)$, choose a sieve as
	\begin{align*}
	{\cal H}_{n}&=\{(\alpha,\beta,\mu):
	\lvert\alpha\rvert^2+\lVert\beta\rVert_2^2+\lVert\mu\rVert_2^2\le n^{2M}\} \times\{\sigma:n^{-M}\le\sigma^2\le e^{Ms_\star \log p}\}\\
	&\quad\times\{\Sigma:n^{-M}\le\rho_{\min}(\Sigma)\le\rho_{\max}(\Sigma)\le e^{Ms_\star  \log p}\}.
	\end{align*}
	Then we have $\rho_{\min}(\Delta_\eta)\ge \sigma^2\wedge\rho_{\min}(\Psi)\ge n^{-M}$ for large $n$, and hence the minimum eigenvalue condition \eqref{c2eqn:thm2c1} is directly met with $\log\gamma_n\asymp\log n$ by the definition of the sieve. To see the entropy condition, observe from \eqref{c2eqn:membo} that for every $\eta_1,\eta_2\in{\cal H}_n$,
	\begin{align*}
	d_n^2(\eta_1,\eta_2)\lesssim n^{4M}e^{2Ms_\star \log p}\big(&\lvert\alpha-\alpha_0\rvert^2+\lVert\beta_1-\beta_2\rVert_2^2 +\lVert\mu_1-\mu_2\rVert_2^2\\
	&+\lVert\Sigma_1-\Sigma_2\rVert_{\rm F}^2+|\sigma_1^2-\sigma_2^2|^2\big).
	\end{align*}
	Therefore, for $\delta_n=1/(6\overline m n^{3M+3/2}e^{Ms_\star \log p})$, the entropy relative to $d_n$ is bounded above by
	\begin{align*}
	&\log N\left(\delta_n,\{(\alpha,\beta,\mu):\lvert\alpha\rvert^2+\lVert\beta\rVert_2^2+\lVert\mu\rVert_2^2\le n^{2M}\},\lVert\cdot\rVert_2\right)\\
	&+\log N\left(\delta_n,\{\sigma:\sigma^2 \le e^{Ms_\star \log p}\},\lVert\cdot\rVert_2\right)\\
	& +\log N\left(\delta_n,\{\Sigma:\lVert\Sigma\rVert_{\rm F}\le \sqrt{q}e^{Ms_\star \log p}\},\lVert\cdot\rVert_{\rm F}\right),
	\end{align*}
	each summand of which is bounded by a multiple of $\log n+s_\star \log p$. 
	This shows that the choice $\epsilon_n=\sqrt{(s_\star \log p)/n}$ satisfies the entropy condition in \eqref{c2eqn:thm2c2}. Further, it is easy to see that condition \eqref{c2eqn:thm2c3} holds using the tail bounds for normal and inverse Wishart distributions as in \eqref{c2eqn:iwtail}.
	\item {\it Verification of \rm\ref{c2con:sep}}: Note that the mean of $Y$ is expressed as $X\theta+Z\xi_{\eta}$ for $Z=1_n\otimes I_{q+1}$.
	Since the condition $\varsigma_{\min}([X_S^\ast,1_n])\gtrsim 1$ implies $\varsigma_{\min}([X_S,Z])\gtrsim 1$, condition \ref{c2con:sep} is satisfied by Remark~\ref{rmk:1}.
\end{itemize}
Therefore we obtain the contraction properties of the posterior distribution as in \eqref{c2eqn:thm3} with $s_\star$ replaced by $s_0$ as $s_0>0$ and $\log n\lesssim \log p$.
The rates for $\eta$ with respect to more concrete metrics than $d_n$ can now be obtained.
Note that for small $\delta>0$, $d_n(\eta,\eta_0)\le\delta$ directly implies $\lVert\mu-\mu_0\rVert_2\le\delta$ and $\lVert\Sigma-\Sigma_0\rVert_{\rm F}\le\delta$ by the definition of $d_n$. For $\beta$, observe that
\begin{align*}
\lVert\beta-\beta_0\rVert_2&\le\lVert\Sigma^{-1}\rVert_{\rm sp}\lVert\Sigma(\beta-\beta_0)\rVert_2\\
& \le\lVert\Sigma^{-1}\rVert_{\rm sp}(\lVert\Sigma\beta-\Sigma_0\beta_0\rVert_2+\lVert\Sigma-\Sigma_0\rVert_{\rm F}\lVert\beta_0\rVert_2)\\
& \lesssim\lVert\Sigma^{-1}\rVert_{\rm sp} \delta.
\end{align*}
Since $\lVert\Sigma^{-1}\rVert_{\rm sp}$ is bounded as $\lVert\Sigma-\Sigma_0\rVert_{\rm F}\le\delta$, the preceding display implies $\lVert\beta-\beta_0\rVert_2\lesssim\delta$. Moreover, we have
\begin{align*}
|\alpha-\alpha_0|&\le|\mu^T\beta-\mu_0^T\beta_0|+\delta\\
&\lesssim\lVert\mu\rVert_2\lVert\beta-\beta_0\rVert_2+\lVert\beta_0\rVert_2\lVert\mu-\mu_0\rVert_2+\delta\\
&\lesssim(\lVert\mu\rVert_2+1)\delta,
\end{align*}
and
\begin{align*}
|\sigma^2-\sigma_0^2|&\le|\beta^T\Sigma\beta-\beta_0^T\Sigma_0\beta_0|+|(\beta^T\Sigma\beta+\sigma^2)-(\beta_0^T\Sigma_0\beta_0+\sigma_0^2)|\\
&\le\lVert\beta\rVert_2\lVert\Sigma\beta-\Sigma_0\beta_0\rVert_2+\lVert\beta_0\rVert_2\lVert\Sigma_0\rVert_{\rm sp}\lVert\beta-\beta_0\rVert_2+\delta\\
&\lesssim(\lVert\beta\rVert_2+1)\delta.
\end{align*}
These show that $|\alpha-\alpha_0|+|\sigma^2-\sigma_0^2|\lesssim\delta$ as $\lVert\mu\rVert_2$ and $\lVert\beta\rVert_2$ are bounded. We finally conclude that $\lvert\alpha-\alpha_0\rvert+\lVert\beta-\beta_0\rVert_2+\lVert\mu-\mu_0\rVert_2+\lvert\sigma^2-\sigma_0^2\rvert+\lVert\Sigma-\Sigma_0\rVert_{\rm F}$ contracts at the same rate of $d_n$.
The optimal posterior contraction is directly obtained by Corollary~\ref{c2cor:opt}. Thus assertions \ref{c2thm:memcon} and \ref{c2thm:memcon2} hold.

Next,	we verify conditions \ref{c2con:thm4asm0}--\ref{c2con:thm4asm2}  and \ref{c2con:sepsp} to apply Theorems~\ref{c2thm:bvm}--\ref{c2thm:sel} and Corollaries~\ref{c2cor:strsel}--\ref{c2cor:strbvm}. The orthogonal projection defined by $H=\tilde Z(\tilde Z^T\tilde Z)^{-1}\tilde Z^T$ with $\tilde Z=1_n\otimes \Delta_{\eta_0}^{-1/2}$ is used to check the conditions.
\begin{itemize}[leftmargin=0.4cm, itemsep=0.0cm]
	\item {\it Verification of \rm\ref{c2con:thm4asm0}}: 	For $H$ defined above, it is easy to see that the first condition of \ref{c2con:thm4asm0} is satisfied. The second condition is directly satisfied by Remark~\ref{rmk:2}.
	\item {\it Verification of \rm\ref{c2con:thm4asm1}}: Choose a map $(\alpha,\beta,\mu,\sigma^2,\Sigma)\mapsto(\alpha+n^{-1}1_n^T X^\ast(\theta-\theta_0),\beta,\mu,\sigma^2,\Sigma)$ for  $\eta\mapsto\tilde\eta_n(\theta,\eta)$.
	To check \ref{c2con:thm4asm1}, we shall verify that this map induces $\Phi(\tilde\eta_n(\theta,\eta))=(\tilde\xi_\eta+H\tilde X(\theta-\theta_0),\tilde\Delta_\eta)$ as follows. 
	Note that for matrices $R_k$, $k=1,\dots,6$, we have the properties of the Kronecker product that $(R_1\otimes R_2)(R_3\otimes R_4)=(R_1R_2\otimes R_3R_4)$ and $(R_5\otimes R_6)^{-1}=R_5^{-1}\otimes R_6^{-1}$ if the matrices allow such operations. Using these properties, we see that $H$ satisfies
	\begin{align*}
	H&=(1_n\otimes \Delta_{\eta_0}^{-1/2}) (1_n^T1_n\otimes\Delta_{\eta_0}^{-1})^{-1}(1_n\otimes \Delta_{\eta_0}^{-1/2})^T\\
	&=\frac{1}{n}(1_n\otimes \Delta_{\eta_0}^{-1/2}) \Delta_{\eta_0}(1_n\otimes \Delta_{\eta_0}^{-1/2})^T\\
	&=\frac{1}{n}(1_n\otimes I_{q+1})(1_n^T\otimes I_{q+1})\\
	&=\frac{1}{n} (1_n1_n^T \otimes I_{q+1}).
	\end{align*}
	Hence,
	\begin{align*}
	Z(\tilde Z^T\tilde Z)^{-1}\tilde Z^T\tilde X(\theta-\theta_0)&=(I_n\otimes\Delta_{\eta_0}^{1/2})H(I_n\otimes\Delta_{\eta_0}^{-1/2}) X(\theta-\theta_0)\\
	&=H X(\theta-\theta_0)=1_n\otimes\begin{pmatrix}
	n^{-1} 1_n^T X^\ast(\theta-\theta_0)\\0_{q\times 1}
	\end{pmatrix},
	\end{align*}
	which implies that the shift only for $\alpha$ as in the given map provides  $\Phi(\tilde\eta_n(\theta,\eta))=(\tilde\xi_\eta+H\tilde X(\theta-\theta_0),\tilde\Delta_\eta)$. 
	Without loss of generality, we assume that the standard normal prior is used for $\alpha$.
	Now, observe that
	\begin{align*}
	\left\lvert\log \frac{d\Pi_{n,\theta}} {d\Pi_{n,\theta_0}}(\eta)\right\rvert&\lesssim\left\lvert \alpha^2-(\alpha+ n^{-1}1_n^TX^\ast(\theta-\theta_0))^2\right\rvert\\
	&\le 2|\alpha||n^{-1}1_n^TX^\ast(\theta-\theta_0)|+(n^{-1}1_n^TX^\ast(\theta-\theta_0))^2,
	\end{align*}
	since the priors for the other parameters cancel out due to invariance.
	One can note that
	\begin{align*}
	\sup_{\eta\in\widehat{\cal H}_n}|\alpha|\lesssim s_\star \sqrt{(\log p)/n}+|\alpha_0|\lesssim 1,
	\end{align*}
	and 
	\begin{align*}
	\frac{1}{\sqrt{n}}\sup_{\theta\in\widehat\Theta_n}\lVert X (\theta-\theta_0)\rVert_2\lesssim s_\star \sqrt{(\log p) /n}.
	\end{align*}
Thus, condition \ref{c2con:thm4asm1} is satisfied. 
	\item {\it Verification of \rm\ref{c2con:thm4asm2}}: Note again that $d_{B,n}(\eta,\eta_0)\lesssim \lVert\Sigma-\Sigma_0\rVert_{\rm F}+ |\sigma^2-\sigma_0^2| +\lVert\beta-\beta_0\rVert_2$  for every $\eta\in\widehat{\cal H}_n$. The inequality also holds for the other direction for every $\eta\in\widehat{\cal H}_n$, by the same argument used for the recovery in the proof of Theorem~\ref{c2thm:exmmem}, \ref{c2thm:memcon}--\ref{c2thm:memcon2}. 
	Hence, for some constants $C_1,C_2>0$, the entropy in \ref{c2con:thm4asm2} is bounded above by 
	\begin{align*}
	&\log N\left(C_1\delta,\left\{\beta:\lVert\beta-\beta_0\rVert_2\le C_2 \hat M_2 \epsilon_n\right\},\lVert\cdot\rVert_2\right) \\
	&+\log N\left(C_1\delta,\left\{\sigma^2:\lvert\sigma^2-\sigma_0^2\rvert\le C_2 \hat M_2 \epsilon_n\right\},\lvert\cdot\rvert\right)\\
	&+\log N\left(C_1\delta,\left\{\Sigma:\lVert\Sigma-\Sigma_0\rVert_{\rm F}\le C_2 \hat M_2 \epsilon_n\right\},\lVert\cdot\rVert_{\rm F}\right).
	\end{align*}
	Since all nuisance parameters are of fixed dimensions, the last display is bounded  by a multiple of $0\vee\log(3C_2\hat M_2\epsilon_n/C_1\delta)$ for every $\delta>0$, so that \ref{c2con:thm4asm2} is bounded by $(s_\star ^5\log^3 p/n)^{1/2}$ by Remark~\ref{rmk:3}. Since $s_\star\lesssim s_0$ in this case, the condition is verified.
	\item {\it Verification of \rm\ref{c2con:sepsp}}: Note that by \eqref{c2eqn:membo}, $d_{B,n}(\eta_1,\eta_2)\lesssim \lVert\Sigma_1-\Sigma_2\rVert_{\rm F}+ |\sigma_1^2-\sigma_2^2| +\lVert\beta_1-\beta_2\rVert_2$  for every $\eta_1,\eta_2\in\widehat{\cal H}_n$. Since each of the parameter spaces of $\Sigma$, $\sigma^2$, and $\beta$ is a separable metric space with each of these norms, \ref{c2con:sepsp} is satisfied.
\end{itemize}
Therefore, under \ref{c2con:combo1r}, Theorem~\ref{c2thm:bvm} implies that the distributional approximation in \eqref{c2eqn:bvm} holds.
Under \ref{c2con:combo1r} and \ref{c2con:selpri}, we obtain the no-superset result in \eqref{c2eqn:sel}.
The remaining assertions in the theorem are direct consequences of Corollary~\ref{c2cor:strsel} and Corollary~\ref{c2cor:strbvm} if the beta-min condition \ref{c2con:betamin} is also satisfied.
These prove  \ref{c2thm:membvm}--\ref{c2thm:membvm3}.

We complete the proof by showing that the covariance matrix of the nonzero part can be written as in the theorem. For given $S$, we obtain
\begin{align*}
&\tilde X_S^T(I_{n(q+1)}-H)\tilde X_S \\
&\quad=X_S^{\ast T}\left(I_n\otimes \{\Delta_{\eta_0}^{-1/2}\}_{\cdot 1}^T\right)\left(I_n\otimes I_{q+1}- H\right) \left(I_n\otimes \{\Delta_{\eta_0}^{-1/2}\}_{\cdot 1}\right)X_S^\ast\\
&\quad = \{\Delta_{\eta_0}^{-1/2}\}_{\cdot 1}^T \{\Delta_{\eta_0}^{-1/2}\}_{\cdot 1} X_S^{\ast T} H^\ast X_S^\ast,
\end{align*}
where $\{\Delta_{\eta_0}^{-1/2}\}_{\cdot 1}$ is the first column of $\Delta_{\eta_0}^{-1/2}$. Note that $\{\Delta_{\eta_0}^{-1/2}\}_{\cdot 1}^T \{\Delta_{\eta_0}^{-1/2}\}_{\cdot 1}=\{\Delta_{\eta_0}^{-1}\}_{1,1}$, where $\{\Delta_{\eta_0}^{-1}\}_{1,1}$ is the top-left element of $\Delta_{\eta_0}^{-1}$, which is equal to 
$(\beta_0^T\Sigma_0\beta_0+\sigma_0^2-\beta_0^T\Sigma_0(\Sigma_0+\Psi)^{-1}\Sigma_0\beta_0)^{-1}=(\sigma_0^2+\beta_0^T\Sigma_0(\Sigma_0+\Psi)^{-1}\Psi\beta_0)^{-1}$ by direct calculations. For the mean $\hat\theta_S$, observe that
\begin{align*}
&\tilde X_S^T(I_{n(q+1)}-H)(U+\tilde X\theta_0)\\
&\quad=X_S^{\ast T}\left(I_n\otimes \{\Delta_{\eta_0}^{-1/2}\}_{\cdot 1}^T\right)\left(I_n\otimes I_{q+1}-\frac{1}{n}1_n1_n^T\otimes I_{q+1}\right)\\ 
& \qquad\times\left\{\left(I_n\otimes \{\Delta_{\eta_0}^{-1/2}\}_{\cdot 1}\right) \Big(Y^\ast-(\alpha_0+\mu_0^T\beta_0)1_n\right)\\
& \quad\qquad\quad+\left(I_n\otimes \{\Delta_{\eta_0}^{-1/2}\}_{\cdot (-1)}\right)\left(W-1_n\otimes\mu_0\right)\Big\},
\end{align*}
where $\{\Delta_{\eta_0}^{-1/2}\}_{\cdot (-1)}$ is the submatrix of $\Delta_{\eta_0}^{-1/2}$ consisting of columns except for  $\{\Delta_{\eta_0}^{-1/2}\}_{\cdot 1}$ the first column. Since $\{\Delta_{\eta_0}^{-1/2}\}_{\cdot 1}^T \{\Delta_{\eta_0}^{-1/2}\}_{\cdot (-1)}=\{\Delta_{\eta_0}^{-1}\}_{1,(-1)}$, where $\{\Delta_{\eta_0}^{-1}\}_{1,(-1)}$ is the first row of $\Delta_{\eta_0}^{-1}$ with the top-left element excluded, the last display is equal to
\begin{align*}
X_S^{\ast T}\Big\{H^\ast\Big[&\{\Delta_{\eta_0}^{-1}\}_{1,1} \left(Y^\ast-(\alpha_0+\mu_0^T\beta_0)1_n\right)\\ &+\left(I_n\otimes\{\Delta_{\eta_0}^{-1}\}_{1,(-1)}\right)\left(W-1_n\otimes\mu_0\right)\Big]\Big\}.
\end{align*}
As we have $\{\Delta_{\eta_0}^{-1}\}_{1,(-1)}=-\{\Delta_{\eta_0}^{-1}\}_{1,1}^{-1}\beta_0^T\Sigma_0(\Sigma_0+\Psi)^{-1}$ by direct calculations, it follows that
\begin{align*}
\hat\theta_S=\left(X_S^{\ast T} H^\ast X_S^\ast\right)^{-1}X_S^{\ast T}\Big\{H^\ast\Big[& \left(Y^\ast-(\alpha_0+\mu_0^T\beta_0)1_n\right)\\
&\!-\left(I_n\otimes(\beta_0^T\Sigma_0(\Sigma_0+\Psi)^{-1})\right)\left(W-1_n\otimes\mu_0\right)\Big]\Big\}.
\end{align*}
This completes the proof.

\subsection{Proof of Theorem~\rm\ref{c2thm:exmparcor}}
We shall verify the conditions for the posterior contraction in Theorem~\ref{c2thm:thetacon} to prove \ref{c2thm:parcorcon}--\ref{c2thm:parcorcon2}. 
First we give the bounds for the eigenvalues of each correlation matrix. It can be shown that
\begin{align}
1-\alpha=\rho_{\min}\left(G_i^{\rm CS}(\alpha)\right)&\le\rho_{\max}\left(G_i^{\rm CS}(\alpha)\right)= 1+(m_i-1)\alpha,\label{c2eqn:eigencs}\\
\frac{1-\alpha^2}{(1+|\alpha|)^2}\le\rho_{\min}\left(G_i^{\rm AR}(\alpha)\right)&\le\rho_{\max}\left(G_i^{\rm AR}(\alpha)\right)\le \frac{1-\alpha^2}{(1-|\alpha|)^2},\label{c2eqn:eigenar}\\
1-2|\alpha|\le\rho_{\min}\left(G_i^{\rm MA}(\alpha)\right)&\le\rho_{\max}\left(G_i^{\rm MA}(\alpha)\right)\le 1+2|\alpha|.\label{c2eqn:eigenma}
\end{align}
The first assertion in \eqref{c2eqn:eigencs} follows directly from the identity $\rho_k(G_i^{\rm CS}(\alpha))=\rho_k(\alpha1_{m_i}1_{m_i}^T)+1-\alpha$ for every $k\le m_i$. For \eqref{c2eqn:eigenar}, see Theorem 2.1 and Theorem 3.5 of \citet{fikioris2018spectral}. The assertion in \eqref{c2eqn:eigenma} is due to Theorem 2.2 of \citet{kulkarni1999eigenvalues}.
\begin{itemize}[leftmargin=0.4cm, itemsep=0.0cm]
	\item {\it Verification of \rm\ref{c2con:maxfrob}}: 	For the autoregressive correlation matrix, note that 
	\begin{align*}
	&\max_{1\le i \le n}\left\lVert \sigma^2 G_i^{\rm AR}(\alpha)-\sigma_0^2 G_i^{\rm AR}(\alpha_0)\right\rVert_{\rm F}^2\\
	&\quad =\overline m (\sigma^2-\sigma_0^2)^2+2\sum_{k=1}^{\overline m -1}(\overline m -k)(\sigma^2\alpha^k-\sigma_0^2\alpha_0^k)^2.
	\end{align*}
	Using $\overline m n\asymp n_\ast$, we have that
	\begin{align*}
	\sum_{k=1}^{\overline m -1}(\overline m -k)(\sigma^2\alpha^k-\sigma_0^2\alpha_0^k)^2
	&\lesssim \frac{1}{n}\sum_{k=1}^{\overline m -1}(\sigma^2\alpha^k-\sigma_0^2\alpha_0^k)^2\sum_{i=1}^n \{(m_i -k)\vee 0\}\\
	&=\frac{1}{n}\sum_{i=1}^n\sum_{k=1}^{m_i -1}(m_i -k)(\sigma^2\alpha^k-\sigma_0^2\alpha_0^k)^2,
	\end{align*}
	and hence 
	\begin{align*}
	\max_{1\le i\le n}\lVert \sigma^2 G_i^{\rm AR}(\alpha)-\sigma_0^2 G_i^{\rm AR}(\alpha_0)\rVert_{\rm F}^2\lesssim \frac{1}{n}\sum_{i=1}^n\lVert \sigma^2 G_i^{\rm AR}(\alpha)-\sigma_0^2 G_i^{\rm AR}(\alpha_0)\rVert_{\rm F}^2.
	\end{align*}
	This gives us $a_n\asymp 1$ for the autoregressive matrices.
	Similarly, we can also show that $a_n\asymp 1$ satisfies \ref{c2con:maxfrob} for the compound-symmetric and the moving average correlation matrices. Also, we have $e_n=0$ for \ref{c2con:maxfrob}  as the true parameter values $\alpha_0$ and $\sigma_0^2$ are in the support of the prior.
	\item {\it Verification of \rm\ref{c2con:thm1prcon}}: 
	Since the nuisance parameters are of fixed dimensions, condition \ref{c2con:thm1prcon} is satisfied with $\bar\epsilon_n=\sqrt{(\log n)/n}$ due to the restricted range of the true parameters, $\sigma_0^2\asymp 1$ and $\alpha_0\in[b_1+\epsilon,b_2-\epsilon]$ for some fixed $\epsilon>0$. 
	\item {\it Verification of \rm\ref{c2con:trueparcon1}}: The assumption $\lVert\theta_0\rVert_\infty\lesssim \lambda^{-1}\log p$ given in the theorem directly satisfies \ref{c2con:trueparcon1}.
	\item {\it Verification of \rm\ref{c2con:trueparcon2}}: 
	Using \eqref{c2eqn:eigencs}--\eqref{c2eqn:eigenma}, we see that for the compound-symmetric correlation matrix, condition \ref{c2con:trueparcon2} is satisfied with the bounded range of the true parameters provided that $\overline m$ is bounded. For the other correlation matrices, condition \ref{c2con:trueparcon2} is satisfied even with increasing $\overline m$.
	\item {\it Verification of \rm\ref{c2con:thm2mod1}}:  	For a sufficiently large $M>0$ and $s_\star =s_0\vee (\log n/\log p)$, choose a sieve ${\cal H}_n=\{\sigma^2:n^{-M}\le \sigma^2\le e^{Ms_\star \log p}\}\times\{\alpha:b_1+n^{-M}\le\alpha\le b_2-n^{-M}\}$. Then using \eqref{c2eqn:eigencs}--\eqref{c2eqn:eigenma}, it is easy to see that the minimum eigenvalue of each correlation matrix is bounded below by a polynomial in $n$, which implies that condition \eqref{c2eqn:thm2c1} is satisfied with $\log \gamma_n\asymp \log n$.  For the entropy calculation, note that for every type of correlation matrix,
	\begin{align}
	\begin{split}
	d_n^2(\eta_1,\eta_2)&=\frac{1}{n}\sum_{i=1}^n\lVert\sigma_1^2 G_i(\alpha_1)-\sigma_2^2 G_i(\alpha_2)\rVert_{\rm F}^2\\
	&\le\frac{1}{n}\sum_{i=1}^n\left\{(\sigma_1^2-\sigma_2^2)^2\lVert G_i(\alpha_1)\rVert_{\rm F}^2+\sigma_2^4\lVert G_i(\alpha_1)-G_i(\alpha_2)\rVert_{\rm F}^2\right\}.
	\end{split}
	\label{c2eqn:parconbo}
	\end{align}
	From the identity
	$\alpha_1^k-\alpha_2^k=(\alpha_1-\alpha_2)\sum_{j=0}^{k-1}\alpha_1^j \alpha_2^{k-1-j}$ for every integer $k\ge 1$, we have that
	$|\alpha_1^k -\alpha_2^k|\lesssim k  |\alpha_1-\alpha_2|$ for every $\alpha_1,\alpha_2\in(b_1,b_2)$. By this inequality we obtain $\lVert G_i(\alpha_1)-G_i(\alpha_2)\rVert_{\rm F}^2\lesssim \overline m^4|\alpha_1-\alpha_2|^2$ for every correlation matrix. Then, the last display is bounded by a multiple of $\overline m^2(\sigma_1^2-\sigma_2^2)^2+e^{2Ms_\star \log p}\overline m^4(\alpha_1-\alpha_2)^2$ for every $\eta_1,\eta_2\in{\cal H}_n$.
	The entropy in \eqref{c2eqn:thm2c2} is thus bounded by
	\begin{align*}
	\log N\big(\delta_n,\{\sigma^2:0< \sigma^2\le e^{Ms_\star \log p}\},|\cdot|\big)+
	\log N\big(\delta_n,\{\alpha: 0<\alpha<1\},|\cdot|\big),
	\end{align*}
	for $\delta_n=(6\overline m^3 n^{3/2+C_1} e^{Ms_\star \log p})^{-1}$ with some constant $C_1>0$. It can be easily checked that each term in the last display is bounded by a multiple of $s_\star \log p$, by which the entropy condition in \eqref{c2eqn:thm2c2} is satisfied with $\epsilon_n=\sqrt{(s_\star \log p)/n}$. Using the tail bounds of inverse gamma distributions and properties of the density $\Pi(d\alpha)$ near the boundaries, condition \eqref{c2eqn:thm2c3} is satisfied as long as $M$ is chosen sufficiently large. 
	\item {\it Verification of \rm\ref{c2con:sep}}: The separation condition is trivially satisfied as there is no nuisance mean part.
\end{itemize}
Therefore, we obtain the posterior contraction properties of $\theta$ with $s_\star =s_0\vee (\log n/\log p)$ by Theorem~\ref{c2thm:thetacon}.
The term $s_\star$ can be replaced by $s_0$ since $s_0>0$ and $\log n\lesssim \log p$.
Since we have $m_i(\sigma^2-\sigma_0^2)^2\le\lVert\sigma^2 G_i(\alpha)-\sigma_0^2 G_i(\alpha_0)\rVert_{\rm F}^2$ by the diagonal entries of each matrix, the contraction rate $\sqrt{(s_0 \log p)/(\overline m n)}$ is obtained for $\sigma^2$ with respect to the $\ell_2$-norm, for every correlation matrix, as $\overline m n\asymp n_\ast$.
In particular, for the compound-symmetric correlation matrix, this rate is reduced to $\sqrt{(s_0 \log p)/n}$ since $\overline m$ is bounded in that case.
We also have $m_i(\sigma^2\alpha-\sigma_0^2\alpha_0)^2\le\lVert\sigma^2 G_i(\alpha)-\sigma_0^2 G_i(\alpha_0)\rVert_{\rm F}^2$ for every correlation matrix, as there are more than $m_i$ entries that is equal to $\sigma^2\alpha-\sigma_0^2\alpha_0$. Hence, by the relation
$|\alpha-\alpha_0|\lesssim |\sigma^2\alpha-\sigma_0^2\alpha_0|+|\alpha||\sigma^2-\sigma_0^2|$, the same rate is also obtained for $\alpha$ relative to the $\ell_2$-norm.
The optimal posterior contraction directly follows from Corollary~\ref{c2cor:opt}.
Thus assertions \ref{c2thm:parcorcon}--\ref{c2thm:parcorcon2} hold.

Next, we verify conditions \ref{c2con:thm4asm0}--\ref{c2con:thm4asm2} and \ref{c2con:sepsp} to apply Theorems~\ref{c2thm:bvm}--\ref{c2thm:sel} and Corollaries~\ref{c2cor:strsel}--\ref{c2cor:strbvm}.
\begin{itemize}[leftmargin=0.4cm, itemsep=0.0cm]
	\item {\it Verification of \rm\ref{c2con:thm4asm0}--\ref{c2con:thm4asm1}}: These conditions are trivially satisfied with the zero matrix $H$ since there is no nuisance mean part.
	\item {\it Verification of \rm\ref{c2con:thm4asm2}}: 	Using the results of contraction rates of $\sigma^2$ and $\alpha$, note that there exists a constant $C_2>0$ such that $\{\eta\in {\cal H}:d_{B,n}(\eta,\eta_0)\le \hat M_2\epsilon_n\} \subset \{\sigma^2 : |\sigma^2-\sigma_0^2|\le C_2\epsilon_n/\sqrt{\overline m}\}\times \{\alpha:|\alpha-\alpha_0|\le C_2\epsilon_n/\sqrt{\overline m}\}$. Thus the entropy in \ref{c2con:thm4asm2} is bounded by $0\vee2\log (3C_2\epsilon_n/\sqrt{\overline m}\delta)$. By Remark~\ref{rmk:3}, \ref{c2con:thm4asm2} is bounded by a multiple of $\{(s_\star ^5\log^3 p)/n\}^{1/2}$, which goes to zero by the assumption since $s_\star\lesssim s_0$.
	\item {\it Verification of \rm\ref{c2con:sepsp}}: Using \eqref{c2eqn:parconbo}, we have $d_{B,n}(\eta_1,\eta_2)\lesssim \overline m |\sigma_1^2-\sigma_2^2| +\overline m^2|\alpha_1-\alpha_2|$ for every $\eta_1,\eta_2\in \widehat{\cal H}_n$. Since the parameter spaces of $\alpha$ and $\sigma^2$ are Euclidean and hence separable under the $\ell_2$-metric, condition \ref{c2con:sepsp} is satisfied.
\end{itemize}
Therefore, under \ref{c2con:combo1r}, the distributional approximation in \eqref{c2eqn:bvm} holds with the zero matrix $H$ by Theorem \ref{c2thm:bvm}.
Under \ref{c2con:combo1r} and \ref{c2con:selpri}, Theorem \ref{c2thm:sel} implies that the no-superset result in \eqref{c2eqn:sel} holds. The strong results in Corollary \ref{c2cor:strsel} and Corollary \ref{c2cor:strbvm} follow explicitly from the beta-min condition \ref{c2con:betamin}. These prove \ref{c2thm:parcorbvm}--\ref{c2thm:parcorbvm3}.

\subsection{Proof of Theorem~\rm\ref{c2thm:exmmm}}
We verify the conditions for the posterior contraction in Theorem~\ref{c2thm:thetacon} to show \ref{c2thm:mmcon}--\ref{c2thm:mmcon2}. 
\begin{itemize}[leftmargin=0.4cm, itemsep=0.0cm]
	\item {\it Verification of \rm\ref{c2con:maxfrob}}: 	Using the assumption $\max_i\lVert Z_i\rVert_{\rm sp}\lesssim1$, note that
	\begin{align} 
	\begin{split}
	&\max_{1\le i\le n}\lVert Z_i(\Psi-\Psi_0)Z_i^T \rVert_{\rm F}^2 \\
	&\quad\le\lVert\Psi-\Psi_0\rVert_{\rm F}^2 \max_{1\le i\le n}\lVert Z_i\rVert_{\rm sp}^4\\
	&\quad\lesssim \frac{1}{\sum_{i=1}^n\mathbbm{1}(m_i\ge q)} \sum_{i:m_i\ge q}\lVert Z_i(\Psi-\Psi_0) Z_i^T\rVert_{\rm F}^2 \lVert (Z_i^TZ_i)^{-1} Z_i^T\rVert_{\rm sp}^4\\
	&\quad\lesssim \frac{1}{n}\sum_{i=1}^n\lVert Z_i(\Psi-\Psi_0)Z_i^T \rVert_{\rm F}^2,
	\end{split}
	\label{c2eqn:maxavebomm}
	\end{align}
	where the last inequality holds since
	$\min_i\{\varsigma_{\min}( Z_i) :m_i\ge q\}\gtrsim 1$ and  $\sum_{i=1}^n\mathbbm{1}(m_i\ge q)\asymp n$. Thus we have $a_n\asymp 1$ and $e_n=0$.
	\item {\it Verification of \rm\ref{c2con:thm1prcon}}: The  condition is satisfied with $\bar\epsilon_n=\sqrt{(\log n)/n}$ as $\Psi$ is fixed dimensional and we have $1\lesssim\rho_{\min}(\Psi_0)\le\rho_{\max}(\Psi_0)\lesssim 1$.
	\item {\it Verification of \rm\ref{c2con:trueparcon1}}: The assumption $\lVert\theta_0\rVert_\infty\lesssim \lambda^{-1}\log p$ given in the theorem directly satisfies \ref{c2con:trueparcon1}.
	\item {\it Verification of \rm\ref{c2con:trueparcon2}}: By Weyl's inequality, we obtain that
	\begin{align}
	\min_{1\le i\le n}\rho_{\min}(\sigma^2 I_{m_i}+Z_i \Psi_0 Z_i^T)&\ge\sigma^2+\min_{1\le i\le n}\rho_{\min}(Z_i \Psi_0 Z_i^T),\label{c2eqn:eigmm1}\\
	\max_{{1\le i\le n}}\rho_{\max}(\sigma^2 I_{m_i}+Z_i \Psi_0 Z_i^T)&\le \sigma^2 +\rho_{\max}(\Psi_0)\max_{1\le i\le n}\lVert Z_i\rVert_{\rm sp}^2.\label{c2eqn:eigmm2}
	\end{align}
	Since $Z_i \Psi_0 Z_i^T$ is nonnegative definite, the right hand side of \eqref{c2eqn:eigmm1} is further bounded below by $\sigma^2$, while the right hand side of \eqref{c2eqn:eigmm2} is bounded. The condition \ref{c2con:trueparcon2} is thus satisfied. 
	\item {\it Verification of \rm\ref{c2con:thm2mod1}}: For a sufficiently large $M$ and $s_\star =s_0\vee (\log n/\log p)$, define a sieve as ${\cal H}_n=\{\Psi:n^{-M}\le\rho_{\min}(\Sigma)\le\rho_{\max}(\Sigma)\le e^{Ms_\star \log p}\}$, so that the minimum eigenvalue condition \eqref{c2eqn:thm2c1} can be satisfied with $\log\gamma_n\asymp\log n$. 
	Similar to the proof of Theorem~\ref{c2thm:exmmissing}, it can be easily shown that conditions \eqref{c2eqn:thm2c2} and \eqref{c2eqn:thm2c3} are satisfied with $\epsilon_n=\sqrt{(s_\star \log p)/n}$. 
	\item {\it Verification of \rm\ref{c2con:sep}}: The separation condition is trivially satisfied as there is no nuisance mean part.
\end{itemize}
Therefore, the posterior contraction rates for $\theta$ are given by Theorem~\ref{c2thm:thetacon} with $s_\star$ replaced by $s_0$ since $s_0>0$ and $\log n\lesssim \log p$. The contraction rate for $\Sigma$ relative to the Frobenius norm is a direct consequence of \eqref{c2eqn:maxavebomm}. 	The optimal posterior contraction easily follows from Corollary~\ref{c2cor:opt}.
Thus assertions \ref{c2thm:parcorcon}--\ref{c2thm:parcorcon2} hold.

Now, we verify conditions \ref{c2con:thm4asm0}--\ref{c2con:thm4asm2} and \ref{c2con:sepsp} to apply Theorems~\ref{c2thm:bvm}--\ref{c2thm:sel} and Corollaries~\ref{c2cor:strsel}--\ref{c2cor:strbvm}.
\begin{itemize}[leftmargin=0.4cm, itemsep=0.0cm]
	\item {\it Verification of \rm\ref{c2con:thm4asm0}--\ref{c2con:thm4asm1}}: These conditions are trivially satisfied with the zero matrix $H$ since there is no nuisance mean part.
	\item {\it Verification of \rm\ref{c2con:thm4asm2}}: 	For some $C_1>0$, the entropy in \ref{c2con:thm4asm2} is bounded above by a multiple of $\log N(\delta,\{\Sigma:\lVert\Sigma-\Sigma_0\rVert_{\rm F}\le \hat M_2 C_1\epsilon_n\},\lVert\cdot\rVert_{\rm F})\lesssim0\vee\log(3\hat M_2C_1\epsilon_n/\delta)$ by \eqref{c2eqn:maxavebomm}. The expression in \ref{c2con:thm4asm2} is thus bounded by a constant multiple of $s_\star ^5\log^3 p$ by Remark~\ref{rmk:3}. This tends to zero since $s_\star\lesssim s_0$.
	\item {\it Verification of \rm\ref{c2con:sepsp}}: It is easy to see that $d_{B,n}(\eta,\eta_0)\lesssim \lVert \Psi-\Psi_0\rVert_{\rm F}$ since $\max_i\lVert Z_i\rVert_{\rm sp}\lesssim1$. The separability of the space is thus trivial.
\end{itemize}
Hence, under \ref{c2con:combo1r}, Theorem \ref{c2thm:bvm} can be applied to obtain the distributional approximation in \eqref{c2eqn:bvm} with the zero matrix $H$. Under \ref{c2con:combo1r} and \ref{c2con:selpri}, we obtain the no-superset result in \eqref{c2eqn:sel} by Theorem \ref{c2thm:sel}. The strong results in Corollary \ref{c2cor:strsel} and Corollary \ref{c2cor:strbvm} follow explicitly from the beta-min condition \ref{c2con:betamin}. These establish \ref{c2thm:parcorbvm}--\ref{c2thm:parcorbvm3}.

\subsection{Proof of Theorem~\rm\ref{c2thm:exmgraph}}
We verify the conditions for the posterior contraction in Theorem~\ref{c2thm:thetacon}. 
\begin{itemize}[leftmargin=0.4cm, itemsep=0.0cm]
	\item {\it Verification of \rm\ref{c2con:maxfrob}}: Since $\Delta_{\eta,i}=\Omega^{-1}$ for every $i\le n$ and $\Omega_0\in{\cal M}_0^+(cL)$ for some $0<c<1$, $a_n=1$ and $e_n=0$ satisfy \ref{c2con:maxfrob}.
	\item {\it Verification of \rm\ref{c2con:thm1prcon}}: 	Using (i) of Lemma~\ref{c2lmm:2} and the relation $1-x\asymp1-x^{-1}$ as $x\rightarrow 1$, observe that
	$\lVert \Omega^{-1}-\Omega_0^{-1}\rVert_{\rm F}\lesssim\lVert \Omega-\Omega_0\rVert_{\rm F}\lesssim \bar\epsilon_n$ if the right hand side is small enough.
	Thus, there exists a constant $C_1>0$ such that $\{\Omega: \lVert\Omega^{-1}-\Omega_0^{-1}\rVert_{\rm F}\le\bar\epsilon_n\}\supset\{\Omega: \lVert\Omega-\Omega_0\rVert_{\rm F}\le C_1\bar\epsilon_n\}$.
	Furthermore, although the components of $\Omega$ are not a priori independent as the prior is truncated to ${\cal M}_0^+(L)$, the truncation can only increase prior concentration since  $\Omega_0\in{\cal M}_0^+(cL)$ for some $0<c<1$. Hence, for some $C_2>0$, 
	\begin{align*}
	\Pi\left(\lVert\Omega^{-1}-\Omega_0^{-1}\rVert_{\rm F}\le\bar\epsilon_n\right)\ge \Pi\left(\lVert\Omega-\Omega_0\rVert_{\infty}\le C_2\bar\epsilon_n/\overline m\right)\gtrsim \left(\frac{C_2\bar\epsilon_n}{\overline m}\right)^{\overline m +d},
	\end{align*}
	which justifies the choice $\bar\epsilon_n\asymp \sqrt{(\overline m+d)(\log n)/n}$ for \ref{c2con:thm1prcon}. 
	\item {\it Verification of \rm\ref{c2con:trueparcon1}}: The assumption $\lVert\theta_0\rVert_\infty\lesssim \lambda^{-1}\log p$ given in the theorem directly satisfies \ref{c2con:trueparcon1}.
	\item {\it Verification of \rm\ref{c2con:trueparcon2}}: This is trivially met as $\Omega_0\in{\cal M}_0^+(cL)$ for some $0<c<1$.
	\item {\it Verification of \rm\ref{c2con:thm2mod1}}: Note that the minimum eigenvalue condition \eqref{c2eqn:thm2c1} is trivially satisfied with $\gamma_n=1$ since the prior is put on ${\cal M}_0^+(L)$.
	Now, for $\bar r_n=Ms_\star \log p/\log n$ with $s_\star =s_0\vee (n\bar\epsilon_n^2/\log p)$ and sufficiently large $M$, choose a sieve as ${\cal H}_n=\{\Omega\in{\cal M}_0^+(L) : \sum_{j,k}\mathbbm{1}{\{\omega_{jk}\ne0\}}\le \bar r_n\}$, that is,  the maximum number of edges of $\Omega$ does not exceed $\bar r_n$. Then, for $\delta_n=1/6\overline m n^{3/2}$, the entropy in \eqref{c2eqn:thm2c2} is bounded by
	\begin{align*}
	\log N(\delta_n/\overline m,{\cal H}_n,\lVert\cdot\rVert_{\infty})&\le \log\left\{\left(\frac{\overline m L}{\delta_n}\right)^{\overline m+\bar r_n} \binom{\binom{\overline m}{2}}{\bar r_n} \right\}\\
	& \le (\overline m+\bar r_n)\log(\overline m L/\delta_n) + 2\bar r_n\log \overline m,
	\end{align*}
	where in the second term, the factor $(\overline m L/\delta_n)^{\overline m}$ comes from the diagonal elements of $\Omega$, while the rest is from the off-diagonal entries. It is easy to see that the last display is bounded by a multiple of $s_\star \log p$ with chosen $\bar r_n$, and hence the entropy condition in \eqref{c2eqn:thm2c2} is satisfied. Lastly, note that for some $C_3>0$,
	\begin{align*}
	\log\Pi({\cal H}\setminus{\cal H}_n)&=\log\Pi(|\Upsilon|>\bar r_n)\lesssim-\bar r_n\log \bar r_n \le -C_3Ms_\star \log p.
	\end{align*}
	Therefore, condition \eqref{c2eqn:thm2c3} is satisfied with sufficiently large $M$. 
	\item {\it Verification of \rm\ref{c2con:sep}}: The separation condition is trivially met as there is no nuisance mean part.
\end{itemize}
Therefore, we obtain the posterior contraction properties for $\theta$ by Theorem~\ref{c2thm:thetacon}. The theorem also implies that the posterior distribution of $\Omega^{-1}$ contracts to $\Omega_0^{-1}$ at the rate $\epsilon_n=\sqrt{(s_0\log p \vee (\overline m+d)\log n )/n}$ with respect to the Frobenius norm. This is also translated as convergence of $\Omega$ to $\Omega_0$ at the same rate, since we obtain
\begin{align}
\lVert \Omega-\Omega_0 \rVert_{\rm F}^2\lesssim \lVert \Omega^{-1}-\Omega_0^{-1}\rVert_{\rm F}^2\lesssim\epsilon_n^2,
\label{c2eqn:graphrec}
\end{align}
by (i) of Lemma~\ref{c2lmm:2} and the inequality $1-x\asymp1-x^{-1}$ as $x\rightarrow 1$.
The assertion for the optimal posterior contraction is directly justified by Corollary~\ref{c2cor:opt}.
These prove \ref{c2thm:graphcon}--\ref{c2thm:graphcon2}.

Next, we verify conditions \ref{c2con:thm4asm00}--\ref{c2con:sepsp}  to obtain the optimal posterior contraction by applying Theorem~\ref{c2thm:thmmar}.
\begin{itemize}[leftmargin=0.4cm, itemsep=0.0cm]
	\item {\it Verification of \rm\ref{c2con:thm4asm00}--\ref{c2con:thm4asm10}}: These conditions are trivially satisfied with the zero matrix $H$ since there is no nuisance mean part.
	\item {\it Verification of \rm\ref{c2con:thm4asm20}}: 
	Note that by \eqref{c2eqn:graphrec}, there exists a constant $C_4>0$ such that the entropy in \ref{c2con:thm4asm20} is bounded by $\log N(\delta,\{\Omega:\lVert \Omega-\Omega_0 \rVert_{\rm F}\le C_4\bar\epsilon_n\},d_{B,n})$ for every $\delta>0$. Using \eqref{c2eqn:graphrec2},
	the entropy is further bounded by $\log N(C_5\delta,\{\Omega:\lVert \Omega-\Omega_0 \rVert_{\rm F}\le C_4\bar\epsilon_n\},\lVert\cdot\rVert_{\rm F})$ 
	for some $C_5>0$. This is  clearly bounded by a multiple of $0\vee\overline m^2\log(3C_4\bar\epsilon_n/C_5\delta)$, and hence using Remark~\ref{rmk:3} we bound \ref{c2con:thm4asm20} by a multiple of $(\sqrt{\bar s_\star} \vee\overline m)\sqrt{(\bar s_\star \log p)/n}$ which goes to zero by assumption. 
	\item {\it Verification of \rm\ref{c2con:sepsp}}: 
	For every $\Omega_1,\Omega_2\in\widehat{\cal H}_n$, note that 
	\begin{align}
	\lVert \Omega_1^{-1}-\Omega_2^{-1}\rVert_{\rm F}\lesssim \lVert \Omega_1-\Omega_2 \rVert_{\rm F}\lesssim \lVert \Omega_1^{-1}-\Omega_2^{-1}\rVert_{\rm F}\lesssim\epsilon_n,
	\label{c2eqn:graphrec2}
	\end{align}
	using (i) of Lemma~\ref{c2lmm:2} and the inequality $1-x\asymp1-x^{-1}$ as $x\rightarrow 1$ again.
	By the first inequality, it suffices to show that $\cal H$ is separable metric space with the Frobenius norm. This is trivial as the parameter space is Euclidean. 
\end{itemize}
Hence, under condition \ref{c2con:combo1}, Theorem~\ref{c2thm:thmmar} verifies \ref{c2thm:graphcon3}.

Now, we verify conditions \ref{c2con:thm4asm0}--\ref{c2con:thm4asm2} to apply Theorems~\ref{c2thm:bvm}--\ref{c2thm:sel} and Corollaries~\ref{c2cor:strsel}--\ref{c2cor:strbvm}.
\begin{itemize}[leftmargin=0.4cm, itemsep=0.0cm]
	\item {\it Verification of \rm\ref{c2con:thm4asm0}--\ref{c2con:thm4asm1}}: These are trivially satisfied for the same reason as \ref{c2con:thm4asm00}--\ref{c2con:thm4asm10}.
	\item {\it Verification of \rm\ref{c2con:thm4asm2}}: 
	Similar to the verification of \ref{c2con:thm4asm20},
	the entropy in \ref{c2con:thm4asm2} is bounded by a multiple of $0\vee\overline m^2\log(3C_6\epsilon_n/\delta)$ for some $C_6>0$. Hence using Remark~\ref{rmk:3} we bound \ref{c2con:thm4asm2} by a multiple of $(s_\star \vee\overline m)\sqrt{(s_\star \log p)^3/n}$ which goes to zero by assumption. 
\end{itemize}
Therefore, under \ref{c2con:combo1r}, we obtain the distributional approximation in \eqref{c2eqn:bvm} with the zero matrix $H$ by Theorem \ref{c2thm:bvm}. Under \ref{c2con:combo1r} and \ref{c2con:selpri}, the no-superset result in \eqref{c2eqn:sel} holds by Theorem \ref{c2thm:sel}. Lastly, we obtain the strong results in Corollary \ref{c2cor:strsel} and Corollary \ref{c2cor:strbvm} if the beta-min condition \ref{c2con:betamin} is also met. These prove \ref{c2thm:graphbvm}--\ref{c2thm:graphbvm3}.

\subsection{Proof of Theorem~\rm\ref{c2thm:exmhet}}
To verify the conditions for Theorem~\ref{c2thm:thetacon}, we will use the following properties of $B$-splines.

For any $f\in \mathfrak{C}^\alpha[0,1]$,   there exists $\beta_\ast\in \mathbb{R}^J$ with  $\lVert\beta_\ast\rVert_\infty < \lVert f\rVert_{\mathfrak{C}^\alpha}$ such that
\begin{align}
\lVert \beta_\ast^T B_J- f\rVert_\infty \lesssim J^{-\alpha}\lVert f\rVert_{\mathfrak{C}^\alpha},
\label{c2eqn:bsperr}
\end{align}
by the well-known approximation theory of B-splines \citep[page 170]{de1978practical}. Writing $f_\beta=\beta^T B_J$, this gives
\begin{align}
\lVert f_\beta-f\rVert_{2,n}\le\lVert f_\beta-f\rVert_\infty\lesssim J^{-\alpha}\lVert f\rVert_{\mathfrak{C}^\alpha}+\lVert f_\beta-f_{\beta_\ast}\rVert_\infty.
\label{c2eqn:bspapp}
\end{align}
We also use the following inequalities: for every $\beta\in\mathbb{R}^J$,
\begin{align}
\lVert\beta\rVert_\infty \lesssim \lVert f_{\beta}\rVert_\infty\le \lVert\beta\rVert_\infty,\quad \lVert\beta\rVert_2 \lesssim \sqrt{J}\lVert f_{\beta}\rVert_{2,n}\lesssim \lVert\beta\rVert_2.
\label{c2eqn:spbound}
\end{align}
See Lemma E.6 of \citet{ghosal2017fundamentals} for proofs with respect to $L_\infty$- and $L_2$-norms. Hence the first relation can be formally justified. For the second relation with respect to the empirical $L_2$-norm, we assume that $z_i$ are sufficiently regularly distributed as in (7.12) of \citet{ghosal2007convergence}.
\begin{itemize}[leftmargin=0.4cm, itemsep=0.0cm]
	\item {\it Verification of \rm\ref{c2con:maxfrob}}: 	If $v_0$ is strictly positive on $[0,1]$, then $v_0$ satisfies the same approximation rule in \eqref{c2eqn:bsperr} for some $\beta_\ast\in(0,\infty)^J$ with $\lVert\beta_\ast\rVert_\infty<\lVert v_0\rVert_{\mathfrak{C}^\alpha}$ (see Lemma E.5 of \citet{ghosal2017fundamentals}). Therefore the approximation in \eqref{c2eqn:bspapp} also holds for $v_0$ even if $\beta$ is restricted to have positive entries only, and thus by \eqref{c2eqn:bsperr} and \eqref{c2eqn:spbound},
	\begin{align*}
	\lVert v_{\beta_\ast}-v_0 \rVert_\infty &\lesssim J^{-\alpha},  \quad\text{for some $\beta_\ast\in(0,\infty)^J$},\\
	\lVert v_{\beta_1}-v_{\beta_2} \rVert_\infty &\lesssim  \sqrt{J} \lVert v_{\beta_1}-v_{\beta_2} \rVert_{2,n}, \quad\beta_1,\beta_2\in(0,\infty)^J,
	\end{align*}	
	which tells us that we have $a_n\asymp J$ and $e_n\asymp J^{1-2\alpha}$ for \ref{c2con:maxfrob}. 
	\item {\it Verification of \rm\ref{c2con:thm1prcon}}: Note that if $J^{-\alpha}\lesssim \bar\epsilon_n$, it follows that for some $C_1>0$,
	\begin{align*}
	\log\Pi (\beta:\lVert v_\beta-v_0\rVert_{2,n}\le \bar\epsilon_n )&\ge\log\Pi (\beta:\lVert \beta-\beta_\ast\rVert_\infty\le C_1\bar\epsilon_n )\gtrsim J\log \bar\epsilon_n.
	\end{align*}
	This implies that condition \ref{c2con:thm1prcon} is satisfied with $\bar\epsilon_n=\sqrt{(J\log n)/n}$. 
	\item {\it Verification of \rm\ref{c2con:trueparcon1}}: The assumption $\lVert\theta_0\rVert_\infty\lesssim \lambda^{-1}\log p$ given in the theorem directly satisfies \ref{c2con:trueparcon1}.
	\item {\it Verification of \rm\ref{c2con:trueparcon2}}: Since $v_0$ is strictly positive on $[0,1]$ and belongs to a fixed multiple of the unit ball of $\mathfrak{C}^{\alpha}[0,1]$, we have that
	\begin{align*}
	1\lesssim\inf_{z\in[0,1]} v_0(z)\le\sup_{z\in[0,1]} v_0(z)\lesssim 1.
	\end{align*}
	The condition \ref{c2con:trueparcon2} is thus satisfied.
	\item {\it Verification of \rm\ref{c2con:thm2mod1}}:  
	For a sufficiently large $M$, choose a sieve as ${\cal H}_n=\prod_{j=1}^J\{\beta_j :n^{-M}\le \beta_j \le n^M\}$. Then the minimum eigenvalue condition \eqref{c2eqn:thm2c1} is satisfied with $\log\gamma_n\asymp\log n$ because for every $i\le n$,
	\begin{align*}
	\inf_{\beta\in{\cal H}_n} v_\beta(z_i)=\inf_{\beta\in{\cal H}_n}\sum_{j=1}^J B_{J,j}(z_i)\beta_j \ge\inf_{\beta\in{\cal H}_n}\min_{1\le j\le J}\beta_j\sum_{j=1}^J B_{J,j}(z_i)\ge n^{-M},
	\end{align*} where $B_{J,j}$ and $\beta_j$ denote the $j$th components of $B_J$ and $\beta$, respectively.
	To check the entropy condition in \eqref{c2eqn:thm2c2}, note that for every $\eta_1,\eta_2\in{\cal H}_n$, we have $d_n(\eta_1,\eta_2)\lesssim \lVert\beta_1-\beta_2\rVert_\infty$ by \eqref{c2eqn:spbound}. Hence, for some $C_2>0$, the entropy in \eqref{c2eqn:thm2c2} is bounded above by a multiple of
	\begin{align*}
	&\log N\left(\frac{1}{C_2\overline m n^{M+3/2}},\{\beta:\lVert\beta\rVert_\infty\le n^M\},\lVert\cdot\rVert_\infty\right)\lesssim J\log n .
	\end{align*}
	The condition \eqref{c2eqn:thm2c3} holds since an inverse Gaussian prior on each $\beta_j$ produces $\Pi({\cal H}\setminus{\cal H}_{n})\lesssim J e^{-C_3n^M}$ for some constant $C_3$, by its exponentially small bounds for tail probabilities on both sides.
	By matching $J^{-\alpha}\asymp\bar\epsilon_n$ and $n\bar\epsilon_n^2 \asymp J \log n$, we obtain $J\asymp (n/\log n)^{1/(2\alpha+1)}$ and $\bar\epsilon_n=(\log n/n)^{\alpha/(2\alpha+1)}$.
	Note that the conditions $a_n\epsilon_n^2\rightarrow 0$ and $e_n\rightarrow0$ hold only if $\alpha>1/2$.
	\item {\it Verification of \rm\ref{c2con:sep}}: The separation condition  holds as there is no additional mean part.
\end{itemize}	
Hence, we obtain the posterior contraction rates for $\theta$ by Theorem~\ref{c2thm:thetacon}. The contraction rate for $v$ is also obtained by the same theorem. The assertion for the optimal posterior contraction is directly justified by Corollary~\ref{c2cor:opt}. Hence we have verified \ref{c2thm:hetcon}--\ref{c2thm:hetcon2}.

Now, we verify \ref{c2con:thm4asm00}--\ref{c2con:sepsp} for the optimal posterior contraction in Theorem~\ref{c2thm:thmmar}.
\begin{itemize}[leftmargin=0.4cm, itemsep=0.0cm]
	\item {\it Verification of \rm\ref{c2con:thm4asm00}--\ref{c2con:thm4asm10}}: These conditions are trivially satisfied as there is no nuisance mean part.
	\item {\it Verification of \rm\ref{c2con:thm4asm20}}: Note that by the inequality $\lVert v_\beta-v_0\rVert_{2,n}\lesssim \lVert v_\beta-v_{\beta_\ast}\rVert_{2,n}+\bar\epsilon_n$,
	the entropy in the integrand is bounded by
	\begin{align*}
	\log N\left(\delta\sqrt{J}, \left\{\beta:\lVert \beta-\beta_\ast\rVert_2\le C_4\sqrt{J}\bar\epsilon_n\right\},\lVert\cdot\rVert_2\right)\lesssim 0\vee J\log \left(\frac{3C_4\bar\epsilon_n}{\delta}\right),
	\end{align*}
	for some $C_4>0$. Thus, the second term of \ref{c2con:thm4asm20} is bounded by $J\bar\epsilon_n$ by Remark~\ref{rmk:3}, while the first term is bounded by $\sqrt{J \bar s_\star^2 (\log p) /n}$.
	Since $\bar s_\star = (J\log n)/\log p\lesssim J$, \ref{c2con:thm4asm20} is bounded by
	$J \bar\epsilon_n = (n/\log n)^{(1-\alpha)/(2\alpha+1)}$, which tends to zero as $\alpha>1$.
	\item {\it Verification of \rm\ref{c2con:sepsp}}: For every $v_{\beta_1},v_{\beta_2}\in\widehat{\cal H}_n$,
	note that $d_{B,n}(\eta_1,\eta_2)=\lVert v_{\beta_1}-v_{\beta_2}\rVert_{2,n}\lesssim \lVert\beta_1-\beta_2\rVert_2$ by \eqref{c2eqn:spbound}. Since we put a prior for $v$ using the B-splines through a Euclidean parameter $\beta$, the separability is trivially satisfied.
\end{itemize}
Therefore, since \ref{c2con:combo1} is satisfied the assumption, assertion \ref{c2thm:hetcon3} holds by Theorem~\ref{c2thm:thmmar}.

Next, we verify conditions \ref{c2con:thm4asm0}--\ref{c2con:thm4asm2} to apply Theorems~\ref{c2thm:bvm}--\ref{c2thm:sel} and Corollaries~\ref{c2cor:strsel}--\ref{c2cor:strbvm}.
\begin{itemize}[leftmargin=0.4cm, itemsep=0.0cm]
	\item {\it Verification of \rm\ref{c2con:thm4asm0}--\ref{c2con:thm4asm1}}: These are trivially satisfied for the same reason as before.
	\item {\it Verification of \rm\ref{c2con:thm4asm2}}: Similar to the verification of \ref{c2con:thm4asm20}, the entropy of interest is bounded by a constant multiple of $0\vee J\log ({3 C_5\epsilon_n}/{\delta})$ for some $C_5>0$.
	Thus, \ref{c2con:thm4asm2} is bounded above by a multiple of $\{(s_\star ^2\vee J)J(s_\star \log p)^3 /n\}^{1/2}$ by Remark~\ref{rmk:3}, and hence goes to zero by the assumption.
	The condition $\alpha>2$ is seen to be necessary by the inequality 
	\begin{align*}
	(s_\star ^2\vee J)J(s_\star \log p)^3 /n\ge J^2 n^2\bar\epsilon_n^6 =n^{2(-\alpha+2)/(2\alpha+1)}{\log n}^{2(3\alpha-1)/(2\alpha+1)}.
	\end{align*}
\end{itemize}
Under \ref{c2con:combo1r}, the distributional approximation in \eqref{c2eqn:bvm} holds with the zero matrix $H$ by Theorem \ref{c2thm:bvm}. Under \ref{c2con:combo1r} and \ref{c2con:selpri}, the no-superset result in \eqref{c2eqn:sel} holds by Theorem \ref{c2thm:sel}. We also obtain the strong results in Corollary \ref{c2cor:strsel} and Corollary \ref{c2cor:strbvm} if the beta-min condition \ref{c2con:betamin} is also met. These prove \ref{c2thm:hetbvm}--\ref{c2thm:hetbvm3}.

\subsection{Proof of Theorem~\rm\ref{c2thm:exmplm}}
We verify the conditions for the posterior contraction in Theorem~\ref{c2thm:thetacon}. 
	\begin{itemize}[leftmargin=0.4cm, itemsep=0.0cm]
		\item {\it Verification of \rm\ref{c2con:maxfrob}}: Since $\Delta_{\eta,i}=\sigma^2$ for every $i\le n$ and $\sigma_0^2$ belongs to the support of the prior, we have $a_n=1$ and $e_n=0$.
		\item {\it Verification of \rm\ref{c2con:thm1prcon}}: Note that we write $d_n^2(\eta,\eta_0) = \lvert\sigma^2-\sigma_0^2\rvert^2+\lVert g_\beta-g_0\rVert_{2,n}^2$.
		To verify the prior concentration condition, observe that
		\begin{align*}
		&\log\Pi\left(\eta\in{\cal H}:d_n(\eta,\eta_0)\le \bar\epsilon_n\right)\\
		&\quad\ge\log\Pi\left(\beta:\lVert g_\beta-g_0\rVert_{2,n}\le\frac{\bar\epsilon_n}{\sqrt{2}}\right)+\log\Pi\left(\sigma:|\sigma^2-\sigma_0^2|\le\frac{\bar\epsilon_n}{\sqrt{2}}\right),
		\end{align*}
		where the second term on the right hand side is trivially bounded below by a constant multiple of $\log \bar\epsilon_n$. Using \eqref{c2eqn:bsperr}--\eqref{c2eqn:spbound}, it is easy to see that if $J^{-\alpha}\lesssim \bar\epsilon_n$, 
		\begin{align*}
		\log\Pi\left(\beta:\lVert g_\beta-g_0\rVert_{2,n}\le\frac{\bar\epsilon_n}{\sqrt{2}}\right)\ge\log\Pi (\beta:\lVert \beta-\beta_\ast\rVert_\infty\le C_1\bar\epsilon_n ) \gtrsim J\log \bar\epsilon_n,
		\end{align*}
	for some $C_1>0$.
	Since $\bar\alpha\le \alpha$, this implies that \ref{c2con:thm1prcon} is satisfied with $\bar\epsilon_n=\sqrt{(J\log n)/n}$.
		\item {\it Verification of \rm\ref{c2con:trueparcon1}}: The assumption $\lVert\theta_0\rVert_\infty\lesssim \lambda^{-1}\log p$ given in the theorem directly satisfies the condition.
		\item {\it Verification of \rm\ref{c2con:trueparcon2}}: This is directly satisfied by $\sigma_0^2\asymp1$.
		\item {\it Verification of \rm\ref{c2con:thm2mod1}}: For a sufficiently large constant $M$ and  $s_\star =s_0\vee (J\log n/\log p)$, choose
		${\cal H}_n=\{g_\beta:\lVert\beta\rVert_\infty\le n^M\}\times\{\sigma:n^{-M}\le\sigma^2\le e^{Ms_\star \log p}\}$, from which the minimum eigenvalue condition \eqref{c2eqn:thm2c1} is directly satisfied with $\log\gamma_n\asymp\log n$. 
		To check the entropy condition in \eqref{c2eqn:thm2c2}, note that for every $\eta_1,\eta_2\in{\cal H}_n$, we have $d_n^2(\eta_1,\eta_2)\lesssim \lVert\beta_1-\beta_2\rVert_\infty^2+\lvert\sigma_1^2-\sigma_2^2\rvert^2$ by \eqref{c2eqn:spbound}. Hence, for some $C_3>0$, the entropy in \eqref{c2eqn:thm2c2} is bounded above by a multiple of
		\begin{align*}
		&\log N\left(\frac{1}{C_3\overline m n^{M+3/2}},\{\beta:\lVert\beta\rVert_\infty\le n^M\},\lVert\cdot\rVert_\infty\right)\\
		&\quad+\log N\left(\frac{1}{C_3\overline m n^{M+3/2}},\{\sigma:\sigma^2\le e^{Ms_\star \log p}\},\lvert\cdot\rvert\right).
		\end{align*}
		The display is further bounded by a multiple of $J\log n+s_\star \log p$, and hence \eqref{c2eqn:thm2c2} is satisfied with $\epsilon_n=\sqrt{(s_\star \log p)/n}$. Using the tail bounds of normal and inverse gamma distributions, condition \eqref{c2eqn:thm2c3} is also satisfied. 
		\item {\it Verification of \rm\ref{c2con:sep}}: The separation condition holds by Remark~\ref{rmk:1} as we have $d_{A,n}(\eta_\ast,\eta_0)=\lVert g_{\beta_\ast}-g_0 \rVert_{2,n}\lesssim \bar\epsilon_n$ for $\eta_\ast=(g_{\beta_\ast},\sigma_0^2)$ in view of \eqref{c2eqn:bsperr}.
	\end{itemize}
	Therefore, the contraction rates for $\theta$ are given by Theorem~\ref{c2thm:thetacon}. The rate for $g$ is also obtained by the same theorem. The assertion for the optimal posterior contraction is directly justified by Corollary~\ref{c2cor:opt}. We thus see \ref{c2thm:hetcon}--\ref{c2thm:hetcon2} hold.

Now, we verify \ref{c2con:thm4asm00}--\ref{c2con:sepsp} for Theorem~\ref{c2thm:thmmar}.
\begin{itemize}[leftmargin=0.4cm, itemsep=0.0cm]
	\item {\it Verification of \rm\ref{c2con:thm4asm00}}: Observe that the left hand side of the first line of \ref{c2con:thm4asm00} is equal to
	\begin{align*}
	\frac{1}{(s_0\vee 1)\log p}\lVert \tilde \xi_{\eta_0}-H \tilde \xi_{\eta_0} \rVert_2^2&= \frac{n}{\sigma_0^2 (s_0\vee 1)\log p}\lVert g_0-\hat\beta_J^T B_J \rVert_{2,n}^2,
	\end{align*}
	where $\hat\beta_J=( W_J^T  W_J)^{-1}  W_J^T (g_0(z_1),\dots,g_0(z_n))^T$ is the the least squares solution. Since $\hat\beta_J$ is the solution minimizing $\lVert g_0-\hat\beta_J^T B_J \rVert_{2,n}^2$, for some $\beta_\ast\in \mathbb{R}^J$, the last display is bounded above by
	\begin{align}
	\frac{n}{\sigma_0^2\log p} \lVert g_0-\beta_\ast^T B_J\rVert_\infty^2\lesssim \frac{n}{J^{2\alpha}\log p},
	\label{c2eqn:subopt}
	\end{align}
	by \eqref{c2eqn:bsperr}, where $s_0\vee 1$ is replaced by $1$ as $s_0$ is unknown.
	Plugging in $J\asymp (n/\log n)^{1/(2\bar\alpha+1)}$, it is easy to see that the right hand side of \eqref{c2eqn:subopt} is the same order of
	$(\log n)^{2\alpha/(2\bar\alpha+1)}n^{(-2\alpha+2\bar\alpha+1)/(2\bar\alpha+1)}/\log p$. This tends to zero by the given boundedness assumption. The necessary condition $\bar\alpha<\alpha$ is implied by this, because $\log p=o(n)$.
	 The second condition of \ref{c2con:thm4asm00} is satisfied by Remark~\ref{rmk:2}.
	\item {\it Verification of \rm\ref{c2con:thm4asm10}}: Let $\tilde\eta_n(\theta,\eta)=(g_\beta(\cdot)+ B_J^T(\cdot)(W_J^T W_J)^{-1} W_J^T X(\theta-\theta_0),\sigma^2)$ for a given $\theta$, where $\eta=(g_\beta(\cdot),\sigma^2)$. This setting satisfies  $\Phi(\tilde\eta_n(\theta,\eta))=(\tilde\xi_\eta+H\tilde X(\theta-\theta_0),\tilde\Delta_\eta)$. 
	Since each entry of $\beta$ has the standard normal prior, $g_\beta(\cdot)$ is a zero mean Gaussian process with the covariance kernel $K(t_1,t_2)=B_J(t_1)^T B_J(t_2)$, and thus its reproducing kernel Hilbert space (RKHS) $\mathbb{K}$ is the set of all functions of the form $\sum_k \zeta_k B_J(t_k)^T B_J(\cdot)$ with coefficients $\zeta_k$, $k\in\{1,2,\dots\}$. It is easy to see that the shift $(\theta-\theta_0)^T X^T W_J (W_J^T W_J)^{-1}  B_J(\cdot)$ is in the RKHS $\mathbb{K}$ since it is expressed as $(\theta-\theta_0)^T X^T W_J(W_J^T W_J)^{-1}\tilde W_J^{-1}\tilde W_J B_J(\cdot)$ using an invertible matrix $\tilde W_J\in\mathbb{R}^{ J\times J}$ with rows $B_J(t_k)$ evaluated by some $t_k$, $k=1,\dots,J$.
	Hence, by the Cameron-Martin theorem, for $\nu=(\nu_1,\dots,\nu_J)^T= (\tilde W_J^T)^{-1}(W_J^T W_J)^{-1}W_J^T X(\theta-\theta_0)$ and $\lVert\cdot\rVert_{\mathbb{K}}$ the RKHS norm, we see that
	\begin{align*}
	\log \frac{d\Pi_{n,\theta}}{d\Pi_{n,\theta_0}}(\eta)&=\sum_{k=1}^J\nu_k g_\beta(t_k)-\frac{1}{2} \lVert \nu^T \tilde W_J B_J \rVert_{\mathbb{K}}^2 =\nu^T \tilde W_J \beta -\frac{1}{2}\lVert \tilde W_J^T \nu\rVert_2^2,
	\end{align*}
	almost surely.
	This gives that
	\begin{align}
	\begin{split}
	\left\lvert\log \frac{d\Pi_{n,\theta}}{d\Pi_{n,\theta_0}}(\eta)\right\rvert
&\lesssim\lVert\beta\rVert_2\lVert(W_J^T W_J)^{-1} W_J^T  X(\theta-\theta_0)\rVert_2 \\
&\quad+\lVert(W_J^T W_J)^{-1} W_J^T  X(\theta-\theta_0)\rVert_2^2.
	\end{split}
	\label{c2eqn:aaqwe}
	\end{align}
	Note that we have
	\begin{align*}
	\sup_{\eta\in\widetilde{\cal H}_n}\lVert\beta\rVert_2&\le\sup_{\eta\in\widetilde{\cal H}_n}\lVert\beta-\beta_\ast\rVert_2+\lVert\beta_\ast\rVert_2\\
	&\lesssim \sqrt J\sup_{\eta\in\widetilde{\cal H}_n}\lVert g_\beta-g_{\beta_\ast}\rVert_{2,n}+1 \lesssim \sqrt{J}\bar\epsilon_n+1,
	\end{align*} and
	\begin{align*}
	\sup_{\theta\in\widetilde\Theta_n}\lVert( W_J^T W_J)^{-1}  W_J^T  X(\theta-\theta_0)\rVert_2&\lesssim \frac{\lVert W_J\rVert_{\rm sp}\sup_{\theta\in\widetilde\Theta_n}\lVert X(\theta-\theta_0)\rVert_2}{\rho_{\min}(W_J^T W_J)}\lesssim \sqrt{J}\bar \epsilon_n,
	\end{align*}
	using \eqref{c2eqn:spbound}. Since $\sqrt{J}\bar\epsilon_n$ is bounded due to $\bar\alpha\ge1/2$, \eqref{c2eqn:aaqwe} is bounded.
	\item {\it Verification of \rm\ref{c2con:thm4asm20}}: Since the entropy in the integral in \ref{c2con:thm4asm20} is bounded above by a multiple of $0\vee\log(3\tilde M_2\bar\epsilon_n/\delta)$ for every $\delta>0$, the second term of \ref{c2con:thm4asm20} is bounded by a constant multiple of $\bar\epsilon_n$ due to Remark~\ref{rmk:3}. The first term is $\bar\epsilon_n^2\sqrt{n/\log p}=(\log n)^{2\bar\alpha/(2\bar\alpha+1)}n^{(-\bar\alpha+1/2)/(2\bar\alpha+1)}/\sqrt{\log p}$ that tends to zero by the boundedness assumption.
	\item {\it Verification of \rm\ref{c2con:sepsp}}: Since we have $d_{B,n}(\eta_1,\eta_2)=|\sigma_1^2-\sigma_2^2|$ for every $\sigma_1^2,\sigma_2^2\in(0,\infty)$ and the parameter space of $\sigma^2$ is Euclidean, the condition is trivially satisfied.
\end{itemize}
Therefore, assertion \ref{c2thm:plmcon3} holds by Theorem~\ref{c2thm:thmmar} since \ref{c2con:combo1} is also satisfied by the given assumption.

Lastly, we verify conditions \ref{c2con:thm4asm0}--\ref{c2con:thm4asm2} to apply Theorems~\ref{c2thm:bvm}--\ref{c2thm:sel} and Corollaries~\ref{c2cor:strsel}--\ref{c2cor:strbvm}.
\begin{itemize}[leftmargin=0.4cm, itemsep=0.0cm]
	\item {\it Verification of \rm\ref{c2con:thm4asm0}}: Similar to the verification of \ref{c2con:thm4asm00}, the first line of \ref{c2con:thm4asm0} is equal to
	\begin{align*}
	s_\star^2\log p\lVert \tilde \xi_{\eta_0}-H \tilde \xi_{\eta_0} \rVert_2^2&\lesssim \frac{ns_\star^2 \log p}{J^{2\alpha}}.
	\end{align*}
	Plugging in $J\asymp (n/\log n)^{1/(2\bar\alpha+1)}$, it is easy to see that this tends to zero by the given boundedness condition, which requires that $\bar\alpha<\alpha-1/2$.
	\item {\it Verification of \rm\ref{c2con:thm4asm1}}: Similar to the verification of \ref{c2con:thm4asm10}, we now have
	\begin{align*}
		\sup_{\eta\in\widehat{\cal H}_n}\lVert\beta\rVert_2& \lesssim s_\star \sqrt{(J\log p)/n}+1,\\
	\sup_{\theta\in\widehat\Theta_n}\lVert( W_J^T W_J)^{-1}  W_J^T  X(\theta-\theta_0)\rVert_2&\lesssim \frac{\lVert W_J\rVert_{\rm sp}\sup_{\theta\in\widehat\Theta_n}\lVert X(\theta-\theta_0)\rVert_2}{\rho_{\min}(W_J^T W_J)} \\
	&\lesssim s_\star \sqrt{(J\log p)/n}.
	\end{align*}
	Since $J\log n =n\bar\epsilon_n^2 \le s_\star\log p$, \eqref{c2eqn:aaqwe} tends to zero since $s_\star ^5\log^3 p = o(n)$.
	\item {\it Verification of \rm\ref{c2con:thm4asm2}}: By the similar calculations as before, we see that \ref{c2con:thm4asm2} is bounded by $(s_\star ^5\log^3 p/n)^{1/2}$ which tends to zero.
	The condition $\bar\alpha>1$ is necessary since $(s_\star^5\log^3 p)/n\ge n^2\bar\epsilon_n^6 = (\log n)^{6\bar\alpha/(2\bar\alpha+1)} n^{-2(\bar\alpha-1)/(2\bar\alpha+1)}$.
\end{itemize}
	Therefore, under \ref{c2con:combo1r}, we have the distributional approximation in \eqref{c2eqn:bvm} by Theorem~\ref{c2thm:bvm}. Under \ref{c2con:combo1r} and \ref{c2con:selpri}, Theorem~\ref{c2thm:sel} implies that the no-superset result in \eqref{c2eqn:sel} holds. The stronger assertions in \eqref{c2eqn:strsel} and \eqref{c2eqn:strbvm} are explicitly derived from Corollary~\ref{c2cor:strsel} and Corollary~\ref{c2cor:strbvm} if the beta-min condition \ref{c2con:betamin} is also met.
	

\section{Auxiliary results}
Here we provide some auxiliary results used to prove the main results.
\begin{lemma}
	Let $p_k$ be the density of ${\rm N}_r(\mu_k,\Sigma_k)$ for $k=1,2$. Then,
	\begin{align*}
	K(p_1,p_2)=&\frac{1}{2}\left\{\log\frac{\det\Sigma_2}{\det \Sigma_1}+{\rm tr}(\Sigma_1\Sigma_2^{-1})-r +\lVert\Sigma_2^{-1/2}(\mu_1-\mu_2)\rVert_2^2\right\},\\
	V(p_1,p_2)=&\frac{1}{2}\Big\{{\rm tr}(\Sigma_1\Sigma_2^{-1}\Sigma_1\Sigma_2^{-1})-2{\rm tr}(\Sigma_1\Sigma_2^{-1})+r\Big\}+\lVert\Sigma_1^{1/2}\Sigma_2^{-1}(\mu_1-\mu_2)\rVert_2^2.
	\end{align*}
	\label{c2lmm:kullback}
\end{lemma}

\begin{proof}
	Let $Z=\Sigma_1^{-1/2}(X-\mu_1)\sim {\rm N}_r(0,I)$ for $X\sim p_1$ and $A=\Sigma_1^{1/2}\Sigma_2^{-1}\Sigma_1^{1/2}$. Then by direct calculations, we have
	\begin{align*}
	K(p_1,p_2)&=\mathbb{E}_{p_1} \left\{\log \frac{p_1}{p_2}(X)\right\}\\
	&=\frac{1}{2} \left\{\log\frac{\det \Sigma_2}{\det \Sigma_1} +\mathbb{E}_{p_1}  Z^TAZ-r+(\mu_1-\mu_2)^T \Sigma_2^{-1}(\mu_1-\mu_2)\right\},
	\end{align*}
	which verifies the first assertion because $\mathbb{E}_{p_1}  Z^TAZ={\rm tr}A$. After some algebra, we also obtain
	\begin{align*}
	V(p_1,p_2)&=\mathbb{E}_{p_1} \left\{\log \frac{p_1}{p_2}(X)-K(p_1,p_2)\right\}^2\\
	&=\frac{1}{4}\mathbb{E}_{p_1} \left\{-Z^TZ +Z^T A Z+2(\mu_1-\mu_2)^T \Sigma_2^{-1}\Sigma_1^{1/2} Z-{\rm tr}(A)+r\right\}^2.
	\end{align*}
	The rightmost side involves forms of $\mathbb{E}_{p_1}(Z Z^T Q_1 Z)$ and $\mathbb{E}_{p_1}(Z^T Q_1 Z Z^T Q_2 Z)$ for two positive definite matrices $Q_1$ and $Q_2$. 
	It is easy to see that the former is zero, while it can be shown the latter equals $2{\rm tr}(Q_1 Q_2)+{\rm tr}(Q_1){\rm tr}(Q_2)$; for example, see Lemma 6.2 of \citet{magnus1978moments}. Plugging in this for the expected values of the products of quadratic forms, it is easy (but tedious) to verify the second assertion.
\end{proof}	

\begin{lemma}
	For $r\times r$ positive definite matrices $\Sigma_1$ and $\Sigma_2$, let $d_1,\dots,d_r$ be the eigenvalues of $\Sigma_2^{1/2}\Sigma_1^{-1}\Sigma_2^{1/2}$. Then the following assertions hold:
	\begin{enumerate}[label={\rm (\roman*)}]
		\item \label{lmm2:1} $\rho_{\max}^{-2}(\Sigma_2)\lVert\Sigma_1-\Sigma_2\rVert_{\rm F}^2\le\sum_{k=1}^r (d_k^{-1}-1)^2\le \rho_{\min}^{-2}(\Sigma_2)\lVert\Sigma_1-\Sigma_2\rVert_{\rm F}^2$,
		\item \label{lmm2:2} $\max_k|d_k-1|$ can be made arbitrarily small if $g^2(\Sigma_1,\Sigma_2)$ is chosen sufficiently small, where $g$ is defined in \eqref{c2eqn:defg}.
	\end{enumerate}
	\label{c2lmm:2}
\end{lemma}

\begin{proof}
	Let $A=\Sigma_2^{-1/2}\Sigma_1\Sigma_2^{-1/2}$.
	Since the eigenvalues of $A-I_r$ are $d_1^{-1}-1,\dots,d_r^{-1}-1$, we can see that $	\lVert\Sigma_1-\Sigma_2\rVert_{\rm F}^2$ is equal to
	\begin{align*}
\lVert\Sigma_2^{1/2}(A-I_r)\Sigma_2^{1/2}\rVert_{\rm F}^2 \le\rho_{\max}^2(\Sigma_2)\lVert A-I_r\rVert_{\rm F}^2 =
	\rho_{\max}^2(\Sigma_2)\sum_{k=1}^r(d_k^{-1}-1)^2.
	\end{align*}
	Conversely, using the sub-multiplicative property of the Frobenius norm, $\lVert BC \rVert_{\rm F}\le\lVert B\rVert_{\rm sp}\lVert C \rVert_{\rm F}$, it can be seen that $	\sum_{k=1}^r(d_k^{-1}-1)^2$ is equal to
	\begin{align*}
\lVert A-I_r\rVert_{\rm F}^2 =\lVert\Sigma_2^{-1/2}(\Sigma_1-\Sigma_2)\Sigma_2^{-1/2}\rVert_{\rm F}^2 \le\rho_{\max}^2(\Sigma_2^{-1})\lVert\Sigma_1-\Sigma_2\rVert_{\rm F}^2.
	\end{align*}
	These verify \ref{lmm2:1}.
	Now, note that by direct calculations,
	\begin{align*}
	\frac{(\det\Sigma_1)^{1/4}(\det\Sigma_2)^{1/4}}{\det((\Sigma_1+\Sigma_2)/2)^{1/2}}&=\left\{\frac{1}{2^r}\det(A^{1/2}+A^{-1/2})\right\}^{-1/2} \\
	&=\left\{\prod_{k=1}^r\frac{1}{2}(d_k^{1/2}+d_k^{-1/2})\right\}^{-1/2}.
	\end{align*}
	Hence, $g^2(\Sigma_1,\Sigma_2)<\delta$ for a sufficiently small $\delta>0$ implies that
	\begin{align*}
	\prod_{k=1}^r \frac{1}{2}(d_k^{1/2}+d_k^{-1/2})<(1-\delta^2/2)^{-2}.
	\end{align*}
	Since every term in the product of the last display is greater than or equal to 1, we have $(d_k^{1/2}+d_k^{-1/2})/2<(1-\delta^2/2)^{-2}$ for every $k$.
	As a function of $d_k$, $(d_k^{1/2}+d_k^{-1/2})/2$ has the global minimum at $d_k=1$, and hence $\delta$ can be chosen sufficiently small to make $|d_k-1|$ small for every $k=1,\dots, r$, which establishes \ref{lmm2:2}.
\end{proof}


\bibliographystyle{apalike}
\bibliography{ref}

\begin{thebibliography}{}

\bibitem[Atchad{\'e}, 2017]{atchade2017contraction}
Atchad{\'e}, Y.~A. (2017).
\newblock On the contraction properties of some high-dimensional
  quasi-posterior distributions.
\newblock {\em The Annals of Statistics}, 45(5):2248--2273.

\bibitem[Bai et~al., 2020]{bai2019spike}
Bai, R., Moran, G.~E., Antonelli, J., Chen, Y., and Boland, M.~R. (2020).
\newblock Spike-and-slab group lassos for grouped regression and sparse
  generalized additive models.
\newblock {\em Journal of the American Statistical Association}, to appear.

\bibitem[Belitser and Ghosal, 2020]{belitser2017empirical}
Belitser, E. and Ghosal, S. (2020).
\newblock Empirical {B}ayes oracle uncertainty quantification for regression.
\newblock {\em The Annals of Statistics}, 48(6):3113--3137.

\bibitem[Bickel and Kleijn, 2012]{bickel2012semiparametric}
Bickel, P.~J. and Kleijn, B.~J. (2012).
\newblock The semiparametric {B}ernstein--von {M}ises theorem.
\newblock {\em The Annals of Statistics}, 40(1):206--237.

\bibitem[Bondell and Reich, 2012]{bondell2012consistent}
Bondell, H.~D. and Reich, B.~J. (2012).
\newblock Consistent high-dimensional {B}ayesian variable selection via
  penalized credible regions.
\newblock {\em Journal of the American Statistical Association},
  107(500):1610--1624.

\bibitem[Carroll et~al., 2006]{carroll2006measurement}
Carroll, R.~J., Ruppert, D., Crainiceanu, C.~M., and Stefanski, L.~A. (2006).
\newblock {\em Measurement Error in Nonlinear Models: A Modern Perspective}.
\newblock Chapman and Hall/CRC.

\bibitem[Castillo, 2012]{castillo2012semiparametric}
Castillo, I. (2012).
\newblock A semiparametric {B}ernstein--von {M}ises theorem for {G}aussian
  process priors.
\newblock {\em Probability Theory and Related Fields}, 152(1-2):53--99.

\bibitem[Castillo et~al., 2015]{castillo2015bayesian}
Castillo, I., Schmidt-Hieber, J., and van~der Vaart, A. (2015).
\newblock Bayesian linear regression with sparse priors.
\newblock {\em The Annals of Statistics}, 43(5):1986--2018.

\bibitem[Castillo and van~der Vaart, 2012]{castillo2012needles}
Castillo, I. and van~der Vaart, A. (2012).
\newblock Needles and straw in a haystack: Posterior concentration for possibly
  sparse sequences.
\newblock {\em The Annals of Statistics}, 40(4):2069--2101.

\bibitem[Chae et~al., 2019]{chae2016bayesian}
Chae, M., Lin, L., and Dunson, D.~B. (2019).
\newblock Bayesian sparse linear regression with unknown symmetric error.
\newblock {\em Information and Inference: A Journal of the IMA}, 8(3):621--653.

\bibitem[De~Boor, 1978]{de1978practical}
De~Boor, C. (1978).
\newblock {\em A Practical Guide to Splines}.
\newblock New York: Springer.

\bibitem[Fikioris, 2018]{fikioris2018spectral}
Fikioris, G. (2018).
\newblock Spectral properties of {K}ac--{M}urdock--{S}zeg{\"o} matrices with a
  complex parameter.
\newblock {\em Linear Algebra and its Applications}, 553:182--210.

\bibitem[Fuller, 1987]{fuller1987measurement}
Fuller, W.~A. (1987).
\newblock {\em Measurement Error Models}.
\newblock John Wiley \& Sons.

\bibitem[Gao et~al., 2020]{gao2015general}
Gao, C., van~der Vaart, A.~W., and Zhou, H.~H. (2020).
\newblock A general framework for bayes structured linear models.
\newblock {\em Annals of Statistics}, 48(5):2848--2878.

\bibitem[Ghosal et~al., 2000]{ghosal2000convergence}
Ghosal, S., Ghosh, J.~K., and van~der Vaart, A.~W. (2000).
\newblock Convergence rates of posterior distributions.
\newblock {\em Annals of Statistics}, 28(2):500--531.

\bibitem[Ghosal and van~der Vaart, 2007]{ghosal2007convergence}
Ghosal, S. and van~der Vaart, A. (2007).
\newblock Convergence rates of posterior distributions for noniid observations.
\newblock {\em The Annals of Statistics}, 35(1):192--223.

\bibitem[Ghosal and van~der Vaart, 2017]{ghosal2017fundamentals}
Ghosal, S. and van~der Vaart, A. (2017).
\newblock {\em Fundamentals of Nonparametric {B}ayesian Inference}.
\newblock Cambridge University Press.

\bibitem[Jeong, 2020]{jeong2020posterior-b}
Jeong, S. (2020).
\newblock Posterior contraction in group sparse logit models for categorical
  responses.
\newblock {\em arXiv preprint arXiv:2010.03513}.

\bibitem[Jeong and Ghosal, 2020]{jeong2020posterior}
Jeong, S. and Ghosal, S. (2020).
\newblock Posterior contraction in sparse generalized linear models.
\newblock {\em Biometrika}, to appear.

\bibitem[Johnson and Rossell, 2012]{johnson2012bayesian}
Johnson, V.~E. and Rossell, D. (2012).
\newblock Bayesian model selection in high-dimensional settings.
\newblock {\em Journal of the American Statistical Association},
  107(498):649--660.

\bibitem[Kulkarni et~al., 1999]{kulkarni1999eigenvalues}
Kulkarni, D., Schmidt, D., and Tsui, S.-K. (1999).
\newblock Eigenvalues of tridiagonal pseudo-{T}oeplitz matrices.
\newblock {\em Linear Algebra and its Applications}, 297:63--80.

\bibitem[Magnus, 1978]{magnus1978moments}
Magnus, J.~R. (1978).
\newblock The moments of products of quadratic forms in normal variables.
\newblock {\em Statistica Neerlandica}, 32(4):201--210.

\bibitem[Martin et~al., 2017]{martin2017empirical}
Martin, R., Mess, R., and Walker, S.~G. (2017).
\newblock Empirical {B}ayes posterior concentration in sparse high-dimensional
  linear models.
\newblock {\em Bernoulli}, 23(3):1822--1847.

\bibitem[Narisetty and He, 2014]{narisetty2014bayesian}
Narisetty, N.~N. and He, X. (2014).
\newblock Bayesian variable selection with shrinking and diffusing priors.
\newblock {\em The Annals of Statistics}, 42(2):789--817.

\bibitem[Ning et~al., 2020]{ning2018bayesian}
Ning, B., Jeong, S., and Ghosal, S. (2020).
\newblock Bayesian linear regression for multivariate responses under group
  sparsity.
\newblock {\em Bernoulli}, 26(3):2353--2382.

\bibitem[Ro{\v{c}}kov{\'a}, 2018]{rockova2018bayesian}
Ro{\v{c}}kov{\'a}, V. (2018).
\newblock Bayesian estimation of sparse signals with a continuous
  spike-and-slab prior.
\newblock {\em The Annals of Statistics}, 46(1):401--437.

\bibitem[Rothman et~al., 2008]{rothman2008sparse}
Rothman, A.~J., Bickel, P.~J., Levina, E., and Zhu, J. (2008).
\newblock Sparse permutation invariant covariance estimation.
\newblock {\em Electronic Journal of Statistics}, 2:494--515.

\bibitem[Song and Liang, 2017]{song2017nearly}
Song, Q. and Liang, F. (2017).
\newblock Nearly optimal {B}ayesian shrinkage for high dimensional regression.
\newblock {\em arXiv preprint arXiv:1712.08964}.

\bibitem[van~der Vaart and Wellner, 1996]{van1996weak}
van~der Vaart, A.~W. and Wellner, J.~A. (1996).
\newblock {\em Weak Convergence and Empirical Processes}.
\newblock Springer.

\end{thebibliography}

\end{document}